\documentclass[11pt,a4paper,oneside]{book}
\usepackage{mathrsfs,amsthm,amssymb,amsmath,url}

\theoremstyle{plain}
    \newtheorem{thm}{Theorem}
    
    \newtheorem{lem}[thm]{Lemma}
    \newtheorem{prop}[thm]{Proposition}
    \newtheorem{cor}[thm]{Corollary}
    \newtheorem{fact}[thm]{Fact}
\theoremstyle{definition}
    \newtheorem{defn}[thm]{Definition}
    \newtheorem{nota}[thm]{Notation}
\theoremstyle{remark}
    \newtheorem{rem}[thm]{Remark}
\theoremstyle{remark}

\newcommand{\To}{\rightarrow}
\newcommand{\gdw}{\leftrightarrow}
\newcommand{\mult}{\times}
\newcommand{\ow}{\text{otherwise}}

\newcommand{\nin}{\notin}
\newcommand{\3}{\ss}
\newcommand{\sm}{\setminus}
\renewcommand{\c}{X \setminus}
\renewcommand{\P}{{\mathscr P}}
\newcommand{\Pow}{{\mathscr P}}

\newcommand{\nf}{_{n_f}}

\newcommand{\cl}[1]{\langle #1 \rangle}
\newcommand{\clf}[1]{\langle \{ #1 \}\rangle}

\newcommand{\OO}{{\mathscr O}}
\renewcommand{\O}{{\mathscr O}}
\newcommand{\On}{{\mathscr O}^{(n)}}
\newcommand{\Ok}{{\mathscr O}^{(k)}}
\newcommand{\Oo}{{\mathscr O}^{(1)}}
\newcommand{\Ot}{{\mathscr O}^{(2)}}

\DeclareMathOperator{\pol}{Pol}
\newcommand{\p}[1]{\pol(#1)}
\newcommand{\inv}{^{-1}}
\newcommand{\C}{{\mathscr C}}

\newcommand{\F}{{\mathscr F}}
\newcommand{\J}{{\mathscr J}}
\newcommand{\I}{{\mathscr I}}
\newcommand{\A}{{\mathscr A}}
\newcommand{\B}{{\mathscr B}}
\newcommand{\D}{{\mathscr D}}
\newcommand{\E}{{\mathscr E}}
\newcommand{\G}{{\mathscr G}}
\newcommand{\T}{{\mathscr T}}
\renewcommand{\S}{{\mathscr S}}
\newcommand{\M}{{\mathscr M}}
\newcommand{\N}{{\mathscr N}}
\newcommand{\U}{{\mathscr U}}
\renewcommand{\H}{{\mathscr H}}
\renewcommand{\L}{{\mathscr L}}

\newcommand{\un}{^{(n)}}
\newcommand{\uk}{^{(k)}}
\newcommand{\uo}{^{(1)}}
\newcommand{\ut}{^{(2)}}

\newcommand{\rest}{\upharpoonright}

\DeclareMathOperator{\Const}{Const} \DeclareMathOperator{\mi}{min}
\DeclareMathOperator{\ma}{max} \DeclareMathOperator{\med}{med}
\DeclareMathOperator{\maj}{maj}

\newcommand{\clfto}[1]{\langle \{ #1 \}\cup T_1 \rangle}
\newcommand{\cto}{\langle T_1 \rangle}
\newcommand{\pto}{\pol(T_1)}
\newcommand{\pton}{\pto^{(n)}}

\newcommand{\niceint}{[\U,\O]}
\newcommand{\lw}{\leq_W}
\newcommand{\aw}{\sim_{W}}
\newcommand{\fto}[1]{\langle #1\rangle_{T_1}}

\begin{document}
\pagenumbering{roman} \setlength{\parindent}{0cm}
\setlength{\parskip}{2ex minus 0.5ex}
\setlength{\normalbaselineskip}{\baselineskip}
\setlength{\baselineskip}{1.07 \normalbaselineskip}

\addcontentsline{toc}{chapter}{Deutsche Titelseite}
\thispagestyle{empty}

\null\vskip0.5in

\begin{center}
    \hyphenpenalty=10000\Large\uppercase\expandafter{\textbf{D i s s e r t a t i o n}}
\end{center}

\vfill

\begin{center}
    \hyphenpenalty=10000\Large\expandafter\textbf{Clones on infinite sets}
\end{center}

\vfill \vfill
\begin{center}
    ausgef\"{u}hrt zum Zwecke der Erlangung des akademischen
    Grades eines Doktors der technischen Wissenschaften unter der
    Leitung von
\end{center}
\begin{center}
    A.o. Prof. Dr. Martin Goldstern\\
    Institut f\"{u}r Diskrete Mathematik und Geometrie (104)\\
\end{center}

\vfill

\begin{center}
    eingereicht an der Technischen Universit\"{a}t Wien\\
    Fakult\"{a}t f\"{u}r Mathematik und Geoinformation
\end{center}

\vfill

\begin{center}
    von
\end{center}

\begin{center}
    Dipl. Ing. Michael Pinsker\\
    Matrikelnummer 9750160\\
    \url{marula@gmx.at}\\
\end{center}

\vfill

\vfill \vfill \vfill \vfill \vfill \vfill \vfill \vfill

 {\parbox[b]{1in}{\smash{Wien, am 3. August 2004}\vskip12pt}\hfill
\parbox[t]{3in}{\shortstack{\vrule width 3in height
0.4pt\\\small Michael Pinsker}}}
\newpage
\addcontentsline{toc}{chapter}{Vorwort mit deutscher Kurzfassung}
\chapter*{Vorwort mit deutscher Kurzfassung}

Sei $X$ eine Menge. Eine fundamentale Frage des mathematischen
Gebietes der universellen Algebra ist

\hspace{10mm}\emph{Beschreibe die Menge aller universellen
Algebren auf $X$}.

Sei nun $\A$ eine universelle Algebra auf $X$. Viele Eigenschaften
von $\A$, wie beispielsweise die Kongruenzen, die Unteralgebren,
und die Automorphismen, h\"{a}ngen nicht von den fundamentalen
Operationen von $\A$ ab, sondern von den Termoperationen, also
jenen Operationen, welche von den fundamentalen Operationen und
den Projektionen durch Funktionskomposition generiert werden. Aus
diesem Grunde bezeichnen wir zwei universelle Algebren als
\emph{\"{a}quivalent} genau dann, wenn sie dieselben
Termoperationen erzeugen. Modulo dieser \"{A}quivalenz k\"{o}nnen
wir obige Frage wie folgt formulieren:

\hspace{10mm}\emph{Beschreibe die Menge aller
\"{A}quivalenzklassen von universellen Algebren auf $X$.}

Ein \emph{Klon} ist eine Menge von Termoperationen einer
universellen Algebra auf $X$. Ebenso kann man einen Klon als Menge
endlichstelliger Funktionen auf $X$, die alle Projektionen
enth\"{a}lt und die unter Funktionskomposition abgeschlossen ist,
definieren. Die Klone entsprechen also den
Term-\"{A}quivalenzklassen von universellen Algebren auf $X$.
Ordnet man die Klone entsprechend der mengentheoretischen
Inklusion, so erh\"{a}lt man einen vollst\"{a}ndigen algebraischen
Verband $Cl(X)$. Das Ziel der Klontheorie ist die Beantwortung
eines bestimmten Aspektes obiger Frage, n\"{a}mlich

\hspace{10mm}\emph{Beschreibe $Cl(X)$.}

Diese Dissertation behandelt Teile dieser Frage, haupts\"{a}chlich
auf unendlichem $X$, und resultiert in einigen Strukturs\"{a}tzen
\"{u}ber $Cl(X)$. Die Dissertation ist in eine Einleitung und drei
Kapitel unterteilt, die unabh\"{a}ngig voneinander gelesen werden
k\"{o}nnen. Die Kapitel entsprechen den Publikationen
\cite{Pin031}, \cite{Pin032}, \cite{Pin033} des Autors.

Das Thema des ersten Kapitels sind Klone auf einer linear
geordneten Grundmenge $X$, die endlich oder unendlich sein kann.
Mithilfe der linearen Ordnung lassen sich gewisse nat\"{u}rliche
Funktionen definieren, von denen wohl die nat\"{u}rlichsten die
Maximum-, die Minimum-, und die \emph{Medianfunktionen} sind, mit
ihren offensichtlichen Definitionen. W\"{a}hrend man leicht
einsieht, da\3 eine Maximumfunktion mindestens zweier
Ver\"{a}nderlicher auch die Maximumfunktionen anderer Stelligkeit
erzeugt, und da\3 dasselbe f\"{u}r die Minimumfunktionen gilt, ist
es nicht klar, ob beispielsweise der dreistellige Median die
Medianfunktionen gr\"{o}{\3}erer Stelligkeit generiert. Unter
Verwendung kombinatiorischer Methoden zeigen wir, da\3 dies
tats\"{a}chlich der Fall, da\3 also alle Medianfunktionen
denselben Klon generieren.

Das zweite Kapitel behandelt Klone auf unendlichen Grundmengen $X$
regul\"{a}rer Kardinalit\"{a}t. Eine Funktion hei\3t fast
un\"{a}r, falls eine ihrer Variablen den Funktionswert schon bis
auf eine Menge bestimmt, deren Kardinalit\"{a}t kleiner als die
von $X$ ist. Die Menge aller fast un\"{a}ren Funktionen bildet
einen Klon, der alle (echt) un\"{a}ren Funktionen enth\"{a}lt;
dieser Klon spielt eine zentrale Rolle in der Struktur des
Klonverbandes oberhalb der un\"{a}ren Funktionen. Wir bestimmen
alle Klone, die den Klon der fast un\"{a}ren Funktionen enthalten.
Es stellt sich heraus, da\3 diese Klone unabh\"{a}ngig von der
Gr\"{o}\3e der Grundmenge eine abz\"{a}hlbar unendliche
absteigende Kette bilden, deren Durchschnitt gerade der Klon der
fast un\"{a}ren Funktionen ist.

Im dritten Kapitel wenden wir uns maximalen Klonen auf unendlichen
Mengen zu. Dabei nennen wir einen Klon \emph{maximal}, wenn er ein
Dualatom des Klonverbandes ist. Es ist bekannt, da\3 die Menge der
maximalen Klone auf unendlichem $X$ schon so gro\3 ist wie der
gesamte Klonverband; daher gibt es wenig Hoffnung, alle maximalen
Klone zu finden. Wir schr\"{a}nken die Menge der betrachteten
Klone ein und erhalten auf unendlichem $X$ regul\"{a}rer
Kardinalit\"{a}t eine explizite Liste aller maximalen Klone, die
alle Permutationen, nicht aber alle un\"{a}ren Funktionen
enthalten. Zudem bestimmen wir auf allen unendlichen Mengen $X$
alle maximalen Submonoide des Transformationsmonoids, die die
Permutationen von $X$ enthalten.
\newpage
\addcontentsline{toc}{chapter}{English title page}
\thispagestyle{empty}

\null\vskip0.5in

\begin{center}
    \hyphenpenalty=10000\Large\uppercase\expandafter{\textbf{Clones on infinite sets}}
\end{center}

\vfill

\begin{center}
    {\large\rm By\\
    {\sc Michael Pinsker}\\}
    \vskip0.3cm
    \url{marula@gmx.at}
\end{center}

\vfill

\begin{center}
    \footnotesize
    \uppercase{Dissertation} \\
    \uppercase{at the}\\
    \uppercase\expandafter{Vienna university of technology} \\
    \uppercase\expandafter{August 2004}
\end{center}

\vskip0.75in

\newpage
\addcontentsline{toc}{chapter}{Preface}
\chapter*{Preface}

Let $X$ be a set. A fundamental problem of the field of universal
algebra is

\hspace{10mm}\emph{Describe the set of all universal algebras on
$X$}.

Consider a universal algebra $\A$ on $X$. Many properties of $\A$,
such as its congruences, its subalgebras, and its automorphisms,
do not depend on the fundamental operations of $\A$, but on its
term operations, that is, the operations which are generated from
its fundamental operations and the projections by function
composition. We therefore call two universal algebras
\emph{equivalent} if and only if they have the same term
operations. Up to this equivalence, we can reformulate our problem
as follows:

\hspace{10mm}\emph{Describe the set of all term equivalence
classes of universal algebras on $X$}.

A \emph{clone} is a set of term operations of a universal algebra
on $X$. Equivalently, a clone can be defined as a set of finitary
operations on $X$ which contains the projections and which is
closed under composition. The set of all clones on $X$ thus
corresponds to the set of term equivalence classes of universal
algebras on $X$. Ordering this set by set-theoretical inclusion,
one obtains a complete algebraic lattice $Cl(X)$. The aim of clone
theory is the solution of a certain aspect of the above-mentioned
problem, namely

\hspace{10mm}\emph{Describe $Cl(X)$.}

This thesis treats instances of the latter question, mainly for
infinite $X$, resulting in several structure theorems on $Cl(X)$.
We divide this thesis into an introduction plus three chapters,
all of which can be read independently. The chapters correspond to
the author's publications \cite{Pin031}, \cite{Pin032},
\cite{Pin033}.

The first chapter deals with clones on a linearly ordered base set
$X$ (finite or infinite). Using the linear order, certain natural
functions can be defined, the most natural ones being the maximum,
the minimum, and the \emph{median functions}, with their obvious
definitions. Whereas it is easily seen that any maximum function
of at least two variables generates the maximum functions of all
arities, and that the same is true for the minimum functions, it
is not clear that the median of, say, three variables generates
the median functions of larger arities. Using combinatorial
methods, we show that this is indeed the case, that is, all median
functions generate the same clone.

In the second chapter, we turn to base sets $X$ of infinite
regular cardinality. A function is called \emph{almost unary} iff
one of its variables determines the value of the function up to a
set of cardinality smaller than the cardinality of $X$. The set of
all almost unary functions forms a clone which contains all
(really) unary functions; this clone is of importance for the
structure of the clone lattice above the unary functions. We
determine all clones containing all almost unary functions; it
turns out that independently of the size of $X$, these clones are
a countably infinite descending chain with the almost unary
functions as its intersection.

Chapter 3 is devoted to maximal clones on infinite sets. A clone
is called \emph{maximal} iff it is a dual atom in $Cl(X)$. Because
the number of maximal clones on an infinite set equals the size of
the whole clone lattice, there is little hope to find all of them.
We restrict the set of clones under consideration and provide on
all infinite $X$ of regular cardinality an explicit list of all
maximal clones which contain all permutations of $X$ but not all
unary functions. Moreover, we determine on all infinite $X$ the
maximal submonoids of the full transformation monoid which contain
the permutations.
\newpage
\addcontentsline{toc}{chapter}{Contents}
\setlength{\baselineskip}{0.7 \normalbaselineskip}
\tableofcontents\newpage \pagenumbering{arabic}
\setlength{\baselineskip}{1.07 \normalbaselineskip}

\pagestyle{myheadings}\markright{\uppercase{Introduction}\hfill}
\addcontentsline{toc}{chapter}{Introduction}
\setlength{\parindent}{0.5cm} \setlength{\parskip}{0.5ex plus
0.5ex}
\chapter*{Introduction}

Let $X$ be a set of size $|X|=\kappa$ and denote by $\On$ the set
of all $n$-ary functions on $X$. Then
$\OO=\bigcup_{n=1}^{\infty}\On$ is the set of all finitary
functions on $X$. A \emph{clone} $\C$ over $X$ is a subset of
$\OO$ which contains the projections, i.e. the functions of the
form $\pi^n_k(x_1,\ldots,x_n)=x_k$ ($1\leq k\leq n$), and which is
closed under composition. Since arbitrary intersections of clones
are obviously again clones, the set of all clones over $X$ forms a
complete lattice $Cl(X)$ with respect to inclusion. This lattice
is a subset of the power set of $\OO$. The clone lattice is
countably infinite if $X$ has only two elements, and has been
completely determined in that case by E. Post \cite{Pos41}. If $X$
is finite and has at least three elements, $Cl(X)$ is already of
size $2^{\aleph_0}$. For infinite $X$ we have
$|Cl(X)|=2^{2^{\kappa}}$. Because the clone lattice is so large in
the latter two cases, it is unlikely that it will ever be fully
described. The approach of clone theory is to investigate
interesting parts of the lattice, such as the maximal clones, the
minimal clones, or natural intervals in the lattice.

A clone is called \emph{maximal} iff it is a dual atom in $Cl(X)$.
On finite $X$ there exist finitely many maximal clones and an
explicit list of those clones has been provided by I. Rosenberg
\cite{Ros70} (see also the diploma thesis \cite{Pin02} for a
self-contained proof of Rosenberg's Theorem). Moreover, the clone
lattice is dually atomic in that case, that is, every clone is
contained in a maximal one. If $X$ is infinite, then the number of
maximal clones equals the size of the whole clone lattice
(\cite{Ros76}, see also \cite{GS022}), so that it seems impossible
to determine all of them. It has also been shown \cite{GS03} that
if the continuum hypothesis holds, then the clone lattice on a
countably infinite base set is not dually atomic. We will deal
with maximal clones in Chapters 2 and 3: In the second chapter, we
obtain on all $X$ of infinite regular cardinality a simple
description of a certain maximal clone above $\Oo$ which is
important for the structure of the interval $[\Oo,\O]$ of the
clone lattice. In the third chapter, we give an explicit list of
the maximal clones which contain the set $\S$ of all permutations
on $X$ but which do not contain $\Oo$.

A \emph{minimal} clone on $X$ is an atom in the lattice $Cl(X)$,
i.e. a minimal element in $Cl(X)\setminus \{\J\}$, where $\J$ is
the trivial clone containing only the projections. Clearly every
minimal clone is generated by a single nontrivial function.
Functions which generate minimal clones are called \emph{minimal}
as well. On finite $X$, the minimal clones are finite in number
and every clone contains a minimal one. Surprisingly, there is no
characterization of minimal clones even on finite $X$. If we take
the base set $X$ to be infinite, then the number of minimal clones
is $2^{\kappa}$, and it is easy to see that not every clone
contains a minimal one. The first chapter deals with a certain
minimal clone on a linearly ordered base set $X$, namely the clone
generated by the median functions.

Because the clone lattice is too large to completely understand
it, it makes sense to pick feasible \emph{intervals} of it and try
to determine them. For example, there exist a number of results on
the interval $[\Oo,\O]$ of clones containing all unary functions.
One such result due to G. Gavrilov \cite{Gav65} is that on
countably infinite $X$, there exist only two maximal clones in
this interval. M. Goldstern and S. Shelah \cite{GS03} proved that
the same is true on $X$ of weakly compact cardinality, but showed
in the same article that on most other cardinals, in particular on
all successors of regulars, there exist $2^{2^\kappa}$ such
clones. We will prove another structure theorem for clones above
$\Oo$ in Chapter 2, determining the interval $[\U,\O]$ of clones
containing all almost unary functions.\\
Another example of an interesting interval in the interval
$[\S,\O]$ of clones containing all permutations of $X$. L.
Heindorf \cite{Hei02} determined on countably infinite $X$ all
maximal clones in this interval. We will extend his result to all
infinite $X$ of regular cardinality in Chapter 3, obtaining an
explicit list of all maximal clones which contain the permutations
but not $\Oo$.\\
The interval $[\J,\Oo]$ consists of those clones which contain
only \emph{essentially unary} functions, i.e. functions that
depend only on one of their variables. Such clones are essentially
submonoids of the full transformation monoid $\Oo$. It is known
that the number of dual atoms in this interval is
$2^{2^{\kappa}}$, so there is no hope to determine them. However,
G. Gavrilov \cite{Gav65} found all dual atoms of this interval
which contain $\S$ (so he found the dual atoms of $[\S,\Oo]$), on
countably infinite $X$. We will generalize his theorem to all
infinite $X$ in Chapter 3.

For extensive introductions to clone theory (although primarily on
finite base sets), we refer to the monograph \cite{Sze86} by
\'{A}. Szendrei and the textbook \cite{PK79} by R. P\"{o}schel and
L. Kalu\v{z}nin.\newpage
\pagestyle{headings}
\pagestyle{myheadings}\markright{\uppercase{Median
functions}\hfill}

\chapter[The clone generated by the median
functions]{The clone generated by the median functions}

Let $X$ be a linearly ordered set of arbitrary size (finite or
infinite). Natural functions on such a set one can define using
the linear order include maximum, minimum and median functions.
While it is clear what the clone generated by the maximum or the
minimum looks like, this is not obvious for the median functions.
We show that every clone on $X$ contains either no median function
or all median functions, that is, the median functions generate
each other.

\section{The median functions}
    Assume $X$ to be linearly ordered. We
    emphasize that the cardinality of $X$ is not relevant.
    For all $n\geq 1$ and all $1\leq k\leq n$ we define a function
    $$
        m^n_k(x_1,\ldots,x_n)=x_{j_k}\quad \text{if} \,\, x_{j_1} \leq \ldots \leq
        x_{j_n}.
    $$
    In words, the function $m^n_k$ returns the $k$-th smallest element from an
    $n$-tuple. The functions $m^n_k$ are totally symmetric, i.e., invariant unter all permutations of their variables, and
    $m^n_k(x_1,\ldots,x_n)=x_k$ whenever $x_1\leq \ldots \leq x_n$. For example, $m^n_n$ is the maximum function $\max_n$
    and $m^n_1$ the minimum function $\min_n$ in $n$ variables. If
    $n$ is an odd number then we call $m^n_{\frac{n+1}{2}}$ the
    $n$-th median function and denote this function by $\med_n$.

    It is easy to check what the clones
    generated by the functions $\max$ and $\min$ look like:
    $$
        \clf{\ma_n}=\{\ma_k(\pi^j_{i_1},\ldots,\pi^j_{i_k}) : 1\leq i_1,\ldots,i_k\leq j,\,1\leq k\leq j \}
    $$
    and
    $$
        \clf{\mi_n}=\{\mi_k(\pi^j_{i_1},\ldots,\pi^j_{i_k}) : 1\leq i_1,\ldots,i_k\leq j ,\,1\leq k\leq j
        \},
    $$
    where $n\geq 2$ is arbitrary. In particular, the two clones
    are minimal. Now it is natural to ask which of these properties hold for the
    functions ``in between'', that is the $m^n_k$ as defined before, most importantly the median
    functions.
    We will show that for odd $n\geq 3$
    $$
        \clf{\med_n}\supseteq \{\med_k(\pi^j_{i_1},\ldots,\pi^j_{i_k}) : 1\leq i_1,\ldots,i_k\leq j,\,1\leq k\leq j,\, k \text{ odd}
        \},
    $$
    but one readily constructs functions in that clone which are
    not a median function and not a projection. However, R. P\"{o}schel
    and L. Kalu\v{z}nin observed in \cite{PK79},
    Theorem 4.4.5, that the median of three variables (and hence by our result, all medians)
    does generate a minimal clone.
    \begin{thm}\label{medIsMinimal}
        The clone generated by the function $\med_3$ is minimal.
    \end{thm}

    We are going to prove

    \begin{thm}\label{medians}
        Let $k,n \geq 3$ be odd natural numbers. Then $\med_k\in \langle \{\med_n\} \rangle$. In
        other words, a clone contains either no median function or
        all median functions.
    \end{thm}
\subsection{Notation}

    For a set of functions $\F$ we shall denote the smallest
    clone containing $\F$ by $\langle \F \rangle$. If $1\leq k\leq n$, we write $\pi^n_k$ for the $n$-ary projection
    on the $k$-th component.\\
    For a positive rational number $q$ we write
    $$
        \lfloor q \rfloor = \max\{n\in\mathbb{N}:\,n\leq q\}
    $$
    and
    $$
        \lceil q \rceil = \min\{n\in\mathbb{N}:\,q\leq n\}.
    $$
    If $a\in X^n$ is an $n$-tuple and $1\leq k\leq n$ we write $a_k$ for the $k$-th component of $a$.
    We will assume $X$ to be linearly ordered by the relation $\leq$ and let $<$ carry the obvious meaning.

\section{The proof of Theorem \ref{medians}}
\subsection{Almost divisibility}

    We split the proof of the theorem into a sequence of
    lemmas.

    \begin{defn}
        Let $k,n\geq 1$ be natural numbers. Denote by $R(\frac{n}{k})$ the remainder of the
        division $\frac{n}{k}$. We say that $n$ is \emph{almost divisible by}
        $k$ iff either $R(\frac{n}{k})\leq \frac{n}{k}$ or $(k-1)-R(\frac{n}{k})\leq
        \frac{n}{k}$.
    \end{defn}

    Note that $n$ is almost divisible by $k$ if it is divisible by
    $k$. The following lemma tells us which medians of smaller
    arity are generated by $\med_n$ by simple identification of
    variables (see also Remark \ref{REM:almostDivisibility}).

    \begin{lem}\label{almdivmed}
        Let $k\leq n$ be odd natural numbers. If $n$ is almost
        divisible by $k$, then $\med_k\in \langle \{\med_n\}\rangle$.
    \end{lem}

    \begin{proof}
        We claim that
        $$
            \med_k(x_1,\ldots,x_k)=\med_n(x_1,\ldots,x_1,x_2,\ldots,x_2,\ldots,x_k,\ldots,x_k),
        $$
        where $x_j$ occurs in the $n$-tuple $\lfloor\frac{n}{k}\rfloor +1$ times if
        $j\leq R(\frac{n}{k})$ and $\lfloor\frac{n}{k}\rfloor$
        times otherwise. Assume $\med_k(x_1,\ldots,x_k)=x_j$. Then
        there are at most $\frac{k-1}{2}$ components smaller than
        $x_j$ and at most $\frac{k-1}{2}$ components larger than
        $x_j$. Thus in our $n$-tuple, there are at most
        \begin{eqnarray}{\label{miss}}
        \frac{k-1}{2}\lfloor\frac{n}{k}\rfloor+\min(R(\frac{n}{k}),\frac{k-1}{2})
        \end{eqnarray}
        elements smaller (larger) than $x_j$.\\
        \textit{Case 1.}  $R(\frac{n}{k})\leq \frac{k-1}{2}$.\\ Since
        $n$ is almost divisible by $k$, we have either $R(\frac{n}{k})\leq
        \frac{n}{k}$ or $(k-1)-R(\frac{n}{k})\leq \frac{n}{k}$. In
        the latter case,
        $$
            R(\frac{n}{k})\leq \frac{k-1}{2} \quad \wedge \quad (k-1)-R(\frac{n}{k})\leq \frac{n}{k}
        $$
        and so
        $$
            R(\frac{n}{k})\leq \frac{n}{k}.
        $$
        Thus in either of the cases, we can calculate from
        (\ref{miss})
        \begin{eqnarray*}
        \begin{aligned}
            &\frac{k-1}{2}\lfloor\frac{n}{k}\rfloor+R(\frac{n}{k})\\
            &=\frac{1}{2}(k \lfloor\frac{n}{k}\rfloor +
            R(\frac{n}{k}))+\frac{1}{2}(R(\frac{n}{k})-\lfloor\frac{n}{k}\rfloor)\\
            &=\frac{n}{2}+\frac{1}{2}(R(\frac{n}{k})-\lfloor\frac{n}{k}\rfloor)\\
            &\leq\frac{n}{2}\\
        \end{aligned}
        \end{eqnarray*}
        and so $\med_n$ yields $x_j$.\\
        \textit{Case 2.} $\frac{k-1}{2}< R(\frac{n}{k})$.\\
        Again we know that either $R(\frac{n}{k})\leq
        \frac{n}{k}$ or $(k-1)-R(\frac{n}{k})\leq \frac{n}{k}$. In
        the first case, we see that
        $$
             \frac{k-1}{2}< R(\frac{n}{k}) \quad \wedge \quad R(\frac{n}{k})\leq \frac{n}{k}
        $$
        implies
        $$
            (k-1)-R(\frac{n}{k})\leq \frac{n}{k}
        $$
        and so (\ref{miss}) yields at most
        \begin{eqnarray*}
        \begin{aligned}
            &\frac{k-1}{2}\lfloor\frac{n}{k}\rfloor+\frac{k-1}{2}\\
            &=\frac{k-1}{2}\lfloor\frac{n}{k}\rfloor+\frac{1}{2}R(\frac{n}{k})+\frac{k-1}{2}-\frac{1}{2}R(\frac{n}{k})\\
            &=\frac{1}{2}(k \lfloor\frac{n}{k}\rfloor +
            R(\frac{n}{k}))-\frac{1}{2}\lfloor\frac{n}{k}\rfloor+\frac{k-1}{2}-\frac{1}{2}R(\frac{n}{k})\\
            &\leq\frac{n}{2}+\frac{1}{2}(-\lfloor\frac{n}{k}\rfloor+(k-1)-R(\frac{n}{k}))\\
            &\leq\frac{n}{2}\\
        \end{aligned}
        \end{eqnarray*}
        components which are smaller (larger) than $x_j$. This finishes the
        proof.
    \end{proof}
    \begin{cor}
        Let $k,n\geq 1$ be odd natural numbers. If $k\leq
        \sqrt{n}$, then $\med_k$ is generated by $\med_n$.
    \end{cor}
    \begin{proof}
        Trivially, $R(\frac{n}{k})\leq k-1$ and $k-1 \leq \frac{n}{k}$ as $k\leq
        \sqrt{n}$. Hence, $n$ is almost divisible by $k$.
    \end{proof}
    \begin{cor}\label{medn2med3}
        Let $n\geq 3$ be odd. Then $\med_3\in\clf{\med_n}$.
    \end{cor}
    \begin{proof}
        Simply observe that all $n\geq 4$ are almost divisible by $3$.
    \end{proof}
\subsection{Majority functions}
    We have seen that we can get small (that is, of small arity) median functions out of
    large ones. The converse inclusion is shown with the help of majority
    functions.
    \begin{defn}
        Let $f\in\On$. We say that $f$ is a \emph{majority
        function} iff $f(x_1,\ldots,x_n)=x$ whenever the value $x$
        occurs at least $\lceil\frac{n+1}{2}\rceil$ times among
        $(x_1,\ldots,x_n)$.
    \end{defn}

    Note that $\med_n$ is a majority function for all odd $n$.
    We observe now that we can build a ternary majority function
    from most larger ones by identifying variables.

    \begin{lem}\label{LEM:ternaryMaj}
        Let $n\geq 5$ and let $\maj_n\in\On$ be a majority
        function. Then $\maj_n$ generates a majority function of
        three arguments.
    \end{lem}
    \begin{proof}
        Set
        $$
            \maj_3=\maj_n(x_1,\ldots,x_1,x_2,\ldots,x_2,x_3,\ldots,x_3),
        $$
        where $x_j$ occurs in the $n$-tuple $\lfloor\frac{n}{3}\rfloor +1$ times if
        $j\leq R(\frac{n}{3})$ and $\lfloor\frac{n}{3}\rfloor$
        times otherwise. It is readily verified that $\maj_3$ is a
        majority function.
    \end{proof}

    The following lemma tells us that we can generate majority
    functions of even arity from majority functions of odd arity.

    \begin{lem}\label{lem:evenmaj}
        Let $n\geq 2$ be an even natural number. Then we can get
        an $n$-ary majority function $\maj_n$ out of any $(n+1)$-ary majority
        function $\maj_{n+1}$.
    \end{lem}
    \begin{proof}
        Set
        $$
            \maj_n(x_1,\ldots,x_n)=\maj_{n+1}(x_1,\ldots,x_n,x_n)
        $$
        and let $x\in X$ have a majority among $(x_1,\ldots,x_n)$.
        Since $n$ is even, $x$ occurs $\frac{n}{2}+1$ times in the
        $n$-tuple which is enough for a majority in the $(n+1)$-tuple
        $(x_1,\ldots,x_n,x_n)$.
    \end{proof}
    We now show that we can construct large majority functions out of
    small ones. This has already been known but we include our own proof
    here.
    \begin{lem} \label{maj}
        Let $n\geq 5$ be a natural number. Then we can
        construct an $n$-ary majority function out of any $(n-2)$-ary
        majority function $\maj_{n-2}$.
    \end{lem}
    \begin{proof}
        For $2\leq j\leq n-1$ and $1\leq i\leq n-1$ with $i\neq j$
        we define functions
        $$
            \gamma_i^j=
            \begin{cases}
                \maj_{n-2}(x_1,\ldots,x_{i-1},x_{i+2},\ldots,x_n)&j\neq
                i+1\\
                \maj_{n-2}(x_1,\ldots,x_{i-1},x_{i+1},x_{i+3},\ldots,x_n)&j=
                i+1\\
            \end{cases}
        $$
        In words, given an $n$-tuple $(x_1,\ldots,x_n)$, $\gamma_i^j$ ignores
        $x_i$ and the next component of the $n$-tuple which is not
        $x_j$ and calculates $\maj_{n-2}$ from what is left.
        Set
        $$
            z_j=\maj_{n-2}(\gamma_1^j,\ldots,\gamma_{j-1}^j,\gamma_{j+1}^j,\ldots,\gamma_{n-1}^j)
        $$
        and
        $$
            f=\maj_{n-2}(z_2,\ldots,z_{n-1}).
        $$
        The function $f$ is an $n$-ary term of depth three over
        $\{\maj_{n-2}\}$.\\
        \textit{Claim.} $f$ is a majority function.\\
        We prove our claim for the case where $n$ is odd. The same proof works in the
        even case, the only difference being that the counting is slightly different (a majority occurs
        $\frac{n+2}{2}$ times instead of $\frac{n+1}{2}$, and so on). We leave the verification of this to the diligent
        reader.

        Assume $x\in X$ has a majority. If $x$ occurs more than
        $\frac{n+1}{2}$ times, then it is readily verified that all the $\gamma_i^j$ yield
        $x$ and so do all $z_j$ and so does $f$. So say $x$ appears
        exactly $\frac{n+1}{2}$ times among the variables of $f$.

        Next we observe that if $x_j=x$, then $z_j=x$: For if
        $\gamma_i^j\neq x$, then both components ignored in $\gamma_i^j$, that is, $x_i$ and the component
        after $x_i$ which is not $x_j$, have to be equal to $x$. We can count
        $$
            |\{i:\gamma_i^j\neq x\}|\leq|\{i\neq j:x_i=x\}\setminus\{\max(i\neq j:
            x_i=x)\}|\leq\frac{n-1}{2}-1=\frac{n-3}{2}.
        $$
        Thus, $z_j=x$.

        Now we shall count a second time to see
        that if $x_1\neq x$ or $x_n\neq x$, then $f=x$: Say
        without loss of generality
        $x_1\neq x$. Then
        $$
            |\{2\leq j\leq n-1:
            x_j=x\}|\geq\frac{n+1}{2}-1=\frac{n-1}{2}
        $$
        and since we have seen that $z_j=x$ for all such $j$ we
        indeed obtain $f=x$.

        In a last step we consider the case where both $x_1=x$ and
        $x_n=x$. Let
        $$
            k=\min\{i: x_i\neq x\}
        $$
        and
        $$
            l=\max\{i:x_i\neq x\}.
        $$
        Since $n\geq 5$ those two indices are not equal. Count
        $$
            |\{i:\gamma^l_i\neq
            x\}|\leq|\{i:x_i=x\}\setminus\{k-1,n\}|=\frac{n+1}{2}-2=\frac{n-3}{2}.
        $$
        Thus, $z_l=x$ and we count for the last time
        $$
            |\{j: z_j=x\}|\geq|\{2\leq j\leq
            n-1:x_j=x\}\cup\{l\}|=\frac{n-3}{2}+1=\frac{n-1}{2},
        $$
        so that also in this case $f=x$.
    \end{proof}
    We conclude that if a clone contains a majority function, then
    it contains majority functions of all arities.
    \begin{cor}\label{allmaj}
        Let $n,k\geq 3$ be natural numbers. Assume $\maj_n\in\On$
        is any majority function. Then $\maj_n$ generates a majority function in $\Ok$.
    \end{cor}
    \begin{proof}
        If $k\geq n$ and $n,k$ are either both even or both odd, then we can iterate Lemma
        \ref{maj} to generate a majority function of arity $k$. Lemma \ref{lem:evenmaj}
        takes care of the case when $k$ is even but $n$ is odd.

        In all other cases with $n\geq 5$, generate a ternary
        majority function from $\maj_n$ first with the help of Lemma
        \ref{LEM:ternaryMaj} and follow the procedure just
        described for the other case.

        Finally, if $n=4$, we can build a majority function $\maj_6$
        from $\maj_4$ first and are back in one of the other
        cases.
    \end{proof}
    Now we use the large majority functions to obtain large median functions.
    \begin{lem}
        For all odd $n\geq 3$ there exists $b\geq n$ such that
        $\med_n\in\clf{\med_3,\maj_b}$ for an arbitrary $b$-ary majority
        function $\maj_b$.
    \end{lem}
    \begin{proof}
        Let $n$ be given. Our strategy to calculate the median from an $n$-tuple will be the following: We apply
        $\med_3$ to all possible selections of three elements of the $n$-tuple. The results we write to an $n_1$-tuple,
        from which we again take all possible selections of three
        elements. We apply $\med_3$ again to these selections and
        so forth. Now the true median of the original $n$-tuple
        ``wins'' much more often in this procedure than the other
        elements, so that after a finite number of steps (a number we
        can give a bound for) more than half of the components of the then giant tuple
        have the true median as their value. To that tuple we
        apply a majority function and obtain the median.

        In detail, we define two sequences
        $(n_j)_{j\in\omega}$ and $(k_j)_{j\in\omega}$ by
        $$
            n_0=n,\quad n_{j+1}=\binom{n_j}{3}
        $$
        and
        $$
            k_0=1,\quad
            k_{j+1}=\binom{k_j}{3}+\binom{k_j}{2}(n_j-k_j)+\binom{k_j}{1}(\frac{n_j-k_j}{2})^2.
        $$
        The sequences have the following meaning: Given an
        $n_j$-tuple, there are $n_{j+1}$ possible selections of
        three elements of the tuple to which we apply the median $\med_3$. If the median of the $n_j$-tuple
        (which is equal to the median of the $n_0=n$-tuple) appeared at least
        $k_j$ times there, then it appears at least $k_{j+1}$
        times in the resulting $n_{j+1}$-tuple. Read $k_{j+1}$ as follows: We assume the worst case, namely that the median
        occurs only once in the original $n$-tuple, so $k_0=1$. If we pick
        three elements from the $n_j$-tuple and calculate $\med_3$, then the result is the median we are looking for
        if either all three elements are equal to the median ($\binom{k_j}{3}$ possibilities) or two are equal to the median
        ($\binom{k_j}{2}(n_j-k_j)$ possibilities) or one is equal to the median,
        one is smaller, and one is larger ($\binom{k_j}{1}\,(\frac{n_j-k_j}{2})^2$
        possibilities). Set
        $r_j=\frac{k_j}{n_j}$ for $j\geq 0$ to be the relative frequency of the median in the tuple after $j$ steps.
        We claim that $\limsup(r_j)_{j\in\omega}=1$:
        \begin{eqnarray*}
        \begin{aligned}
            r_{j+1}&=\frac{k_{j+1}}{n_{j+1}}\\
                   &=\frac{k_j}{n_j}\,\,\frac{(k_j-1)(k_j-2)+3(k_j-1)(n_j-k_j)+\frac{3}{2}(n_j^2-2n_jk_j+k_j^2)}{(n_j-1)(n_j-2)}\\
                   &=r_j\,\,\frac{3(n_j-1)^2+1-k_j^2}{2(n_j-1)(n_j-2)}
        \end{aligned}
        \end{eqnarray*}
        Further calculation yields
        \begin{eqnarray*}
        \begin{aligned}
            r_{j+1}&\geq r_j\,\,\frac{3(n_j-1)^2+1-k_j^2}{2(n_j-1)^2}\\
            &= r_j\,\,
            (\frac{3}{2}-\frac{(k_j-1)^2}{2(n_j-1)^2}-\frac{k_j-1}{(n_j-1)^2})\\
            &\geq r_j\,\,
            (\frac{3}{2}-\frac{1}{2}r_{j}^2-\frac{r_{j}}{n_j-1})\\
            &\geq r_j\,\,
            (\frac{3}{2}-\frac{1}{2}r_{j}^2-\frac{1}{n_j-1}).
        \end{aligned}
        \end{eqnarray*}
        Suppose towards a contradiction that
        $(r_j)_{j\in\omega}$ is bounded away from $1$ by $p\,$: $r_j<p<1$ for all $j\in \omega$.
        Choose $j$ large enough so that
        $$\frac{1}{n_j-1}<\frac{1-p}{4}.$$
        Then
        \begin{eqnarray*}
        \begin{aligned}
            r_{i+1}&>r_i\,\,(\frac{3}{2}-\frac{p}{2}-\frac{1-p}{4})\\
            &=r_i\,\,(1+\frac{1-p}{4})
        \end{aligned}
        \end{eqnarray*}
        for all $i\geq j$ so that there exists $l>j$
        such that $r_l>p$, in contradiction to our assumption. Hence,
        $\limsup(r_j)_{j\in\omega}=1$.

        Now if we calculate $j$ such that $r_j>\frac{1}{2}$, and
        choose $b=n_j$,
        we can obtain the median with the help of a $b$-ary majority function.
    \end{proof}
    We are ready to prove our main theorem.
    \begin{proof}[Proof of Theorem \ref{medians}]
        Let $k, n$ be given. Corollary \ref{medn2med3} tells us that we can construct $\med_3$ out of $\med_n$.
        Since $\med_3$ is also a majority function, we can get majority
        functions of arbitrary arity with the help of Corollary \ref{allmaj}. Then by the preceding lemma, we can
        generate $\med_k$.
    \end{proof}
    \begin{rem}\label{REM:almostDivisibility}
        In fact, the lemma on almost divisibility is not needed
        for the proof of the theorem, since we only have to get
        $\med_3$ out of $\med_n$ (and $\med_3(x_1,x_2,x_3)=\med_n(x_1,x_2,\ldots,x_2,x_3,\ldots,x_3)$
        where $x_2$ and $x_3$ occur $\frac{n-1}{2}$ times
        in the $n$-tuple) and then apply Lemma \ref{maj} to generate large
        majority functions. Still, the lemma shows what we can
        construct by simple identification of variables.
    \end{rem}
\section{Minimality of the $m^n_k$}
    We mentioned that the clones generated by the maximum, the minimum and the median functions are
    minimal. Anyone who hoped that the same holds for all $m^n_k$ will be
    disappointed by the following lemma.
    \begin{lem}
        Let $n\geq 4$ and $2\leq k \leq \lfloor\frac{n}{2}\rfloor$.
        Then $m^n_k$ is not a minimal function.
    \end{lem}
    \begin{proof}
        It is enough to see that
        $$
            \mi_2(x,y)=m^n_k(x,\ldots,x,y,\ldots,y)\in\clf{m^n_k},
        $$
        where $x$ occurs in the $n$-tuple exactly
        $\lfloor\frac{n}{2}\rfloor$ times. The clone generated
        by $\min_2$ is obviously a nontrivial proper subclone of $\clf{m^n_k}$.
    \end{proof}
    Now comes just another disappointment.
    \begin{lem}
        Let $\lceil\frac{n}{2}\rceil< k < n$. Then $m^n_k$ is not
        a minimal function.
    \end{lem}
    \begin{proof}
        This time we have that
        $$
            \ma_2(x,y)=m^n_k(x,\ldots,x,y,\ldots,y)\in\clf{m^n_k},
        $$
        where $x$ occurs in the $n$-tuple exactly
        $\lfloor\frac{n}{2}\rfloor$ times. The clone generated by
        the maximum functions is obviously a proper subclone of
        $\clf{m^n_k}$.
    \end{proof}
    We summarize our results in the following corollary.
    \begin{cor}
        Let $n\geq 2$ and $1\leq k\leq n$. Then $m^n_k$ is minimal
        iff $k=1$ or $k=n$ or $n$ is odd and $k=\frac{n+1}{2}$.
        That is, the minimal functions among the $m^n_k$ are
        exactly the maximum, the minimum and the median functions.
    \end{cor}
\section{Variations of the median function}
    For even natural numbers $n$ we did not define median
    functions. One could consider the so-called
    ``lower
    median'' instead:
    $$
        \med^{low}_n=m^n_{\frac{n}{2}}
    $$
    But as a consequence of the preceding corollary, $\med^{low}_n$ is
    not generated by the real medians and does therefore not serve
    as a perfect substitute. For the same reason, the ``upper
    median''
    $$
        \med^{upp}_n=m^n_{\frac{n}{2}+1}
    $$
    is not an ideal replacement either.

    However, the other direction almost works: $\med^{low}_n$
    generates the medians if and only if $n\geq 6$. Indeed,
    simple identification of variables suffices:
    $$
        \med_3=\med^{low}_n(x_1,\ldots,x_1,x_2,\ldots,x_2,x_3,\ldots,x_3)
    $$
    where $x_j$ occurs in the $n$-tuple $\lfloor\frac{n}{3}\rfloor+1$ times
    if $j\leq R(\frac{n}{3})$ and $\lfloor\frac{n}{3}\rfloor$
    times otherwise. Of course we can do the same with the upper medians. It is easy to see that $\med^{low}_4$
    cannot generate the medians.

    One could have the idea of using a more general notion of
    median functions: Let $(X,\wedge,\vee)$ be a lattice. Define
    $$
        {\widetilde{m}}_k^n(x_1,\ldots,x_n)=\bigwedge_{(j_1,\ldots,j_k)\in \{1,\ldots,n\}^k}\bigvee_{1\leq
        i\leq k}x_{j_i}.
    $$
    If the order induced by the lattice on $X$ is a chain,
    this definition agrees with our definition of $m^n_k$.
    However, although we can get $\widetilde{\med}_3$ out of $\widetilde{\med}_n$
    just like described in Remark \ref{REM:almostDivisibility}, our proof to
    obtain large medians via majority functions fails. We do
    not know under which conditions on the lattice the same results can be
    obtained.
\newpage
\pagestyle{myheadings}\markright{\uppercase{Almost unary
functions}\hfill}
\begin{chapter}{Clones containing all almost unary functions}

Let $X$ be an infinite set of regular cardinality. We determine
all clones on $X$ which contain all almost unary functions. It
turns out that independently of the size of $X$, these clones form
a countably infinite descending chain. Moreover, all such clones
are finitely generated over the unary functions. In particular, we
obtain an explicit description of the only maximal clone in this
part of the clone lattice. This is especially interesting if $X$
is countably infinite, in which case it is known that such a
simple description cannot be obtained for the second maximal clone
over the unary functions.

\begin{section}{Background}

\subsection{Almost unary functions}

    Let $X$ be of infinite regular cardinality from now on unless otherwise stated. We
    call a subset $S\subseteq X$ \emph{large} iff $|S|=|X|$, and
    \emph{small} otherwise. If $X$ is itself a regular cardinal,
    then the small subsets are exactly the bounded subsets of $X$.
    A function $f(x_1,\ldots,x_n)\in\On$ is \emph{almost unary} iff
    there exists a function $F:X\To \Pow(X)$ and $1\leq k\leq n$ such that $F(x)$ is small for all
    $x\in X$ and such that for all
    $(x_1,\ldots,x_n)\in X^n$ we have $f(x_1,\ldots,x_n)\in F(x_k)$. If
    we assume $X$ to be a regular cardinal itself, this is
    equivalent to the existence of a function $F\in\Oo$ and a
    $1\leq k\leq n$ such that $f(x_1,\ldots,x_n)< F(x_k)$ for all $(x_1,\ldots,x_n)\in
    X^n$. Because this is much more convenient and does not
    influence the properties of the clone lattice, we shall assume
    $X$ to be a regular cardinal throughout this chapter. Let $\U$ be the set of all almost unary functions. It is
    readily verified that $\U$ is a clone.
    We will determine all clones which contain $\U$; in particular, such clones contain $\Oo$.

\subsection{Maximal clones above $\Oo$}

    Although on an infinite set $X$ not every clone must be contained in a maximal one \cite{GS03},
     the sublattice of $Cl(X)$ of functions containing $\Oo$ is dually atomic by Zorn's lemma,
     since $\OO$ is finitely generated over $\Oo$.
    G. Gavrilov proved in \cite{Gav65} that for countably
    infinite $X$ there are only two maximal clones containing all unary
    functions. M. Goldstern and S. Shelah extended this result to
    clones on weakly compact cardinals in the article \cite{GS022}, where an uncountable
    cardinal $X$ is called \emph{weakly compact} iff whenever we
    colour the edges of a complete graph $G$ of size $X$ with two
    colours, then there exists a complete subgraph of $G$ of size
    $X$ on which the colouring is constant. In the same paper, the authors proved
    that on other regular cardinals $X$ satisfying a certain partition
    relation there are even $2^{2^X}$ maximal clones above $\Oo$.

    There exists exactly one maximal clone above $\U$. So far,
    this clone has been defined using the following concept:
    Let $\rho\subseteq X^J$ be a relation on $X$ indexed by $J$ and let
    $f\in\On$. We say that $f$ \emph{preserves} $\rho$ iff for all $r^1=(r^1_i:i\in J),\ldots,r^n=(r^n_i:i\in J)$ in
    $\rho$ we have $(f(r^1_i,\ldots,r^n_{i}):i\in J) \in \rho$.
    We define the set of \emph{polymorphisms} $\pol(\rho)$ of $\rho$ to be the set of all functions in
    $\O$ preserving $\rho$; $\pol(\rho)$ is easily seen to be a clone. In particular, if
    $\rho\subseteq X^{X^k}$ is a set of $k$-ary functions, then a
    function $f\in\On$ preserves $\rho$ iff for all functions
    $g_1,\ldots,g_n$ in $\rho$ the composite $f(g_1,\ldots,g_n)$ is a
    function in $\rho$.

    Write
    $$
        T_1=\U\ut=\{f \in \Ot: f \,\,\text{almost unary}\}.
    $$

    The following was observed by G. Gavrilov \cite{Gav65} for
    countable base sets and extended to all regular $X$ by R. Davies and I.
    Rosenberg \cite{DR85}. Uniqueness on uncountable regular cardinals
    is due to M. Goldstern and S. Shelah \cite{GS022}.
    \begin{thm}\label{THM:poltone}
        Let $X$ have infinite regular cardinality. Then $\pol(T_1)$ is a maximal clone containing all unary
        functions. Furthermore, $\pto$ is the only maximal clone
        containing all almost unary functions.
    \end{thm}

    For $S$ a subset of $X$ we set
    $$
        \Delta_S=\{(x,y)\in S^2 : y < x\},\quad \nabla_S=\{(x,y)\in S^2 : x
        <y\}.
    $$
    We will also write $\Delta$ and $\nabla$ instead of $\Delta_X$ and $\nabla_X$. Now define
    $$
        T_2=\{f\in\Ot : \forall S\subseteq X \,(S \,\,\text{large}
        \To\,
        \text{neither} \, f\upharpoonright_{\Delta_S}\,\text{nor} \,
        f\upharpoonright_{\nabla_S}\, \text{are 1-1})\}.
    $$
    The next result is due to G. Gavrilov \cite{Gav65} for $X$ a countable set and due to M. Goldstern and S. Shelah
    \cite{GS022} for $X$ weakly compact.
    \begin{thm}\label{THM:theTwoMaximal}
         Let $X$ be countably infinite or weakly compact.
         Then $\pol(T_2)$ is a maximal clone which contains $\Oo$. Moreover,
         $\pol(T_1)$ and $\pol(T_2)$ are the only maximal clones above $\Oo$.
    \end{thm}

    The definition of $\pol(T_2)$ not only looks more
    complicated than the one of $\pol(T_1)$. First of all, a result of R. Davies and I.
    Rosenberg in \cite{DR85} shows that assuming the continuum hypothesis, $T_2$ is not closed under
    composition on $X=\aleph_1$ and so it is unclear what $\pol(T_2)$ is.
    Secondly, on countable $X$, if we equip $\OO$ with the natural topology which we shall specify later,
    then $T_2$ is a non-analytic set in that space and so is
    $\pol(T_2)$; in particular, neither $\cl{T_2}$ nor $\pol(T_2)$ are countably generated over
    $\Oo$ (see \cite{Gol02}), where for a set of functions $\F$ we denote by $\cl{\F}$ the clone generated by $\F$.
    The clones $\cl{T_1}$ and $\pto$ on the other hand turn out to be rather simple with respect to this topology,
    and both clones are finitely generated $\Oo$.

    Fix any
    injection $p$ from $X^2$ to $X$; for technical reasons we assume that $0\in X$ is not in the range of
    $p$.
    \begin{fact}
        $\cl{\{p\}\cup\Oo}=\O$, that is, the function $p$ together with $\Oo$ generate $\O$.
    \end{fact}
    For a subset
    $S$ of $X^2$ we write
    $$
        p_S(x_1,x_2)=\begin{cases}p(x_1,x_2)&,(x_1,x_2)\in
        S\\0&,\ow\\ \end{cases}
    $$
    M. Goldstern observed the following \cite{Gol02}. Since the
    result has not yet been published, but is important for our investigations, we include a proof here.
    \begin{fact}\label{FAC:pDelta}
        $\cl{\{p_\Delta\}\cup\Oo}=\cto$.
    \end{fact}
    \begin{proof}
        Set $\C= \cl{\{p_\Delta\}\cup\Oo}$. Since $p_\Delta(x_1,x_2)$ is
        obviously bounded by the unary function
        $\gamma(x_1)=\sup\{p_\Delta(x_1,x_2):x_2\in
        X\}+1=\sup\{p(x_1,x_2):x_2<x_1\}+1$, where by $\alpha+1$ we
        mean the successor of an ordinal $\alpha$,
        we have $p_\Delta\in T_1$ and hence
        $\C\subseteq\cto$.

        To see the other inclusion, note first that the function
        $$
            q(x_1,x_2) =\begin{cases}p_\Delta(x_1,x_2)&, (x_1,x_2)\in\Delta\\
                 x_1&,\ow\end{cases}
        $$
        is in $\C$. Indeed, choose $\epsilon\in\Oo$ strictly
        increasing such that $p_\Delta(x_1,x_2)<\epsilon(x_1)$ for
        all $x_1,x_2\in X$ and consider
        $t(x_1,x_2)=p_\Delta(\epsilon(x_1),p_\Delta(x_1,x_2))$. On
        $\Delta$, $t$ is still one-one, and outside $\Delta$, the term
        is a one-one function of the first component $x_1$. Moreover, the ranges $t[\Delta]$ and
        $t[X^2\setminus\Delta]$ are disjoint. Hence,
        we can write $q=u\circ t$ for some unary $u$.
        By the same argument we see that for arbitrary unary functions $a,b\in\Oo$ the
        function
        $$
                q_{a,b}(x_1,x_2) =\begin{cases}a(p_\Delta(x_1,x_2))&, (x_1,x_2)\in\Delta\\
                 b(x_1)&,\ow\end{cases}
        $$
        is an element of $\C$.

        Now let $f\in T_1$ be given and say $f(x_1,x_2) < \delta(x_1)$ for all $x_1,x_2\in X$,
        where $\delta\in\Oo$ is strictly increasing. Choose $a\in\Oo$ such that $a(p_\Delta(x_1,x_2))=f(x_1,x_2)+1$ for all
        $(x_1,x_2)\in \Delta$. Then set
        $$
        f_1(x_1,x_2) =q_{a,\delta+1}(x_1,x_2)=\begin{cases}
              f(x_1,x_2)+1&,(x_1,x_2)\in\Delta\\
                \delta(x_1)+1           &,\ow
                \end{cases}
        $$
        We construct a second function
        $$
        f_2(x_1,x_2) =\begin{cases} 0&, (x_1,x_2)\in\Delta\\
                                f(x_1,x_2)+1&, \ow\end{cases}
        $$
        It is readily verified that
        $f_2(x_1,x_2) = u (p_\Delta(x_2+1, x_1))$ for some unary $u$.
        Now  $f_2(x_1,x_2) <  f_1(x_1,x_2)$ and $f_1,f_2\in\C$.
        Clearly $$f(x_1,x_2)  = u(p_\Delta( f_1(x_1,x_2) , f_2(x_1,x_2) ))$$  for some
        unary $u$. This shows $f\in\C$ and so $\cto\subseteq\C$ as
        $f\in T_1$ was arbitrary.
    \end{proof}

    We shall see that $\pto$ is also finitely generated
    over $\Oo$. Moreover, for countable $X$ it is a Borel set in the topology yet to be
    defined.
    Our explicit description $\pto$ holds for all infinite $X$ of regular cardinality, but
    is interesting only if there are not too many other maximal
    clones containing $\Oo$. By Theorem \ref{THM:theTwoMaximal}, this
    is at least the case for $X$ countably infinite or weakly
    compact.

    Throughout this chapter, the assumption that the base set $X$ has
    regular cardinality is essential. To give an example, we prove now that
    $\U$ is a clone. Let $f\in\U\un$ and $g_1,\ldots,g_n\in\U^{(m)}$. By
    definition, there exists $F\in\Oo$ and some $1\leq k\leq n$ such
    that $f(x)< F(x_k)$ for all $x\in X^n$. Because $g_k\in\U$, we
    obtain $G_k\in\Oo$ and $1\leq i\leq m$ such that $g_k(x)<
    G_k(x_i)$ for all $x\in X^m$. Therefore
    $f(g_1,\ldots,g_n)(x)<H(x_i)$, where we define
    $H(x_i)=\sup_{y<G_k(x_i)}\{F(y)\}$. Now since $F(y)<X$ for all
    $y\in X$, and since the supremum ranges over a set of size
    $G_k(x_i)<X$, the regularity of $X$ implies that
    $H(x_i)<X$, so that the composite $f(g_1,\ldots,g_n)$ is bounded
    by a unary function and hence an element of $\U$. It is easy to
    see that on singular $X$, neither of the definition of an
    almost unary function by means of small sets nor the one via
    boundedness by a unary function yield a clone. Also, the two
    definitions differ on singulars, whereas on regulars they
    coincide.
\subsection{Notation}

    For a set of functions $\F$ we shall denote the smallest
    clone containing $\F$ by $\langle \F \rangle$. By $\F^{(n)}$ we refer to the set
    of $n$-ary functions in $\F$. \\
    We call the projections which every clone contains $\pi^n_i$
    where $n\geq 1$ and $1\leq i\leq n$.
    If $f\in\On$ is an $n$-ary function, it sends $n$-tuples of elements of $X$ to $X$ and we
    write $(x_1,\ldots,x_n)$ for these tuples unless otherwise stated as in $f(x,y,z)$; this is the only place where
    we do not stick to set-theoretical notation (according to which we would have to write $(x_0,\ldots,x_{n-1})$).
    The set $\{1,\ldots,n\}$ of indices of $n$-tuples will play an important role and we write $N$
    for it. We denote the set-theoretical complement of
    a subset $A\subseteq N$ in $N$ by $-A$. We identify the set $X^n$ of $n$-tuples
    with the set of functions from $N$ to $X$, so that if $A\subseteq N$ and $a: A\To X$ and $b: -A\To X$ are
    partial functions, then $a\cup b$ is an $n$-tuple. Sometimes, if the arity of $f\in\OO$ has not yet been given a name,
    we refer to that arity by $n_f$.\\
    If $a\in X^n$ is an $n$-tuple and $1\leq k\leq n$ we write
    $(a)^n_k$ or only $a_k$ for the $k$-th component of $a$. For $c\in
    X$ and $J$ an index set we write $c^J$ for the $J$-tuple with
    constant value $c$. The
    order relation $\leq$ on $X$ induces the
    pointwise partial order on the set of $J$-tuples of elements of $X$ for any index
    set $J$: For $x,y\in X^J$ we write $x\leq y$ iff $x_j\leq y_j$ for all $j\in J$. Consequently we also
    denote the induced pointwise
    partial order of $\On$ by $\leq$, so that for
    $f,g\in\On$ we have $f\leq g$ iff $f(x)\leq g(x)$ for all $x \in
    X^n$. Whenever we state that a function $f\in\On$ is monotone, we
    mean it is monotone with respect to $\leq$: $f(x)\leq f(y)$
    whenever $x\leq y$. We denote the power set of $X$ by $\Pow(X)$. The
    element $0\in X$ is the smallest element of $X$.
\end{section}

\begin{section}{Properties of clones above $\U$ and the clone $\pto$}

\subsection{What $\cto$ is}
    We start by proving that the almost unary clone $\U$ is a so-called \emph{binary
    clone}, that is, it is generated by its binary part. Thus,
    when investigating $[\U,\pol(T_1)]$, we are in fact dealing
    with an interval of the form $[\cl{\C^{(2)}},\pol(\C^{(2)})]$ for $\C$ a clone.
    \begin{lem}
        The binary almost unary functions generate all almost unary
        functions. That is, $\cto=\U$.
    \end{lem}
    \begin{proof}
        Trivially, $\cto\subseteq \U$. Now we prove by induction
        that $\U^{(n)}\subseteq \cto$ for all $n\geq 1$. This is
        obvious for $n=1,2$. Assume we have $\U^{(k)}\subseteq
        \cto$ for all $k<n$ and take any function $f\in\U^{(n)}$.
        Say without loss of generality that $f(x_1,\ldots,x_n)\leq\gamma
        (x_1)$ for some $\gamma\in\Oo$. We will use the function
        $p_\Delta\in T_1$ to code two variables into one and then
        use the induction hypothesis. Define
        $$
            g_1(x_1,\ldots,x_{n-2},z)=
            \begin{cases}f(x_1,\ldots,x_{n-2},(p_\Delta^{-1}(z))^2_1,(p_\Delta^{-1}(z))^2_2)\quad
            &,z\in p_\Delta[X^2]\setminus\{0\} \\
            x_1\quad &,\ow \end{cases}
        $$
        The function is an element of $\U^{(n-1)}$ as it is bounded by
        $\max(x_1,\gamma(x_1))$. Intuitively, $g_1$ does the following: If $z\neq 0$ and in the range
        of $p_\Delta$, then
        $g_1$ imagines a pair $(x_{n-1},x_n)$ to be coded into $z$
        via $p_\Delta$. It reconstructs the pair $(x_{n-1},x_n)$
        and calculates $f(x_1,\ldots,x_n)$. If $z=0$ or not in the range of $p_\Delta$, then $g$ knows
        there is no information in $z$; it simply forgets about the tuple
        $(x_2,\ldots,x_n)$ and returns $x_1$, relying on the following similar function
        to do the job: Set $\Delta'=\Delta\cup\{(x,x):x\in X\}$
        and define
        $$
            g_2(x_1,\ldots,x_{n-2},z)=
            \begin{cases}f(x_1,\ldots,x_{n-2},(p_{\Delta'}^{-1}(z))^2_2,(p_{\Delta'}^{-1}(z))^2_1)\quad
            &,z\in p_{\Delta'}[X^2]\setminus\{0\}\\
            x_1\quad &,\ow\end{cases}
        $$
        The function $g_2$ does exactly the same as $g_1$ but assumes the
        pair $(x_{n-1},x_n)$ to be coded into $z$ in wrong
        order, namely as $(x_n,x_{n-1})$, plus it cares for the diagonal.
        Now consider
        $$
            h(x_1,\ldots,x_n)=g_2(g_1(x_1,\ldots,x_{n-2},p_\Delta(x_{n-1},x_n)),x_2,\ldots,x_{n-2},p_{\Delta'}(x_n,x_{n-1})).
        $$
        All functions which occur in $h$ are almost unary with at most $n-1$ variables.
        We claim that $h=f$. Indeed, if $x_{n-1}<x_n$, then
        $p_\Delta(x_{n-1},x_n)\neq 0$ and $g_1$ yields $f$. But
        $p_{\Delta'}(x_n,x_{n-1})=0$ and so $g_2$ returns $g_1=f$. If
        on the other hand $x_n\leq x_{n-1}$, then $p_\Delta(x_{n-1},x_n)=
        0$ and $g_1=x_1$, whereas $p_{\Delta'}(x_n,x_{n-1})\neq 0$,
        which implies $h=f(g_1,x_2,\ldots,x_n)=f(x_1,\ldots,x_n)$.
    \end{proof}
    The following lemma will be
    crucial for our investigation of clones containing $T_1$.
    \begin{cor}
        Let $\C$ be a clone containing $T_1$. Then $\C$ is
        downward closed, that is, if $f\in \C$, then also $g\in \C$
        for all $g\leq f$.
    \end{cor}
    \begin{proof}
        If $f\in \C^{(n)}$ and $g\in\On$ with $g\leq f$ are given, define
        $h_g(x_1,\ldots,x_{n+1})=\min(g(x_1,\ldots,x_n),x_{n+1})$. Then
        $h_g\leq x_{n+1}$ and consequently, $h_g\in\cto\subseteq\C$.
        Now $g=h_g(x_1,\ldots,x_n,f(x_1,\ldots,x_n))\in\C$.
    \end{proof}

\subsection{Wildness of functions}
    We have seen in the last section that the interval $\niceint$ is
    about growth of functions as all clones in that interval
    are downward closed. But mind we are not talking about how
    rapidly functions are growing in the sense of polynomial growth,
    exponential growth and so forth since we are considering
    clones modulo $\Oo$ (and so we can make functions as steep as
    we like); the growth of a function will be determined by which
    of its variables are responsible for the function to obtain many values.
    The following definition is due to M. Goldstern and S. Shelah
    \cite{GS022}. Recall that $N = \{1,\ldots,n\}$.
    \begin{defn}
        Let $f\in\On$. We call a set $\emptyset\neq A\subseteq N$ \emph{$f$-strong} iff for
        all $a\in X^A$ the set $\{f(a\cup x):x\in X^{-A}\}$ is
        small. $A$ is \emph{$f$-weak} iff it is not $f$-strong.
        In order to use the defined notions more freely, we define the
        empty set to be $f$-strong iff $f$ has small range.
    \end{defn}
    Thus, a set of indices of variables of $f$ is strong iff $f$ is
    bounded whenever those variables are. For example, a function is almost unary
    iff
    it has a one-element strong set. Here, we shall rather
    think in terms of the complements of weak sets.
    \begin{defn}
        Let $f\in\On$ and let $A\subsetneqq N$ and $a\in X^{-A}$. We
        say $A$ is \emph{$(f,a)$-wild} iff the set $\{f(a\cup
        x):x\in X^A\}$ is large. The set $A$ is called
        \emph{$f$-wild} iff there exists $a\in X^{-A}$ such that $A$ is
        $(f,a)$-wild. We say that $A$ is \emph{$f$-insane} iff $A$
        is $(f,a)$-wild for all $a\in X^{-A}$. The set $N$ itself
        we call $f$-wild and $f$-insane iff $f$ is unbounded.
    \end{defn}
    Observe that if $A\subseteq B\subseteq N$ and $A$ is $f$-wild, then $B$ is $f$-wild as well.
    Obviously, $A\subseteq N$ is $f$-wild iff $-A$ is
    $f$-weak. It is useful to state the following trivial criterion for a
    function to be almost unary.
    \begin{lem}\label{LEM:almostUnaryCriterion}
        Let $n\geq 2$ and $f\in\On$. $f$ is almost
        unary iff there exists a subset of $N$ with $n-1$ elements which is not $f$-wild.
    \end{lem}
    \begin{proof}
        If $f$ is almost unary, then there is a one-element $f$-strong subset of $N$ and
        the complement of that set is not $f$-wild. If on the other hand there exists
        $k\in N$ such that $N\setminus\{k\}$ is
        not $f$-wild, then $\{k\}$ is $f$-strong and so $f$ is almost unary.
    \end{proof}

    We will require the following fact from
    \cite{GS022}.
    \begin{fact} \label{FAC:goldsternShelah}
        If $f\in\pton$ and $A_1,A_2\subseteq N$ are $f$-wild, then
        $A_1\cap A_2\neq \emptyset$.
    \end{fact}
    We observe that the converse of this statement holds as well.
    \begin{lem}\label{charPto}
        Let $f\in\On$ be any $n$-ary function. If all pairs of
        $f$-wild subsets of $N$ have a nonempty intersection, then
        $f\in\pto$.
    \end{lem}
    \begin{proof}
        Let $g_1,\ldots,g_n\in T_1$ be given and set $A_1=\{k\in
        N:\exists\gamma\in\Oo\,(g_k(x_1,x_2)\leq\gamma(x_1))\}$ and
        $A_2=-A_1$. Since $A_1\cap A_2=\emptyset$ either $A_1$
        or $A_2$ cannot be $f$-wild. Thus
        $f(g_1,\ldots,g_n)$ is bounded by a unary function of $x_2$
        in the first case and by a unary function of $x_1$ in the
        second case.
    \end{proof}
    The equivalence yields a first description of
    $\pto$ with an interesting consequence.
    \begin{thm}\label{COR:FirstptoDescription}
        A function $f\in\On$ is an element of $\pto$ iff all pairs
        of $f$-wild subsets of $N$ have a nonempty intersection.
    \end{thm}

\subsection{Descriptive set theory}

    We show now that for countable $X$, this description implies that $\pto$ is a Borel set with respect to
    the natural topology on
    $\OO$. The reader not interested in the topic can skip
    this part and proceed directly to the next section.

    We first explain the very basics of descriptive set theory; for more details consult
    \cite{Kec95}. Let $\T=(T, \Upsilon)$ be a \emph{Polish space}, that is, a
    complete, metrizable, separable topological space. The \emph{Borel sets} of
    $\T$ are the smallest $\sigma$-algebra on $T$ which contains the
    open sets. These sets can be ordered according to their
    complexity: One starts by defining
    $\Sigma_1^0=\Upsilon\subseteq\P(T)$ to consist exactly of the open
    sets and $\Pi^0_1$ of the closed sets. Then one continues inductively for all $1< \alpha<\omega_1$ by
    setting $\Pi_\alpha^0$ to
    contain precisely the complements of $\Sigma_\alpha^0$ sets, and
    $\Sigma_{\alpha}^0$ to consist of all countable unions of sets which are elements of
    $\bigcup_{1\leq \delta <\alpha} \Pi_\delta^0$.
    The sequences
    $(\Sigma_\alpha^0)_{1\leq \alpha<\omega_1}$ and $(\Pi_\alpha^0)_{1\leq
    \alpha<\omega_1}$ are increasing and the union over either of the two
    sequences yields the Borel sets.

    Equip our base set $X=\omega$ with the discrete topology. Then the product space $\N=\omega^\omega=\Oo$
    is the so-called \emph{Baire space}. It is obvious that
    $\On=\omega^{\omega^n}$ is homeomorphic to $\N$. Examples of open sets in $\On$ are
    the $A_{x}^y=\{f\in\On:f(x)=y\}$, where $x\in X^n$ and $y\in
    X$; in fact, these sets form a subbasis of the topology of
    $\On$.
    $\OO=\bigcup_{n=1}^{\infty}\On$
    is the sum space of $\omega$ copies of
    $\N$: The open sets in $\O$ are those whose intersection
    with each $\On$ is open in $\On$. With this topology, $\O$ is a Polish space,
    and in fact again homeomorphic to $\N$.

    Since clones are subsets of $\O$, they can divided into
    \emph{Borel clones} and clones which are no Borel sets. In our
    case, we find that $\pol(T_1)$ is a very simple Borel set.

    \begin{thm}
        Let $X$ be countably infinite. Then $\pto$ is a Borel set
        in $\OO$.
    \end{thm}
    \begin{proof}
        We have to show that $\pto\un$ is Borel in $\On$ for each $n\geq 1$. By the preceding theorem,
        $$
            \pto\un=\{f\in\On: \forall A,B\subseteq N
            (A,B \,f\text{-wild} \To A\cap B\neq \emptyset )\}
        $$
        There are no (only finite) quantifiers in this definition
        except for those which might occur in the predicate of wildness (observe that $\exists$-quantifiers correspond to
        unions and $\forall$-quantifiers to intersections).
        Now
        $$
            A\subseteq N\,f\text{-wild} \gdw \exists a\in X^{-A}\forall
            k\in X\exists b\in X^A (f(a\cup b)>k)
        $$
        For fixed $A\subseteq N$, $a\in X^{-A}$, $k\in X$, and $b\in X^A$, the set of all
        functions in $\On$ for which $(f(a\cup b)>k)$ is open.
        Thus, the set of all $f\in\On$ for which $A$ is $f$-wild
        is of the form $\bigcup\bigcap\bigcup\, open$, and hence
        $\Sigma_3^0$ by counting of unions and negations. Observe that all
        unions which occur in the definition are countable.

        Since the predicate of wildness is negated in the definition of $\pto\un$,
        we conclude that $\pto\un$ is $\Pi_3^0$.
    \end{proof}

    It is readily verified that $\U$ (and hence, $T_1$) is a Borel set as
    well. This is interesting in connection with the following:

    Above the Borel sets of a Polish space, one can continue the
    hierarchy of complexity. The next level, $\Sigma^1_1$, comprises
    the so-called \emph{analytic sets}, which are the continuous
    images of Borel sets; the \emph{co-analytic sets} ($\Pi^1_1$) are
    the complements of analytic sets. It is easy to see that the clone
    generated by a Borel set of functions in $\O$ is an analytic set.
    Since $\Oo$ and all countable sets are Borel, every set which is
    countably generated over $\Oo$ is analytic. M. Goldstern showed in
    \cite{Gol02} that $T_2$ and $\pol(T_2)$ are relatively
    complicated:

    \begin{thm}
        Let $X$ be countably infinite. Then $T_2$ and $\pol(T_2)$ are co-analytic but not analytic in
        $\O$. Hence, neither of the two clones $\cl{T_2}$ and $\pol(T_2)$ is countably
        generated over $\Oo$.
    \end{thm}
\subsection{What wildness means}
    We wish to compare the wildness of functions. Write $S_N$ for the set of all permutations on $N$.
    \begin{defn}
        For $f,g\in\On$ we say that $f$ is \emph{as wild as $g$}
        and write $f \aw g$ iff there exists a permutation $\pi\in
        S_N$ such that $A$ is $f$-wild if and only
        if $\pi[A]$ is $g$-wild for all $A\subseteq N$. Moreover, $g$ is \emph{at least as wild
        as $f$} ($f\lw g$) iff there is a permutation $\pi\in S_N$ such that
        for all $f$-wild subsets $A\subseteq N$ the image $\pi[A]$
        of $A$ under $\pi$ is $g$-wild.
    \end{defn}

    \begin{lem}
        $\aw$ is an equivalence relation and $\lw$ a quasiorder extending $\leq$ on
        the set of $n$-ary functions $\On$.
    \end{lem}
    \begin{proof}
        We leave the verification of this to the reader.
    \end{proof}
    \begin{lem}
        Let $f,g\in\On$. Then $f\aw g$ iff $f\lw g$ and $g\lw f$.
    \end{lem}
    \begin{proof}
        It is clear that $f\lw g$ (and $g\lw f$) if $f\aw g$. Now
        assume $f\lw g$ and $g\lw f$. Then there are
        $\pi_1,\pi_2\in S_N$ which take $f$-wild and $g$-wild
        subsets of $N$ to $g$-wild and $f$-wild sets, respectively.

        Set $\pi=\pi_2\circ\pi_1$. Then $A$ is $f$-wild iff $\pi[A]$
        is $f$-wild for any subset $A$ of $N$: If $A$ is $f$-wild,
        then $\pi_1[A]$ is $g$-wild, then $\pi_2[\pi_1[A]]=\pi[A]$
        is $f$-wild. If on the other hand $\pi[A]$ is $f$-wild,
        then take $k\geq 1$ such that $\pi^k=id_N$ and observe
        that $\pi^{k-1}[\pi[A]]=\pi^k[A]=A$ is $f$-wild.

        Now we see that $A$ is $f$-wild iff $\pi_1[A]$ is $g$-wild for all $A\subseteq N$:
        If $\pi_1[A]$ is $g$-wild, then so is
        $\pi_2\circ\pi_1[A]=\pi[A]$ and so is $A$ by the preceding
        observation. Hence, the permutation $\pi_1$ shows that
        $f\aw g$.
    \end{proof}

    \begin{cor}
        Let $n\geq 1$. Then $\mathord{\leq_W}/\mathord{\aw}$ is a partial order on
        the $\aw$-equivalence classes of $\On$.
    \end{cor}

    \begin{nota}
        Let $f \in\On$. By $\fto{f}$ we mean $\clfto{f}$ from now
        on. $\fto{f}$ is the smallest clone containing $f$ as well
        as all almost unary functions.
    \end{nota}
    We are aiming for the following theorem which tells us why we invented wildness.
    \begin{thm}\label{THM:mainWild}
        Let $f,g\in \On$. If $f\lw g$, then $f\in \fto{g}$. In
        words, if $g$ is at least as wild as $f$, then it
        generates $f$ modulo $T_1$.
    \end{thm}
    \begin{cor}
        Let $f,g\in\On$. If $f\aw g$, then $\fto{f}=\fto{g}$.
    \end{cor}
    We split the proof of Theorem \ref{THM:mainWild} into a
    sequence of lemmas. In the next lemma we see that it does not matter which $a\in
    X^{-A}$ makes a set $A\subseteq N$ wild.
    \begin{lem}\label{zerowild}
        Let $g\in\On$. Then there exists $g'\in\fto{g}\un$ such
        that for all $A\subseteq N$ the following holds:  If $A$ is $g$-wild, then $A$
        is $(g',0^{-A})$-wild.
    \end{lem}
    \begin{proof}
        Fix for all $g$-wild $A\subseteq N$ a tuple $a_A\in
        X^{-A}$ such that $\{g(x\cup a_A):x\in X^A\}$ is large.
        For an $n$-tuple $(x_1,\ldots,x_n)$ write $P=P(x_1,\ldots,x_n)=\{l\in N:x_l\neq
        0\}$ for the set of indices of positive components in the tuple.
        Define for $1\leq i\leq n$ functions
        $$
            \gamma_i(x_1,\ldots,x_n)=\begin{cases}x_i&,x_i\neq 0 \vee
            P(x_1,\ldots,x_n)\,\,\text{not}\,\,
            g\text{-wild}\\(a_{P})_i&,\text{otherwise}\end{cases}
        $$
        In words, if the set $P$ of indices of positive components in $(x_1,\ldots,x_n)$ is a wild set,
        then the $\gamma_i$ leave those positive components alone and send the zero components
        to the respective values making $P$ wild. Otherwise, they act just like projections.
        It is obvious that $\gamma_i$ is almost unary, $1\leq i \leq
        n$. Set $g'=g(\gamma_1,\ldots,\gamma_n)\in\fto{g}$. To
        prove that $g'$ has the desired property, let $A\subseteq N$
        be $g$-wild. Choose any minimal $g$-wild $A'\subseteq A$. Then
        by the definition of wildness the set $\{g(x\cup a_{A'}):x\in
        X^{A'}\}$ is large. Take a large
        $B\subseteq X^{A'}$ such that the sequence $(g(x\cup
        a_{A'}): x\in B)$ is one-one. Select further a large $C\subseteq
        B$ such that each component in the sequence of tuples
        $(x:x\in C)$ is either constant or injective and such that $0$ does not occur
        in any of the injective components (it is a simple
        combinatorial fact that this is possible). If one of
        the components were constant, then $A'$
        would not be minimal $g$-wild; hence, all components are
        injective. Now we have
        $$
        \begin{aligned}
            |X|&=|\{g(x\cup a_{A'}): x\in C\}|\\&= |\{g'(x\cup
            0^{-A'}):x\in C\}|
            &\leq|\{g'(x\cup 0^{-A}):x\in X^A\}|
        \end{aligned}
        $$
        and so $A$ is $(g',0^{-A})$-wild.
    \end{proof}
    We prove that we can assume functions to be monotone.
    \begin{lem}\label{monotone}
        Let $g\in\On$. Then there exists $g''\in\fto{g}\un$ such
        that $g\leq g''$ and $g''$ is monotone with respect to the
        pointwise order $\leq$.
    \end{lem}
    \begin{proof}
    We will define a mapping $\gamma$ from $X^n$ to $X^n$ such that $\gamma_i=\pi^n_i\circ\gamma$ is almost unary
        for $1\leq i\leq n$ and such that
        $g''=g\circ\gamma$ has the desired property. We fix for every
        $g$-wild $A\subseteq N$ a sequence $(\alpha^A_\xi)_{\xi\in X}$
        of elements of $X^n$ so that all components of $\alpha^A_\xi$ which lie not in $A$ are constant and so that
        $(g(\alpha^A_\xi))_{\xi\in X}$ is
        monotone and unbounded.

        Let $x\in X^n$. The \textit{order type} of $x$ is the unique
        $n$-tuple $(j_1,\ldots,j_n)$ of indices in $N$ such that
        $\{j_1,\ldots,j_n\}=\{1,\ldots,n\}$ and such that
        $x_{j_1}\leq\ldots\leq x_{j_n}$ and such that $j_k<j_{k+1}$ whenever
        $x_{j_k}=x_{j_{k+1}}$. Let $1\leq k\leq n$ be the largest element with the property that the set
        $\{j_k,\ldots,j_n\}$ is $g$-wild. We call the set $\{j_k,\ldots,j_n\}$ the \textit{pushing
        set} $Push(x)$
        and $\{j_1,\ldots,j_{k-1}\}$ the \textit{holding set} of $x$ with
        respect to $g$.

        We define by transfinite recursion
        $$
            \gamma:\quad\begin{matrix} X^n &\To& X^n\\
            x&\mapsto& \alpha^{Push(x)}_{\lambda(x)} \end{matrix}
        $$
        where
        $$
            \lambda(x)=\min\{\xi:g(\alpha^{Push(x)}_\xi)\geq\sup(\{g''(y):y<x\}\cup\{g(x)\})\}.
        $$
        This looks worse than it is: We simply map $x$ to the
        first element of the sequence $(\alpha^{Push(x)}_\xi)_{\xi\in X}$ such that all values of $g''$
        already defined as well as $g(x)$ are topped. By
        definition, $g''=g\circ\gamma$ is monotone and $g\leq
        g''$. It only remains to prove that all $\gamma_i$, $1\leq
        i\leq n$, are almost unary to see that $g''\in\fto{g}$.

        Suppose not, and say that $\gamma_k$ is not almost unary for some $1\leq k\leq n$. Then there exists a
        value $c\in X$ and a sequence of $n$-tuples
        $(\beta_\xi)_{\xi\in X}$ with constant value $c$ in the $k$-th component
        such that $(\gamma_k(\beta_\xi))_{\xi\in X}$ is unbounded.
        Since there exist only finitely many order types of
        $n$-tuples, we can assume that all $\beta_\xi$ have the same
        order type $(j_1,\ldots,j_n)$; say without loss of generality $(j_1,\ldots,j_n)=(1,\ldots,n)$.
        Then all $\beta_\xi$ have the same pushing set $Push(\beta)$ of indices.
        If $k$ was an element of the holding set of the
        tuples $\beta_\xi$, then $(\gamma_k(\beta_\xi):\xi\in X)$ would be constant so that $k$ must be in $Push(\beta)$.
        Clearly, $(\lambda(\beta_\xi))_{\xi\in X}$
        has to be unbounded as otherwise $(\gamma_k(\beta_\xi))_{\xi\in X}$ would be bounded. Since by definition
        the value of $\lambda$
        increases only when it is necessary to keep $g\leq g''$,
        the set $\{g(y):\exists \xi\in X(y\leq \beta_\xi)\}$ is unbounded. But because of the order type of the $\beta_\xi$,
        whenever $i\leq k$, then we have $(\beta_\xi)^n_{i}\leq c$ for all
        $\xi\in X$ so that the components of the $\beta_\xi$ with index in the set
        $\{1,\ldots,k\}$ are bounded. Thus, $\{k+1,\ldots,n\}$
        is $g$-wild, contradicting the fact that $k$ is in the
        pushing set $Push(\beta)$.

    \end{proof}
    In a next step we shall see that modulo $T_1$, wildness is insanity.
    \begin{lem}\label{LEM:MonotoneAndInsane}
        Let $g\in\On$. Then there exists $g''\in\fto{g}\un$ such
        that $g''$ is monotone and for all $A\subseteq N$ the following holds:  If $A$ is $g$-wild, then $A$
        is $g''$-insane.
    \end{lem}
    \begin{proof}
        Let $g'\in\fto{g}\un$ be provided by Lemma \ref{zerowild} and make a monotone $g''$ out
        of it with the help of the preceding lemma. We claim that
        $g''$ already has both desired properties. To prove this,
        consider an arbitrary $g$-wild $A\subseteq N$. By construction
        of $g'$, $A$ is $(g',0^{-A})$-wild and so it is also
        $(g'',0^{-A})$-wild as $g'\leq g''$. But $0^{-A}\leq a$
        for all $a\in X^{-A}$; hence the fact that $g''$ is monotone implies
        that $A$ is $(g,a)$-wild for all $a\in X^{-A}$ which means
        exactly that $A$ is $g''$-insane.
    \end{proof}
    \begin{lem}\label{leq}
        Let $f, g\in\On$. If $f\lw g$, then there exists
        $h\in\fto{g}\un$ such that $f\leq h$.
    \end{lem}
    \begin{proof}
        Without loss of generality, we assume that the permutation $\pi\in S_N$ taking
        $f$-wild subsets of $N$ to $g$-wild sets is the identity on $N$.
        We take $g''\in\fto{g}$ according to the preceding lemma.
        We wish to define $\gamma\in\Oo$ with $f\leq\gamma\circ
        g''$. For $x\in X$ write $U_x=g''^{-1}[\{x\}]$ for the
        preimage of $x$ under $g''$. Now set
        $$
            \gamma(x)=\begin{cases}\sup\{f(y):y\in U_x\}&,U_x\neq
            \emptyset\\0&,\text{otherwise}\end{cases}
        $$
        We claim that $\gamma$ is well-defined, that is, the
        supremum in its definition always exists in $X$. For suppose there
        is an $x\in X$ such that the set $\{f(y):y\in U_x\}$ is
        unbounded. Choose a large subset $B\subseteq U_x$
        making the sequence $(f(y):y\in B)$ one-one. Take further
        a large
        $C\subseteq B$ so that all components in the sequence
        $(y:y\in C)$ are either one-one or constant. Set $A\subseteq N$ to consist of the
        indices of the injective components.
        Obviously, $A$ is $f$-wild; therefore it is $g''$-insane.
        Since $g''$ is also monotone, the set $\{g''(y):y\in C\}$
        is large, contradicting the fact that $g''$ is constant
        on $U_x$. Thus, $\gamma$ is well-defined and clearly
        $f\leq h\in\fto{g}$ where $h=\gamma\circ g''$.
    \end{proof}
    \begin{proof}[Proof of Theorem \ref{THM:mainWild}]
        The assertion is an immediate consequence of the preceding
        lemma and the fact that all clones above $\U$ are
        downward closed.
    \end{proof}
    \begin{rem}\label{REM:noConverse}
        Unfortunately, the converse does not hold: If $f,g\in\On$ and $f\in\fto{g}$
        then it need not be true that $f\lw g$. We will see an example at the end of the section.
    \end{rem}
\subsection{$\med_3$ and $T_1$ generate $\pto$}
    We are now ready to prove the explicit description of $\pto$.
    \begin{defn}\label{DEF:Mnk}
        For all $n\geq 1$ and all $1\leq k\leq n$ we define a function
        $$
            m^n_k(x_1,\ldots,x_n)=x_{j_k}\quad ,\text{if} \,\, x_{j_1} \leq \ldots \leq
            x_{j_n}.
        $$
        For example, $m^n_n$ is the maximum function $\max_n$
        and $m^n_1$ the minimum function $\min_n$ in $n$ variables. Note that
        $\min_n\in \pol(T_1)$ (it is even almost unary) but
        $\max_n\nin\pto$ (and hence $\fto{\max_n}=\OO$).
        If
        $n$ is an odd number then we call $m^n_{\frac{n+1}{2}}$ the
        $n$-th median function and denote this function by $\med_n$.
    \end{defn}
    For fixed odd $n$ it is easily verified (check the wild sets and apply Theorem \ref{COR:FirstptoDescription})
    that $\med_n$ it is the largest of
    the $m^n_k$ which still lies in $\pol(T_1)$: $m^n_k\in
    \pol(T_1)$ iff $k\leq \frac{n+1}{2}$. It is for this reason
    that we are interested in the median functions on our quest
    for a nice generating system of $\pol(T_1)$. As a consequence
    of the following theorem from the preceding chapter (Theorem  \ref{medians}) it does not matter which of the median
    functions we consider:
    \begin{thm}
        Let $k,n \geq 3$ be odd natural numbers. Then $\med_k\in \langle \{\med_n\} \rangle$. In
        other words, a clone contains either no median function or
        all median functions.
    \end{thm}

    The following lemma states that within the restrictions of
    functions of $\pto$ (Fact \ref{FAC:goldsternShelah}), we can construct functions of arbitrary
    wildness with the median.
    \begin{lem}\label{LEM:verkettetesSystem}
        Let $n\geq 1$ and let $\A=\{A_1,\ldots,A_k\}\subseteq \P(N)$ be a
        set of subsets of $N$ with the property that $A_i\cap A_j\neq \emptyset$ for all $1\leq i,j\leq k$.
        Then there exists monotone $t_\A\in\clf{\med_3}\un$ such that all
        members of $\A$ are $t_\A$-insane.
    \end{lem}
    \begin{proof}
        We prove this by induction over the size $k$ of $\A$. If
        $\A$ is empty there is nothing to show. If $k=1$, we can
        set $t_\A=\pi^n_i$, where $i$ is an arbitrary element of
        $A_1$. Then $A_1$ is obviously $t_\A$-insane. If $k=2$,
        then define $t_\A=\pi^n_i$, where $i\in A_1\cap A_2$ is
        arbitrary. Clearly, both
        $A_{1}$ and $A_{2}$ are $t_\A$-insane. Finally, assume
        $k\geq 3$. By induction hypothesis, there exist monotone terms
        $t_\B,t_\C,t_\D\in\clf{\med_3}\un$ for the sets $\B=\{A_1,\ldots,A_{k-1}\}$,
        $\C=\{A_1,\ldots,A_{k-2},A_k\}$ and $\D=\{A_{k-1},A_k\}$ such
        that all sets in $\B$ (and $\C,\D$ respectively) are $t_\B$-insane ($t_\C$-insane,
        $t_\D$-insane). Set
        $$
            t_\A=\med_3(t_\B,t_\C,t_\D).
        $$
        Then each $A_i$ is insane for two of the three terms in
        $\med_3$. Thus, if we fix the variables outside $A_{i}$
        to arbitrary values, then at least two of the three subterms
        in
        $\med_3$ are still unbounded and so is $t_\A$ by the
        monotonicity of its subterms.
        Hence, every $A_{i}$ is $t_\A$-insane, $1\leq i\leq k$. Obviously $t_\A$ is
        monotone.
    \end{proof}

    \begin{lem}\label{tf}
        Let $f\in\pton$. Then there exists $t_f\in\clf{\med_3}$
        such that $f\lw t_f$.
    \end{lem}
    \begin{proof}
        Write $\A=\{A_1,\ldots,A_k\}$ for the set of $f$-wild
        subsets of $N$. By Fact \ref{FAC:goldsternShelah}, $A_i\cap A_j\neq \emptyset$ for all
        $1\leq i,j\leq k$. Apply the preceding lemma to $\A$.
    \end{proof}
    \begin{thm}\label{ptoMain}
        $\pto=\fto{\med_3}$.
    \end{thm}
    \begin{proof}
        It is clear that $\pto\supseteq\fto{\med_3}$. On the
        other hand we have just seen that  if $f\in\pto$, then there
        exists $t_f\in\clf{\med_3}$ such that $f\lw t_f$, whence
        $f\in \fto{\med_3}$.
    \end{proof}
    \begin{cor}
        $\pto$ is the $\leq$-downward closure of the clone generated by
        $\med_3$ and the unary functions $\Oo$.
    \end{cor}
    \begin{proof}
        Given $f\in\pto$, by Lemma \ref{tf} there exists $t_f\in\clf{\med_3}$ such that $f\lw t_f$.
        By Lemma \ref{LEM:verkettetesSystem}, $t_f$ is monotone and each $t_f$-wild set is in
        fact even $t_f$-insane. Now one follows the proof of Lemma
        \ref{leq} to obtain $\gamma\in\Oo$ such that $f\leq
        \gamma\circ t_f$.
    \end{proof}
    \begin{cor}
        $\pto=\cl{\{\med_3,p_{\Delta}\}\cup\Oo}$. In particular,
        $\pto$ is finitely generated over the unary functions.
    \end{cor}
    \begin{proof}
        Remember that $\cl{\{p_\Delta\}\cup\Oo}=\cto$ (Fact \ref{FAC:pDelta}) and apply
        Theorem \ref{ptoMain}.
    \end{proof}
        Now we can give the example promised in Remark
        \ref{REM:noConverse}. Set
         $$
            g(x_1,\ldots,x_4)=\med_3(x_1,x_2,x_3)
         $$
         and
         $$
            f(x_1,\ldots,x_4)=\med_5(x_1,x_1,x_2,x_3,x_4).
         $$
         It is obvious that $\fto{g}=\fto{\med_3}=\pto$. Next observe that
         $\fto{f}\subseteq\fto{\med_5}=\pto$ and that
         $f(x_1,x_2,x_3,x_3)=\med_3$ which implies
         $\pto=\fto{\med_3}\subseteq\fto{f}$. Thus,
         $\fto{g}=\fto{f}$.
         Consider on the other hand the 2-element wild sets of the
         two functions: Exactly $\{1,2\},\{1,3\}$ and \{2,3\} are
         $g$-wild, and $\{1,2\},\{1,3\},\{1,4\}$ are the wild sets of two elements for $f$.
         Now the intersection of first group is empty, whereas
         the one of the second group is not; so there is no permutation
         of the set $\{1,2,3,4\}$ which takes the first group to
         the second or the other way. Hence, neither $f\lw g$ nor
         $g\lw f$.
\end{section}

\begin{section}{The interval $[\U,\O]$}
\subsection{A chain in the interval}
    Now we shall show that the open interval $(\cto,\pto)$ is not empty by
    exhibiting a countably infinite descending chain therein with
    intersection $\U$.
    \begin{nota}
        For a natural number $n\geq 2$, we write
        $\M_n=\cl{\{m^n_2\}\cup T_1}$.
    \end{nota}
     Observe that since $m_2^2=\ma_2\nin \pto$, Theorem \ref{THM:poltone} implies that $\M_2=\OO$.
     Moreover, $m_2^3=\med_3$ and hence, $\M_3=\pto$.
    \begin{lem}\label{kLeqn}
        Let $n\geq 2$. Then $\M_n^{(k)}=\U^{(k)}$ for all $1\leq k<n$. That is,
        all functions in $\M_n$ of arity less than $n$ are
        almost unary.
    \end{lem}
    \begin{proof}
        Given $n,k$ we show by induction over terms that if
        $t\in\M_n\uk$, then $t$ is almost unary. To start the
        induction we note that the only $k$-ary functions in the
        generating set of $\M_n$ are almost unary. Now assume
        $t=f(t_1,t_2)$, where $f\in T_1$ and $t_1,t_2\in\M_n\uk$.
        By induction hypothesis, $t_1$ and $t_2$ are almost
        unary and so is $t$ as the almost unary functions are closed under
        composition. Finally, say $t=m_2^n(t_1,\ldots,t_n)$,
        where the $t_i$ are almost unary $k$-ary functions, $1\leq i\leq n$.
        Since $k<n$, there exist $i,j\in N$ with $i\neq j$, $l\in\{1,\ldots,k\}$ and
        $\gamma,\delta\in\Oo$ such that $t_i\leq\gamma(x_l)$ and
        $t_j\leq\delta(x_l)$. Then, $t\leq
        \max(\gamma,\delta)(x_l)$ and so $t$ is almost unary as well.
    \end{proof}
    \begin{cor}\label{MnNinMn1}
        If $n\geq 2$, then $m_2^n\nin\M_{n+1}$. Consequently,
        $\M_n\nsubseteq\M_{n+1}$.
    \end{cor}
    \begin{lem}\label{Mn1InMn}
        If $n\geq 2$, then $m_2^{n+1}\in\M_{n}$. Consequently,
        $\M_{n+1}\subseteq\M_n$.
    \end{lem}
    \begin{proof}
        Set
        $$
            f(x_1,\ldots,x_{n+1})=m_2^n(x_1,\ldots,x_n)\in\M_n.
        $$
        Then every $n$-element subset of $\{1,\ldots,n+1\}$ is
        $f$-wild. Hence, $m_2^{n+1}\lw f$ and so
        $m_2^{n+1}\in\fto{f}\subseteq\M_n$.
    \end{proof}
    \begin{thm}
        The sequence $(\M_n)_{n\geq 2}$ forms a countably infinite descending chain:
        $$
        \OO=\M_2\supsetneqq\M_3=\pol(T_1)\supsetneqq\M_4\supsetneqq
        \ldots\supsetneqq\M_n\supsetneqq\M_{n+1}\supsetneqq \ldots
        $$
        Moreover,
        $$
            \bigcap_{n\geq 2} \M_n=\U.
        $$
    \end{thm}
    \begin{proof}
        The first statement follows from Corollary \ref{MnNinMn1}
        and Lemma \ref{Mn1InMn}. The second statement a direct consequence of
        Lemma \ref{kLeqn}.
    \end{proof}

\subsection{Finally, this is the interval}

    We will now prove that there are no more clones in the
    interval $\niceint$ than the ones we already exhibited. We
    first state a technical lemma.

    \begin{lem}
        Let $f\in\On$ be a monotone function such that all
        $f$-wild subsets of $N$ are $f$-insane. Define for $i,j\in N$ with $i\neq j$ functions
        $$
            f^{(i,j)}(x_1,\ldots,x_n)=f(x_1,\ldots,x_{i-1},x_j,x_{i+1},\ldots,x_n)
        $$
        which replace the $i$-th by the $j$-th component and calculate
        $f$. Then the following implications hold for all $f$-wild $A\subseteq N$ and all $i,j\in N$ with $i\neq j$:
        \begin{itemize}
        \item[(i)]{If $i\nin A$, then $A$ is
        $f^{(i,j)}$-insane.}
        \item[(ii)]{If $j\in A$, then $A$ is
        $f^{(i,j)}$-insane.}
        \end{itemize}
    \end{lem}
    \begin{proof}
        We have to show that if we fix the variables outside $A$ to constant values, then $f^{(i,j)}$ is still
        unbounded; because $f$ is monotone, we can assume all values are fixed to
        $0$. Fix a sequence $(\alpha_\xi:\xi\in X)$ of
        elements of $X^n$ such that all components outside $A$ are
        zero for all tuples of the sequence and such that
        $(f(\alpha_\xi):\xi\in X)$ is unbounded. Define a sequence
        of $n$-tuples $(\beta_\xi:\xi\in X)$ by
        $$
            (\beta_\xi)^n_k=\begin{cases}0&,k\nin
            A\\\xi&,\ow\end{cases}
        $$
        For each $\xi\in X$ there exist a $\lambda\in X$ such that
        $\alpha_\xi\leq\beta_\lambda$. Then $f(\alpha_\xi)\leq
        f(\beta_\lambda)$. In either of the cases (i) or
        (ii), $f(\beta_\lambda)\leq f^{(i,j)}(\beta_\lambda)$.
        Thus,
        $(f^{(i,j)}(\beta_\xi):\xi\in X)$ is unbounded.
    \end{proof}

    \begin{lem}
        Let $f\in\On$ not almost unary. Then there exists $n_0\geq
        2$ such that $\fto{f}=\fto{m^{n_0}_2}$.
    \end{lem}
    \begin{proof}
        We shall prove this by induction over the arity $n$ of
        $f$. If $n=1$, there are no not almost unary functions so there is nothing to
        show. Now assume our assertion holds for all $1\leq k<n$.
        We distinguish two cases:

        First, consider $f$ such that
        all $f$-wild subsets of $N$ have size at least $n-1$. Then
        $f\aw m^n_2$ and so $\fto{f}=\fto{m^{n}_2}$.

        Now assume
        there exists an $f$-wild subset of $N$ of size $n-2$, say
        without loss of generality that $\{2,\ldots,n-1\}$ is such a
        set. By Lemma \ref{LEM:MonotoneAndInsane} and Theorem \ref{THM:mainWild} there exists a monotone $\hat{f}$
        with $\fto{f}=\fto{\hat{f}}$ and with the property that all $f$-wild subsets of $N$ are
        $\hat{f}$-insane. Since we could replace $f$ by $\hat{f}$, we assume
        that $f$ is monotone and that all $f$-wild sets are
        $f$-insane.

        Consider the $f^{(i,j)}$ as defined in the preceding lemma. Formally, these functions
        are still $n$-ary, but in fact they depend only on $n-1$ variables. Thus, all
        of the $f^{(i,j)}$ which are not almost unary satisfy the induction
        hypothesis. Set
        $$
            n_0=\min\{k:\exists i,j\in N\,
            \fto{f^{(i,j)}}=\fto{m^{k}_2}\}.
        $$
        The minimum is well-defined: Because $\{2,\ldots,n-1\}$ is $f$-insane,
        $f^{(n,1)}$ is not almost unary so that it generates the same clone
        as some $m^n_2$ modulo $T_1$; thus, the set is
        not empty. Clearly, $m^{n_0}_2\in\fto{f}$. We show that $m_2^{n_0}$ is strong enough to
        generate $f$. Since $\M_n\subseteq\M_{n_0}$ for all $n\geq
        n_0$ we have $f^{(i,j)}\in\fto{m^{n_0}_2}$ for all $i,j\in
        N$ with $i\neq j$. Now define
        $$
            t(x_1,\ldots,x_n)=f^{(n,1)}(x_1,f^{(1,2)},f^{(1,3)},\ldots,f^{(1,n-1)})\in\fto{m^{n_0}_2}.
        $$
        We claim that $f\lw t$. Indeed, let $A\subseteq
        N$ be $f$-wild and whence $f$-insane by our assumption.

        If $1\nin A$, then $A$ is
        $f^{(1,j)}$-insane for all $2\leq j\leq n-1$ by the preceding lemma. So $A$ is
        insane for all components in the definition of $t$ except
        the first one. Hence, because $f$ is monotone, $A$ must be $t$-insane as
        otherwise $f^{(n,1)}$ would be almost unary.

        If $1\in A$, then by the preceding lemma $A$ is still $f^{(1,j)}$-insane whenever $j\in A$.
        Thus, increasing the components with
        index in $A$ increases the first component in $t$ plus all
        subterms $f^{(1,j)}$ with $j\in A$; but by the definition of $f^{(n,1)}$, that is the same
        as increasing the variables $A\cup\{n\}\supseteq
        A$ in $f$. Whence, $A$ is $t$-insane.

        This proves $f\lw t$ and thus $f\in\fto{m^{n_0}_2}$.
    \end{proof}

    So here it is, the interval and the end of our quest.

    \begin{thm}\label{THM:thisIsTheInterval}
        Let $\C\supsetneqq \U$ be a clone. Then there exists
        $n\geq 2$ such that $\C=\M_n$.
    \end{thm}

    \begin{proof}
        Set
        $$
            n_\C=\min\{n\geq 2:\M_n\subseteq\C\}.
        $$
        Since $\C$ contains a function which is not almost unary,
        the preceding lemma implies that
        the set over which we take the minimum is nonempty.
        Obviously, $\M_{n_\C}\subseteq\C$. Now let $f$ be an
        arbitrary function in $\C$ which is not almost unary. Then
        by the preceding lemma, there exists $n_0$ such that
        $\fto{m^{n_0}_2}=\fto{f}$. Clearly, $n_0\geq n_\C$ so that
        $f\in\M_{n_0}\subseteq\M_{n_\C}$.
    \end{proof}

    We state a lemma describing how the $k$-ary
    parts of the $\M_n$ for arbitrary $k$ relate to each other.
    \begin{lem}
         Let $m>n\geq 2$ and $k\geq 2$. If $k\geq n$ (that is, if $\M_n\uk$ is nontrivial), then
         $\M_n\uk\supsetneqq\M_m\uk$.
    \end{lem}
    \begin{proof}
        We know that $\M_n\uk\supseteq\M_m\uk$. To see the
        inequality of the two sets, observe that
        $$
            f(x_1,\ldots,x_k)=m^n_2(x_1,\ldots,x_n)
        $$
        is an element of $\M_n\uk$ but definitely not one of
        $\M_m\uk$.
    \end{proof}

    \begin{cor}
        Let $k\geq 2$. Then
        $$
        \M_2\uk\supsetneqq\M_3\uk\supsetneqq
        \ldots\supsetneqq\M_k\uk\supsetneqq\M_{k+1}\uk=\U\uk
        $$
        Consequently, there are $k$
        different $k$-ary parts of clones of the interval $[\U,\O]$ for each
        $k$.
    \end{cor}

    In general, if $\C$ is a clone, then
    $$
        \pol(\C^{(1)})\supseteq \pol(\C^{(2)})\supseteq \ldots \supseteq
        \pol(\C^{(n)})\supseteq \ldots
    $$
    Moreover,
    $$
        \pol(\C\un)\un=\C\un \quad\text{ and }\quad \bigcap_{n\geq 1} \pol(\C\un)=\C.
    $$
    It is natural to ask whether or not for $\C=\U$ this chain coincides with
    the chain we discovered.
    \begin{thm}
        Let $n\geq 1$. Then $\M_{n+1}=\pol(\U\un)$.
    \end{thm}
    \begin{proof}
        Clearly, $\M_2=\pol(\U^{(1)})=\O$, so assume $n\geq 2$.
        Consider $m^{n+1}_2$ and let $f_1,\ldots,f_{n+1}$ be
        functions in $\U\un$. Then two of the $f_j$ are bounded by
        unary functions of the same variable. Thus
        $m^{n+1}_2(f_1,\ldots,f_{n+1})$ is bounded by a unary function of
        this variable. This shows $m^{n+1}_2\in \pol(\U\un)$ and
        hence $\M_{n+1}\subseteq \pol(\U\un)$. Now consider $m^n_2$ and observe that
        $m^n_2\nin\U\un=\pol(\U\un)\un$; this proves $\M_{n}\nsubseteq
        \pol(\U\un)$. Whence, $\M_{n+1}=\pol(\U\un)$.
    \end{proof}
\subsection{The $m^n_k$ in the chain}
    As an example, we will show where the clones generated by the $m^n_k$ (as in Definition \ref{DEF:Mnk}) and $T_1$ can be found in
    the chain.
    \begin{nota}
        For $1\leq k\leq n$ we set $\M_n^k=\fto{m_k^n}$.
    \end{nota}
    Note that if $k=1$, then $\M_n^k=\U$, and if $k>\frac{n+1}{2}$,
    then
    $\M_n^k=\OO$. Observe also that $\M_n=\M^2_n$ for all $n\geq 2$.
    \begin{nota}
        For a positive rational number $q$ we write
        $$
          \lfloor q \rfloor = \max\{n\in\mathbb{N}:\,n\leq q\}
        $$
        and
        $$
            \lceil q \rceil = \min\{n\in\mathbb{N}:\,q\leq n\}.
        $$
        The remainder of the division $\frac{n}{k}$ we denote by the
        symbol $R(\frac{n}{k})$.
    \end{nota}
    \begin{lem}
        Let $2\leq k\leq \frac{n+1}{2}$ and let $t\in\M_n^k$ not almost unary. Then all $t$-wild
        subsets of $N_t$ have size at least
        $\frac{n}{k-1}-1$.
    \end{lem}
    \begin{proof}
        Our proof will be by induction over terms. If $t=m_k^n$,
        then all $t$-wild subsets of $N_t=N$ have at least $n-k+1$
        elements in accordance with our assertion. For the
        induction step, assume $t=f(t_1,t_2)$, where $f\in T_1$,
        say $f(x_1,x_2)\leq\gamma(x_1)$ for some $\gamma\in\Oo$.
        Then $t$ inherits the asserted property from $t_1$.
        Finally we consider the case where $t=m_k^n(t_1,\ldots,t_n)$.
        Suppose towards contradiction there exists $A\subseteq
        N_t$ $t$-wild with $|A|<\frac{n}{k-1}-1$.
        There have to be at least $n-k+1$ terms $t_j$ for which
        $A$ is $t_j$-wild so that $A$ can be $t$-wild. By
        induction hypothesis, these $n-k+1$ terms are almost
        unary and bounded by a unary function of a variable with index in $A$.
        From the bound on the size of $A$ we conclude that at least
        $$
            \lceil\frac{n-k+1}{|A|}\rceil>\frac{n-k+1}{\frac{n}{k-1}-1}=k-1
        $$
        of the terms $t_j$ are bounded by an unary function of the
        same variable with index in $A$. But if $k$ of the $t_j$ have the same one-element strong set,
        then $t$ is bounded by a unary function of this variable as
        well in contradiction to the assumption that $t$ is not
        almost unary.
    \end{proof}
    \begin{cor}\label{COR:Mnknsubseteq}
        Let $2\leq k\leq \frac{n+1}{2}$. Then $\M_{\lceil\frac{n}{k-1}\rceil-1}\nsubseteq\M_n^k$.
    \end{cor}
    \begin{proof}
        With the preceding lemma it is enough to observe that
        $m_2^{\lceil\frac{n}{k-1}\rceil-1}\in\M_{\lceil\frac{n}{k-1}\rceil-1}$
        has a wild set of size $\lceil\frac{n}{k-1}\rceil-2$.
    \end{proof}
    So we identify now the $\M_j$ which $\M_n^k$ is
    equal to.
    \begin{lem}\label{LEM:M_quotientFromMnk}
        Let $2\leq k\leq n$. Then $\M_{\lceil\frac{n}{k-1}\rceil}\subseteq\M_n^k$.
    \end{lem}
    \begin{proof}
        It suffices to show that $m_k^n$ generates
        $m_2^{\lceil\frac{n}{k-1}\rceil}$. But this is easy:
        $$
            m_2^{\lceil\frac{n}{k-1}\rceil}=m_k^n(x_1,\ldots,x_1,x_2,\ldots,x_2,\ldots,x_{\lceil\frac{n}{k-1}\rceil},\ldots,x_{\lceil\frac{n}{k-1}\rceil}),
        $$
        where $x_j$ occurs $k-1$ times if $1\leq j\leq
        \lfloor\frac{n}{k-1}\rfloor$ and $R(\frac{n}{k-1})<k-1$ times
        if $j=\lfloor\frac{n}{k-1}\rfloor+1$. For if we evaluate the function for a $\lceil\frac{n}{k-1}\rceil$-tuple with
        $x_{j_1}\leq
        \ldots\leq x_{j_{\lceil\frac{n}{k-1}\rceil}}$, then $x_{j_1}$
        occurs at most $k-1$ times in the tuple, but $x_{j_1}$
        together with $x_{j_2}$ occur more than $k$ times; thus,
        the $k$-th smallest element in the tuple is $x_{j_2}$ and
        $m_k^n$ returns $x_{j_2}$.
    \end{proof}
    \begin{thm}
        $\M_n^k=\M_{\lceil\frac{n}{k-1}\rceil}$ for all $2\leq k\leq
        n$.
    \end{thm}
    \begin{proof}
         By Theorem \ref{THM:thisIsTheInterval}, $\M_n^k$ has to be somewhere in the chain $(\M_n)_{n\geq 2}$.
         Because of Corollary \ref{COR:Mnknsubseteq} and Lemma
         \ref{LEM:M_quotientFromMnk} the assertion follows.
    \end{proof}
\subsection{Further on the chain}
    We conclude by giving one simple guideline for where to search
    the clone $\fto{f}$ in the chain for arbitrary $f\in \O$.

    \begin{lem}\label{kElWildSet}
        Let $1\leq k\leq n$ and let $f\in\On$ be a not almost unary function
        which has a $k$-element $f$-wild subset of $N$. Then
        $\M_{k+1}\subseteq\fto{f}$.
    \end{lem}
    \begin{proof}
        We can assume that $\{1,\ldots,k\}$ and all $A\subseteq N$ with $|A|=n-1$ are $f$-insane
        and that $f$ is monotone. Define
        $$
            g(x_1,\ldots,x_{k+1})=f(x_1,\ldots,x_k,x_{k+1},\ldots,x_{k+1})\in\fto{f}.
        $$
        Let $A\subseteq\{1,\ldots,k+1\}$ with $|A|=k$ be given. If
        $A=\{1,\ldots,k\}$ then $A$ is $f$-wild and so it is
        $g$-wild. Otherwise $A$ contains $k+1$ and so it affects $n-1$ components in the
        definition of $g$. Therefore $A$ is
        $g$-wild by Lemma \ref{LEM:almostUnaryCriterion}. Hence, $m_2^{k+1}\lw g$
        and so $\M_{k+1}\subseteq\fto{g}\subseteq\fto{f}$.
    \end{proof}
    \begin{rem}
        Certainly it is not true that if the smallest wild set of
        a function $f\in\O$ has $k$ elements, then
        $\M_{k+1}=\fto{f}$. The $m^n_k$ are an example.
    \end{rem}
    \begin{cor}\label{COR:2elwildset}
        Let $f\in\pto$ not almost unary and such that there exists a
        2-element $f$-wild subset of $N$. Then $\fto{f}=\pto$.
    \end{cor}
\end{section}

\begin{section}{Summary and a  nice picture}

    We summarize the main results of this chapter: For the interval of
    clones containing the almost unary functions we have
    $[\U,\O]=\{\M_2,\M_3,\ldots,\U\}$, where the
    $\M_n=\cl{\{m^n_2\}\cup\U}=\cl{\{m^n_2,p_\Delta\}\cup\Oo}$ are all
    finitely generated over $\Oo$. Alternatively, $\M_n$ can be
    described as the $\leq$-downward closure of
    $\cl{\{m^n_2\}\cup\Oo}$. The interval is a chain:
    $\M_2=\Oo\supsetneqq \M_3=\pto\supsetneqq \M_4\supsetneqq \ldots$ and
    $\bigcap_{n\geq 2} \M_n=\U$. Together with the fact that
    $\M_{n+1}=\pol(\U\un)$ for all $n\geq 1$ this yields that $\U$ an
    example of a clone $\C$ for which the chain
    $\pol(\C^{(1)})\supseteq \pol(\C\ut)\supseteq \ldots \supseteq\C$ is
    unrefinable and collapses nowhere. $\U$ is a so-called binary
    clone, that is, $\cl{\U\ut}=\U$.

    The $\M_n$ have the property that $\M_n\uk=\U\uk$ whenever $1\leq
    k < n$. Furthermore, $\M_n\uk\supsetneqq \M_m\uk$ whenever $m>n\geq 2$
    and $k\geq n$. Consequently, for each $k\geq 1$ there exist
    exactly
    $k$ different $k$-ary parts of clones of the interval $[\U,\O]$.

    Using wildness, a notion which completely determines a function
    modulo $\U$, it is possible to calculate for all $2\leq k\leq n$
    that
    $\M_n^k=\cl{\{m^n_k\}\cup\U}=\M_{\lceil\frac{n}{k-1}\rceil}$. In
    general, if one knows the wild subsets of $\{1,\ldots,n\}$ of a
    function $f\in\O\un$, he can draw certain conclusions about where
    to find the clone $\cl{\{f\}\cup\U}$ in the chain.

    On countable $X$, if we equip $\O$ with the natural topology, then the sets $T_1$ and $\pol(T_1)$ are
    Borel sets of low complexity, as opposed to the sets $T_2$ and
    $\pol(T_2)$ which have been shown by M. Goldstern to be
    non-analytic. In fact, with the results of this chapter, all
    clones above the almost unary functions can be shown to be
    Borel.

    If $X$ is countably infinite or weakly compact, we can draw the
    situation we ran into like this. \vspace*{4mm}

\newpage
\vspace*{10mm}
\newlength{\normalunitlength}
\setlength{\normalunitlength}{\unitlength}

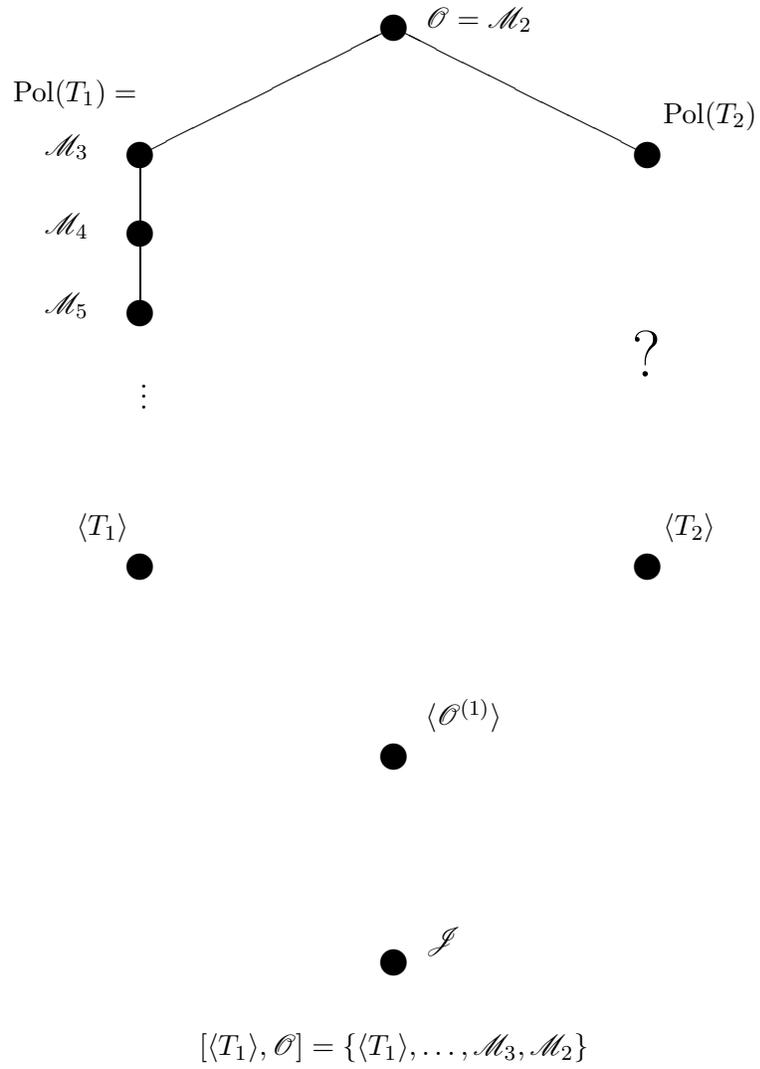
\setlength{\unitlength}{1.2\unitlength}
    \begin{center}
    \begin{picture}(200,300)
    \put(100,5){\circle*{8}}
    \put(110,10){$\J$}

    \put(100,70){\circle*{8}}
    \put(110,80){$\cl{\Oo}$}

    \put(100,300){\circle*{8}}
    \put(110,300){$\OO=\M_2$}

    \put(180,260){\circle*{8}}
    \put(185,270){$\pol(T_2)$}

    \put(20,130){\circle*{8}}
    \put(0,140){$\cto$}

    \put(180,130){\circle*{8}}
    \put(185,140){$\langle T_2\rangle$}

    \put(20,260){\line(2,1){80}}
    \put(180,260){\line(-2,1){80}}

    \put(175,190){\Huge{?}}

    \put(20,260){\circle*{8}}
    \put(-20,277){$\pto=$}
    \put(-10,260){$\M_3$}

    \put(20,210){\line(0,1){50}}

    \put(20,235){\circle*{8}}
    \put(-10,235){$\M_4$}

    \put(20,210){\circle*{8}}
    \put(-10,210){$\M_5$}

    \put(20,180){$\vdots$}

    \end{picture}
    \end{center}
    $$[\cto,\O]=\{\cto,\ldots,\M_3,\M_2\}$$
\setlength{\unitlength}{\normalunitlength}


\end{section}

\end{chapter}\newpage
\pagestyle{myheadings}\markright{\uppercase{Maximal clones above
the permutations}\hfill}
\begin{chapter}
{Maximal clones on uncountable sets that include all permutations}

We first determine the maximal clones on a set $X$ of infinite
regular cardinality $\kappa$ which contain all permutations but
not all unary functions, extending a result of L. Heindorf for
countably infinite $X$. If $\kappa$ is countably infinite or
weakly compact, this yields a list of all maximal clones
containing the permutations since in that case the maximal clones
above the unary functions are known. We then generalize a result
of G. Gavrilov to obtain on all infinite $X$ a list of all maximal
submonoids of the monoid of unary functions which contain the
permutations.

\begin{section}{Background and the results}

\subsection{Clones containing the bijections}

    Although the clone lattice on an infinite base set $X$
    need not be dually atomic by a result of M. Goldstern and S. Shelah \cite{GS03},
    the sublattice of $Cl(X)$ of clones containing the set $\S$ of all permutations of
    $X$ is dually atomic since $\OO$ is finitely generated over
    $\S$: Call a set $A\subseteq X$ \emph{large} iff $|A|=|X|=\kappa$ and \emph{small} otherwise.
    Moreover, $A$ is \emph{co-large} iff $X\sm A$ is large, and
    \emph{co-small} iff $X\sm A$ is small. Set
    $$
        \I=\{f\in\Oo: f \text{ is injective and } f[X]
        \text{ is co-large}\}
    $$
    and
    $$
        \J=\{g\in\Oo: g^{-1}[y] \text{ is large for all } y\in
        X\}.
    $$
    It is readily verified that for arbitrary fixed $f\in\I$ and
    $g\in\J$ we have
    $$
        \I=\{\alpha\circ f: \alpha\in\S\}\text{ and
        } \J=\{\alpha\circ g\circ\beta: \alpha,\beta\in\S\}.
    $$
    Moreover,
    $$
        \Oo=\{j\circ i: j\in \J,i\in\I\}.
    $$
    Together with the well-known fact that $\Oo\cup\{p\}$ generates $\O$ for any
    binary injection $p$ we conclude that $\O$ is generated by $\S\cup\{p,f,g\}$.
    Hence Zorn's lemma implies that the interval $[\S,\O]$ is
    dually atomic.

    We will determine all maximal
    clones $\C$ on a base set of regular cardinality for which $\S\subseteq\C$ but not $\Oo\subseteq\C$. This has
    already been done for countable base sets by L. Heindorf in the article
    \cite{Hei02} using the following concept: Let $\rho\subseteq X^J$ be a relation on $X$ indexed by $J$ and
    let
    $f\in\On$. We say that $f$ \emph{preserves} $\rho$ iff for all $r^1=(r^1_i:i\in J),\ldots,r^n=(r^n_i:i\in J)$ in
    $\rho$ we have $(f(r^1_i,\ldots,r^n_{i}):i\in J) \in \rho$. We define the clone of \emph{polymorphisms} $\pol(\rho)$
    of $\rho\subseteq X^J$ to consist exactly  of the
    functions in $\OO$ preserving $\rho$. In particular, if
    $\rho\subseteq X^{X^k}$ is a set of $k$-ary functions, then the
    polymorphisms of $\rho$ are exactly those $f\in\On$ for which the composite $f(g_1,\ldots,g_n)\in\rho$
    whenever $g_1,\ldots,g_n\in\rho$.
    It is obvious that since
    clones are closed under composition we have
    $\C\subseteq\p{\C\un}$ for any clone $\C$ and for all $n\geq 1$, where $\C\un=\C\cap\On$. Moreover,
    $\p{\C\un}\un=\C\un$. Therefore if $\C$ is a maximal clone such
    that $\S\subseteq\C\uo\subsetneqq\Oo$, then
    $\C\subseteq\p{\C\uo}\subsetneqq\O$ holds. Hence $\C=\p{\C\uo}$ by
    the maximality of $\C$. We conclude that all
    maximal clones with $\S\subseteq\C\uo\subsetneqq\Oo$ are of the form
    $\p{\G}$ where $\S\subseteq\G\subsetneqq\Oo$ is a \emph{submonoid} of
    $\Oo$, that is, a set of unary functions closed under composition and containing the identity map.

    \begin{thm}[L. Heindorf]\label{THM:heindorfThm}
        Let $X$ be a countably infinite set. The maximal clones
        over $X$ which contain all bijections but not all unary
        functions are exactly those of the form $\pol(\G)$, where
        $\G\in\{\A,\B,\D,\E,\F\}\cup\{\G_n: 1\leq n<\aleph_0\}$ is one of the following submonoids of $\Oo$:
        \begin{enumerate}
            \item{$\A=\{f\in\Oo:f^{-1}[\{y\}]$ is finite for
            almost all $y\in X\}$}
            \item{$\B=\{f\in\Oo:f^{-1}[\{y\}]$ is finite for all $y\in X\}$}
            \item{$\D=\{f\in\Oo:f$ is almost injective or not almost surjective$\}$}
            \item{$\E=\{f\in\Oo:f$ is almost surjective$\}$}
            \item{$\F=\{f\in\Oo:f$ is almost surjective or constant$\}$}
            \item{$\G_n=\{f\in\Oo:$ if $A\subseteq X$ has cardinality $n$ then $|X\setminus f[X\setminus A]|\geq n\}$}
        \end{enumerate}
        Consequently the number of such clones is countably infinite.
    \end{thm}
    In the theorem, ``almost all'' means ``all but finitely
    many'', ``almost injective'' means that there exists a finite
    subset $A$ of $X$ such that the restriction of $f$ to
    $X\setminus A$ is injective, and ``almost surjective'' means
    that the range of $f$ is co-finite.

    The restriction in the theorem to clones which do not contain all unary functions is not
    important since G. Gavrilov showed the following in \cite{Gav65}.

    \begin{thm}[G. Gavrilov]\label{THM:gavrilovTwoClones}
        If $X$ is countably infinite, then there
        exist exactly two maximal clones which contain $\Oo$.
    \end{thm}

    The two results imply that the number of maximal clones
    containing the permutations is countably infinite on a
    countably infinite base set.

    We now turn to base sets of any infinite cardinality.
    A property $P(y)$ holds for \emph{almost all} $y\in X$ iff
    the set of all elements for which the property does not hold
    is small. For $\lambda\leq \kappa$ a cardinal define a
    unary function $f$ to be \emph{$\lambda$-surjective} iff $|X\sm f[X]|<\lambda$. Instead of $\kappa$-surjective we also
    say \emph{almost surjective}; this means that the range of $f$ is co-small. $f$ is
    \emph{$\lambda$-injective} iff
    $|\{x\in X:\exists y\neq x\, (f(x)=f(y))\}|<\lambda$. For $\lambda=1$ or infinite, this is
    the case iff there exists a set $A\subseteq X$ such that
    $|A|<\lambda$ and such that the restriction of $f$ to the
    complement of $A$ is injective. \emph{Almost injective} means
    $\kappa$-injective.

    We are going to prove
    \begin{thm}\label{THM:bijections}
        Let $X$ be a set of regular cardinality $\kappa$. The maximal clones
        over $X$ which contain all bijections but not all unary
        functions are exactly those of the form $\pol(\G)$, where
        $\G\in\{\A,\B,\E,\F\}\cup\{\G_\lambda: 1\leq\lambda\leq\kappa,\,\lambda\text{ a cardinal}\}$ is one of the following submonoids of $\Oo$:
        \begin{enumerate}
            \item{$\A=\{f\in\Oo:f^{-1}[\{y\}]$ is small for
            almost all $y\in X\}$}
            \item{$\B=\{f\in\Oo:f^{-1}[\{y\}]$ is small for all $y\in X\}$}
            \item{$\E=\{f\in\Oo:f$ is almost surjective$\}$}
            \item{$\F=\{f\in\Oo:f$ is almost surjective or constant$\}$}
            \item{$\G_\lambda=\{f\in\Oo:$ if $A\subseteq X$ has cardinality $\lambda$ then $|X\setminus f[X\setminus A]|\geq \lambda\}$}
        \end{enumerate}
    \end{thm}

    \begin{cor}\label{COR:numberOfMaximal}
        Let $X$ be a set of regular cardinality $\kappa=\aleph_\alpha$.
        Then there exist $\max(|\alpha|,\aleph_0)$ maximal clones on $X$
        which contain all bijections but not all unary functions.
    \end{cor}

    Unfortunately, we do not know the maximal clones above $\Oo$
    on all regular cardinals; however we do on some. Let $\kappa$ be a cardinal. The \emph{partition symbol}
    $\kappa\To (\kappa)_2^2$
    means: Whenever the edges of a complete graph with
    $\kappa$ vertices are colored with 2 colors, then there
    is a complete subgraph with $\kappa$ vertices, all of
    whose edges have the same color. An uncountable $\kappa$
    for which the partition relation $\kappa\To
    (\kappa)_2^2$ holds is called \emph{weakly compact}.
    For example, the well-known
    theorem of F. Ramsey says that the defined partition relation holds for
    $\aleph_0$: If $G$ is a complete countably infinite graph and we color its
    edges with two colors, then there is an infinite
    complete subgraph of $G$ on which the coloring is constant.

    M. Goldstern and S. Shelah \cite{GS022} extended G. Gavrilov's result on maximal clones containing $\Oo$ to
    weakly compact cardinals.
    \begin{thm}[M. Goldstern and S. Shelah]\label{THM:weaklyCompact}
        If $\kappa=|X|$ is a weakly compact cardinal, then there
        exist exactly two maximal clones on $X$ which contain $\Oo$.
    \end{thm}

    Hence in the case of a weakly compact base set we know all maximal clones containing the
    permutations. It is a fact that weakly compact cardinals $\kappa$
    satisfy $\kappa=\aleph_\kappa$. Thus Corollary \ref{COR:numberOfMaximal} and Theorem \ref{THM:weaklyCompact}
    imply

    \begin{cor}
        Let $X$ be a set of weakly compact cardinality $\kappa$.
        Then there exist $\kappa$ maximal clones
        which contain all bijections.
    \end{cor}

    Unfortunately things are not always that easy.

    \begin{thm}[M. Goldstern and S. Shelah \cite{GS022}]
        For many regular cardinalities of $X$, in particular for all successors of regulars,
        there exist $2^{2^{|X|}}$ maximal clones which contain $\Oo$.
    \end{thm}

    It is interesting that whereas above $\Oo$ the
    number of maximal clones varies heavily with the partition
    properties of the underlying base set (2 for weakly compact
    cardinals, $2^{2^\kappa}$ for many others), the number of
    maximal clones above the permutations but not above $\Oo$ is a
    monotone function of $\kappa$ and always relatively small
    ($\leq\kappa$).

\subsection{Maximal submonoids of $\Oo$}
    Not all monoids appearing in Theorem
    \ref{THM:bijections} are maximal submonoids of $\Oo$.
    More surprisingly, there exist maximal submonoids of $\Oo$ above the
    permutations whose polymorphism clone is not maximal. Observe that submonoids of $\Oo$
    are simply \emph{unary clones}, that is clones consisting
    only of essentially unary functions, and that the lattice of
    monoids which contain the permutations is dually atomic by the
    argument we have seen before.

    \begin{thm}[G. Gavrilov \cite{Gav65}]\label{THM:GavMaximalMonoids}
        On a countably infinite base set $X$ the maximal submonoids of $\Oo$ containing the
        permutations are precisely the monoids $\A$, $\D$, $\G_1$,
        $\M$ and $\N$, where
        $$
            \M=\{f\in\Oo:f \text{ is surjective or not injective}\}
        $$
        and
        $$
            \N=\{f\in\Oo:f \text{ is almost surjective or not almost
            injective}\}.
        $$
    \end{thm}

    We will generalize
    Theorem \ref{THM:GavMaximalMonoids} to arbitrary infinite sets in the last
    section, obtaining
    \begin{thm}\label{THM:allMaximalMonoids}
        Let $X$ be an infinite set. If $X$ has regular cardinality, then the
        maximal submonoids of $\Oo$ which contain the permutations
        are exactly the monoid $\A$ and the monoids $\G_\lambda$ and $\M_\lambda$ for
        $\lambda=1$ and $\aleph_0\leq \lambda\leq\kappa$, $\lambda$ a cardinal, where
        $$
            \M_\lambda=\{f\in\Oo: f \text{ is
            }\lambda\text{-surjective or not
            }\lambda\text{-injective}\}.
        $$
        If $X$ has singular cardinality, then the same is true
        with the monoid $\A$ replaced by
        $$
            \A'=\{f\in\Oo: \exists \lambda < \kappa \,\,(\,|f\inv[\{x\}]|\leq\lambda\text{ for almost all }x\in X\,)\,\}.
        $$
    \end{thm}
    \begin{cor}\label{COR:numberOfMaximalMonoids}
        On a set $X$ of infinite cardinality $\aleph_\alpha$ there
        exist $2\,|\alpha|+5$ maximal submonoids of $\Oo$ that contain
        the permutations. Hence the smallest cardinality on which there
        are infinitely many such monoids is $\aleph_{\omega}$.
    \end{cor}

    Observe that the statement about singular cardinals in Theorem \ref{THM:allMaximalMonoids} differs only
    slightly from the corresponding one for regulars.
    We do not know whether Theorem \ref{THM:bijections} can be generalized to singulars, but
    in our proof we do use the
    regularity condition.
\subsection{Where has $\D$ gone?}
One might ask why in the general Theorems \ref{THM:bijections} and
\ref{THM:allMaximalMonoids} there is no monoid $\D$ as in Theorems
\ref{THM:heindorfThm} and \ref{THM:GavMaximalMonoids}. The answer
to that question is the following: Define for $\lambda=1$ and for
all $\aleph_0\leq\lambda \leq\kappa$ monoids
$$
    \delta(\lambda)=\{f\in\Oo: f \text{ is }\lambda\text{-injective or not }\lambda\text{-surjective}\}
$$
(this definition is due to I. Rosenberg \cite{Ros74}). Then we
have
\begin{lem}\label{LEM:whereHasDgone}
    $\delta(\lambda)=\G_\lambda$ for $\lambda=1$ and
    $\aleph_0\leq\lambda\leq\kappa$. In particular,
    $\D=\delta(\kappa)=\G_\kappa$.
\end{lem}
\begin{proof}
    Note that for $\lambda=1$, $\lambda$-injective simply means
    injective and $\lambda$-surjective means surjective. The lemma
    is easily verified for that case, and we prove it for
    $\lambda$ infinite.\\
    Assuming $f\in\delta(\lambda)$ we show $f\in\G_\lambda$. It is
    clear that if $f$ is not $\lambda$-surjective, then
    $f\in\G_\lambda$. So assume $f$ is $\lambda$-surjective; then
    by the definition of $\delta(\lambda)$, $f$ is
    $\lambda$-injective. Now let $A\subseteq X$ be an arbitrary set of
    size $\lambda$. Assume towards contradiction that $|X\sm f[X\sm
    A]|<\lambda$. Then two things can happen: If $|f[A]\cap f[X\sm
    A]|\geq\lambda$, then $|\{x\in X:\exists y\neq
    x\,(f(x)=f(y))\}|\geq |\{x\in A:\exists y\in X\sm A\,
    (f(x)=f(y))\}|\geq\lambda$, contradicting the
    $\lambda$-injectivity of $f$. Otherwise, $A$ is mapped onto a
    set of size smaller than $\lambda$, again in contradiction to $f$ being
    $\lambda$-injective.\\
    To see the other inclusion, take any $f\nin\delta(\lambda)$.
    Then $f$ is not $\lambda$-injective; thus we can find $A\subseteq
    X$ of size $\lambda$ such that $f[X]=f[X\sm A]$. But then
    $|X\sm f[X\sm A]|$=$|X\sm f[X]|<\lambda$ as $f$ is
    $\lambda$-surjective. Hence, $f\nin\G_\lambda$.
\end{proof}

 Before we start with the proofs we fix some notation.
\subsection{Notation}
    For a set of functions $\F$ we shall denote the smallest
    clone containing $\F$ by $\langle \F \rangle$. We call the projections which every clone contains $\pi^n_i$
    where $n\geq 1$ and $1\leq i\leq n$.
    We write $n_f$ for the arity of a function
    $f\in\O$ whenever that arity has not been given another name. If $a\in X^n$ is an
    $n$-tuple and $1\leq k\leq n$ we write $a_k$ for the $k$-th component of
    $a$. The image of a set $A\subseteq X^n$ under a function $f\in\On$ we
    denote by $f[A]$. Similarly we write $f\inv[A]$ for the
    preimage of $A\subseteq X$ under $f$. If $A=\{c\}$ is a singleton we
    cut short and write $f\inv[c]$ rather than $f\inv[\{c\}]$.
    Occasionally we shall denote the constant function with value
    $c\in X$ also by $c$.
    Whenever we identify $X$ with its cardinality we let $<$ and
    $\leq$ refer to the canonical well-order on $X$.
\end{section}

\begin{section}{The proof of Theorem \ref{THM:bijections}}

    In this section we are going to prove Theorem
    \ref{THM:bijections}; it will be the direct consequence of Propositions \ref{PROP:PolAisMaximal}, \ref{PROP:PolGUnderPolA},
    \ref{PROP:PolGUnderPolB}, \ref{PROP:someIlambdaIsBelow},
    \ref{PROP:PolGUnderPolGKappa},
    \ref{PROP:PolGlambdaIsMax}, \ref{PROP:GunderE}, \ref{PROP:GunderF}, and \ref{PROP:remainingPolsNotMaximal}.
    The first part of the proof (Section \ref{SEC:asInCountable}) is not much more than a translation of
    L. Heindorf's paper \cite{Hei02} to arbitrary regular cardinals; the
    reader familiar with that article should not be surprised to find
    the same constructions here. In Section \ref{SEC:almostUnary} we
    generalize a completeness criterion due to G. Gavrilov from
    countable sets to the uncountable to finish the proof.

\subsection{The core of the proof}\label{SEC:asInCountable}

We start with a general observation which will be useful.
\begin{lem}\label{LEM:unaryFunctionsSuffice}
    Let $\G$ be a proper submonoid of $\Oo$ such that
    $\cl{\pol(\G)\cup\{h\}}=\O$ for all unary $h\nin\G$. Then
    $\pol(\G)$ is maximal.
\end{lem}
\begin{proof}
    Let $f\nin \pol(\G)$ be given. Then there exist
    $h_1,\ldots,h_{n_f}\in\G$ such that $h=f(h_1,\ldots,h_{n_f})\nin\G$.
    Now $h\in\cl{\G\cup\{f\}}\subseteq\cl{\pol(\G)\cup\{f\}}$ and
    $\cl{\pol(\G)\cup\{h\}}=\O$ by assumption so that we conclude
    $\cl{\pol(\G)\cup\{f\}}=\O$.
\end{proof}

\subsubsection{The monoids $\A$ and $\B$}
\begin{prop}\label{PROP:PolAisMaximal}
    The clones $\p{\A}$ and $\p{\B}$ are maximal.
\end{prop}
\begin{proof}
    The maximality of $\p{\A}$ has been proved in
    \cite{Gav65} for the countable case and in \cite{Ros74} (Proposition 4.1) for
    arbitrary infinite sets.

    For the maximality of $\p{\B}$, let a unary $h\nin\B$ be given; by
    Lemma \ref{LEM:unaryFunctionsSuffice}, it suffices to show
    $\cl{\pol(\G)\cup\{h\}}=\O$. By the definition of $\B$ there
    exists $c\in X$ such that the preimage $Y=h\inv [c]$ is large.
    Choose any injection $g:X\To Y$; then $h\circ g(x)=c$ for all $x\in X$.

    Now let $f\in\On$ be an arbitrary function and consider
    $\tilde{f}\in\O^{n+1}$ defined by
    $$
        \tilde{f}(x_1,\ldots,x_n,y)=\begin{cases}f(x_1,\ldots,x_n)&,y=c\\y&,y\neq
        c\\\end{cases}.
    $$
    We claim that $\tilde{f}\in\p{\B}$. For let
    $\alpha_1,\ldots,\alpha_n,\beta\in\B$ and $d\in X$ be given. If
    $\tilde{f}(\alpha_1,\ldots,\alpha_n,\beta)(x)=d$, then by the definition of $\tilde{f}$ either
    $\beta(x)=c$ and $f(\alpha_1(x),\ldots,\alpha_n(x))=d$ or $\beta(x)\neq
    c$ and $\beta(x)=d$. But since $\beta\in\B$, the set of all
    $x\in X$ such that $\beta(x)=c$ or $\beta(x)=d$ is small.
    Hence $\tilde{f}(\alpha_1,\ldots,\alpha_n,\beta)\inv[d]$ is small and so
    $\tilde{f}(\alpha_1,\ldots,\alpha_n,\beta)\in\B$.

    Now to finish the proof it is enough to observe that
    $f(x_1,\ldots,x_n)=\tilde{f}(x_1,\ldots,x_n,c)=
    \tilde{f}(x_1,\ldots,x_n,h\circ g(x_1))\in\cl{\p{\B}\cup\{h\}}$.
\end{proof}
    We will prove now that $\B$ is the only proper submonoid of
    $\A$ whose $\pol$ is maximal. We start with a lemma.

\begin{lem}\label{LEM:dasGammaLemma}
    If $f\nin\p{\A}$, then there exist $\alpha_1,\ldots,\alpha\nf\in\Oo$ constant or injective such that
    $f(\alpha_1,\ldots,\alpha\nf)\nin\A$.
\end{lem}
\begin{proof}
    Since $f\nin\p{\A}$, there exist $\beta_1,\ldots,\beta\nf\in\A$
    such that $f(\beta_1,\ldots,\beta\nf)\nin\A$. We will use induction over
    $n_f$. If $n_f=1$, then $f\nin\p{\A}\uo=\A$ so that
    $f(\pi^1_1)=f \nin\A$ which proves the assertion for that
    case.
    Now assume the lemma holds for all functions of arity at most
    $n_f-1$. Define for $1\leq i\leq n_f$ sets $B_i=\{y\in
    X:\beta_i\inv[y] \text{ is large}\}$. By definition of $\A$, all
    $B_i$ are small. Set
    $$
        \Gamma=(\beta_1,\ldots,\beta\nf)[X]\setminus\prod_{1\leq
        i\leq n_f}B_i\subseteq X^{n_f}
    $$
    \textbf{Claim.} There exists a large set $D\subseteq X$ such that $f\inv [d]\cap \Gamma$ is large for all $d\in
    D$.\\
    To prove the claim, set $D=\{d\in X: f(\beta_1,\ldots,\beta\nf)\inv[d]\text{ large}\}\setminus
    f[\prod_{1\leq i\leq n_f} B_i]$. The set $D$ is large as $f(\beta_1,\ldots,\beta\nf)\nin\A$ and as $\prod_{1\leq i\leq n_f}
    B_i$ is small. Define $A_d=(f(\beta_1,\ldots,\beta\nf))\inv[d]$ for each $d\in
    D$. Then $(\beta_1,\ldots,\beta\nf)[A_d]\subseteq \Gamma$ is
    large for all $d\in D$. Indeed, assume to the contrary that there exists $d\in D$ such that
    $(\beta_1,\ldots,\beta\nf)[A_d]$ is small; then, since $|X|=\kappa$ is regular, there is an $x\in
    (\beta_1,\ldots,\beta\nf)[A_d]$ so that
    $(\beta_1,\ldots,\beta\nf)\inv[x]$ is large. But then we would
    have $x\in\prod_{1\leq i\leq n_f}B_i$, in contradiction to the
    assumption that $d\nin f[\prod_{1\leq i\leq n_f} B_i]$.
    This proves the claim since $f\inv [d]\cap \Gamma=(\beta_1,\ldots,\beta\nf)[A_d]$ is
    large for every $d\in D$.

    Setting $H^i_b=\{x\in X^{n_f}: x_i=b\}$ for all $1\leq i\leq n_f$ and all $b\in X$, we can
    write $\Gamma$ as follows:
    $$
        \Gamma=(\bigcup_{i=1}^{n_f}\bigcup_{b\in
        B_i}\Gamma\cap H^i_b)\cup(\Gamma\setminus \Delta),
    $$
    where $\Delta=\bigcup_{i=1}^{n_f}\bigcup_{b\in
        B_i} H^i_b$. Since $\kappa$ is regular and the union consists only of a small number of sets, we have that either there
    exist $1\leq i\leq n_f$ and some $b\in B_i$ such that $f\inv[d]\cap\Gamma\cap H^i_b$ is large for a large set
    of $d\in D$, or
    $f\inv[d]\cap\Gamma\setminus \Delta$ is large for a large set of $d\in D$. We distinguish
        the two cases:\\
    \textbf{Case 1.} There exist $1\leq i\leq n_f$ and $b\in B_i$
    such that $f\inv [d]\cap\Gamma\cap H^i_b$ is large for many $d\in D$; say
    without loss of generality $i=n_f$. Then
    $f(\beta_1,\ldots,\beta_{n_f-1},b)\nin \A$. By induction
    hypothesis, there exist $\alpha_1,\ldots,\alpha_{n_f-1}$ injective
    or constant such that $f(\alpha_1,\ldots,\alpha_{n_f-1},b)\nin\A$.
    Setting $\alpha\nf(x)=b$ for all $x\in X$ proves the lemma.\\
    \textbf{Case 2.} $f\inv [d]\cap\Delta$ is large for many $d\in D$.
    Observe that for all $a\in X$ and all $1\leq i\leq n_f$, $\Delta\cap
        H^i_a$ is small, for otherwise $\beta_i\inv[a]$ would be
        large and thus $a\in B_i$, contradiction. Set
    $$
        C=\{c\in X: f\inv[c]\cap\Delta\text{ large}\}.
    $$
    By the assumption for this case, $C$ is large. Now fix any $g:
    X\To C$ such that $g\inv[c]$ is large for all $c\in C$. We
    define a function $\alpha:X\To\Delta$ such that
    $f\circ\alpha=g$; moreover, $\alpha_i=\pi^{n_f}_i\circ\alpha$ will be injective, $1\leq i\leq n_f$.
    Identify $X$ with its cardinality
    $\kappa$. Then all $\alpha_i$ are injective iff
    $\alpha_i(x)\neq\alpha_i(y)$ for all $y< x$ and all $1\leq i \leq
    n_f$. This is the case iff
    $$
        (\alpha_1,\ldots,\alpha\nf)(x)\in\Delta\setminus\bigcup_{y<x}\bigcup_{i=1}^{n_f}
        H^i_{\alpha_i(y)}.
    $$
    Using transfinite induction on $\kappa$, we define $(\alpha_1,\ldots,\alpha\nf)$ by picking
    $$
        (\alpha_1,\ldots,\alpha\nf)(x)\in
        (f\inv[g(x)]\cap\Delta)\setminus\bigcup_{y<x}\bigcup_{i=1}^{n_f}
        H^i_{\alpha_i(y)}.
    $$
    This is possible as $f\inv[g(x)]\cap\Delta$ is large for all
    $x\in X$
    whereas $\Delta\cap\bigcup_{y<x}\bigcup_{i=1}^{n_f}
    H^i_{\alpha_i(y)}$ is small. Clearly
    $f(\alpha_1,\ldots,\alpha\nf)=g\nin \A$ and the proof of the
    lemma is complete.
\end{proof}
\begin{prop}\label{PROP:PolGUnderPolA}
    Let $\G\subseteq \A$ be a submonoid of $\Oo$ which contains
    all permutations. Then either $\G\subseteq\B$ or
    $\p{\G}\subseteq\p{\A}$.
\end{prop}
\begin{proof}
    Assume $\G \nsubseteq \B$; we show $\p{\G}\subseteq\p{\A}$. Observe first that
    for all co-large $A\subseteq X$ and all $a\in
    X$ there exists $g\in\G$ such that $g[A]=\{a\}$.
    Indeed, choose any $h\in\G\setminus\B$. There exists $y\in X$
    such that $h\inv[y]$ is large. Choose bijections
    $\alpha,\beta\in\S$ with the property that $\alpha[A]\subseteq
    h\inv[y]$ and that $\beta(y)=a$. Then $g=\beta\circ
    h\circ\alpha$ has the desired property.\\
    Now let $f\nin\p{\A}$ be arbitrary; we show $f\nin\p{\G}$. By
    the preceding lemma there exist $\alpha_1,\ldots,\alpha_{n_f}$
    constant or injective such that
    $f(\alpha_1,\ldots,\alpha\nf)\nin\A$. Choose a large and co-large $A\subseteq
    X$ such that $f(\alpha_1,\ldots,\alpha\nf)\inv[x]\cap A$  is large for
    a large set of
    $x\in X$. We modify the $\alpha_i$ to $\gamma_i\in\G$ in
    such a way that $\alpha_i\rest_A=\gamma_i\rest_A$ for $1\leq i\leq n_f$: If $\alpha_i$ is injective, then we
    can choose $\gamma_i$ to be a bijection. If $\alpha_i$ is
    constant, then $\gamma_i$ is delivered by the observation we just made. Thus, as
    $f(\alpha_1,\ldots,\alpha\nf)\rest_A=f(\gamma_1,\ldots,\gamma\nf)\rest_A$ we have
    $f(\gamma_1,\ldots,\gamma\nf)\nin\A\supseteq\G$.
\end{proof}
\begin{prop}\label{PROP:PolGUnderPolB}
    Let $\G\subseteq \B$ be a submonoid of $\Oo$ which contains
    all permutations. Then $\p{\G}\subseteq\p{\B}$.
\end{prop}
\begin{proof}
    Let $f\nin\p{\B}$ be arbitrary. We show $f\nin\p{\G}$. There
    are $\beta_1,\ldots,\beta\nf\in\B$ such that there exists $c\in
    X$ with the property that $f(\beta_1,\ldots,\beta\nf)\inv[c]$ is
    large. Define $\Gamma=(\beta_1,\ldots,\beta\nf)[X]$. Then since $\beta_i\in\B$,
    $H^i_a\cap\Gamma$ is small for all $1\leq i\leq n_f$ and all
    $a\in X$, where $H^i_a=\{x\in X^{n_f}:x_i=a\}$.
    Moreover, $f\inv[c]\cap\Gamma$ is large. Just like
    at the end of the proof of Lemma
    \ref{LEM:dasGammaLemma}, we can construct injective
    $\alpha_1,\ldots,\alpha\nf$ such that
    $f(\alpha_1,\ldots,\alpha\nf)$ is constant with value $c$. Choose
    $A\subseteq X$ large and co-large and bijections
    $\gamma_1,\ldots,\gamma\nf$ such that $\gamma_i\rest_A=\alpha_i\rest_A$ for $1\leq i\leq n_f$.
    Then, being constant on $A$,
    $f(\gamma_1,\ldots,\gamma\nf)\nin\B\supseteq\G$. Thus,
    $f\nin\p{\G}$.
\end{proof}
\subsubsection{Generous functions}
We now turn to monoids $\G\supseteq\S$ which are not submonoids of
$\A$. Our first goal is Proposition \ref{PROP:someIlambdaIsBelow},
in which we give a positive description of such monoids.

\begin{defn}
    A function $f\in\Oo$ is called \textit{generous} iff
    $f\inv[y]$ is either large or empty for  all $y\in X$.
\end{defn}
\begin{nota}
    Let $0\leq\lambda\leq \kappa$ be a cardinal. We denote by
    $\I_\lambda$ the set of all generous functions $f$ with the
    property that $|X\setminus f[X]|=\lambda$.
\end{nota}
\begin{lem}\label{LEM:basicIlambdaFacts}
    \begin{enumerate}
        \item{If $g\in \Oo$ is generous, then $f\circ g$ is generous for all $f\in\Oo$.}
        \item{$\I_\lambda$ is a subsemigroup and $\I_\lambda\cup\S$ a submonoid of $\Oo$ for all $\lambda\leq\kappa$.}
        \item{If $\lambda<\kappa$ and $f,g\in \I_\lambda$, then there exist $\alpha,\beta\in\S$ such that $f=\alpha\circ g\circ\beta$.}
        \item{$\I_\kappa$ contains all generous functions with small range, in particular the constant functions.}
        \item{If $g\in \I_\kappa$ has large range, then $\cl{\S\cup\{g\}}\supseteq \I_\kappa$.}
    \end{enumerate}
\end{lem}
\begin{proof}
    (1) and (4) are obvious. For (2), let $f,g\in \I_\lambda$; we
    want to show that $f\circ g\in \I_\lambda$. By (1)
    $f\circ g$ is generous, so it remains to show that $|X\sm f\circ g[X]|=\lambda$. We distinguish two
    cases: If $\lambda=\kappa$, then we
    have $\kappa \geq |X\setminus f\circ g[X]|\geq |X\setminus f[X]|=\kappa$ and so we are finished.
    Otherwise we
    claim $f[X]= f\circ g[X]$. Indeed, if $y\in f[X]$, then
    $f\inv[y]$ is large so that $g[X]\cap f\inv[y]\neq\emptyset$
    as $\lambda<\kappa$. Hence, $y\in f\circ g[X]$. Thus $f[X]\subseteq f\circ g[X]$ and the other inclusion is obvious.\\
    We prove (3). Write $f[X]=\{c_i\}_{i\in\kappa}$ and
    $g[X]=\{d_i\}_{i\in\kappa}$. Set further $C_i=f\inv[c_i]$ and
    $D_i=g\inv[d_i]$ and let $\beta_i$ be bijections from $C_i$
    onto $D_i$, $i\in\kappa$. Then $\beta=\bigcup_{i\in\kappa}\beta_i$
    is a bijection on $X$. Define the function $\alpha$ by
    $\alpha(d_i)=c_i$ for all $i<\kappa$ and extend $\alpha$ to $X$
    by an arbitrary bijection from $X\setminus g[X]$ onto $X\setminus
    f[X]$. It is readily verified that $f=\alpha\circ g\circ
    \beta$.\\
    To prove (5), let an arbitrary $f\in \I_\kappa$ be given; we
    show $f\in \cl{\S\cup\{g\}}$.
    Select any bijection $\gamma$ with the property that
    $g\circ\gamma(x)=g\circ\gamma(y)$ implies $f(x)=f(y)$ for all $x,y\in X$. This is
    possible since $g$ has large range and since both $f$ and $g$ are generous. Choose another bijection
    $\beta$ such that for all $x,y\in X$ we have $g\circ\beta\circ g\circ\gamma(x)=g\circ\beta\circ
    g\circ\gamma(y)$ iff $f(x)=f(y)$. Then it is clear that there
    is a bijection $\alpha$ satisfying $f(x)=\alpha\circ g\circ\beta\circ
    g\circ\gamma(x)$ for all $x\in X$.
\end{proof}
\begin{lem}\label{LEM:generousInfiniteRange}
    If $g\nin\A$, then there exists $\alpha\in\S$ such that the
    function $g\circ\alpha\circ g$ is generous and has large
    range.
\end{lem}
\begin{proof}
    There exists a large set $A\subseteq X$ such that $g\inv[a]$ is large for
    all $a\in A$. Set $E=g[X]\setminus A$ and $D=X\setminus g[X]$.
    Choose $B\subseteq A$ with the property that $A\setminus B$ is
    large and that $|B|=|g\inv [E]|$. Fix $a_0\in A\setminus B$.
    Take any function $\gamma: D\To g\inv[A]$ making
    $g\circ\gamma$ injective. $\gamma$ exists as $A$ is large. We
    want $\alpha\in \Oo$ to satisfy the following properties:
    \begin{itemize}
    \item[(i)]{$E$ shall be mapped injectively on a co-large part of $g\inv [a_0]$.}
    \item[(ii)]{$B$ shall be mapped bijectively onto $g\inv[E]$}
    \item[(iii)]{$\alpha\rest_D=\gamma$}
    \end{itemize}
    Since $E,B$ and $D$ are disjoint, we can indeed choose an
    injective partial function $\widetilde{\alpha}$ defined on
    $E\cup B\cup D$ which satisfies (i)-(iii). Because $X\setminus
    (E\cup B\cup D)\supseteq A\setminus B$ the domain of
    $\widetilde{\alpha}$ is co-large. Its range is also co-large
    as at least a large subset of $g\inv[a_o]$ is not in the
    range. Hence we can extend $\widetilde{\alpha}$ to
    $\alpha\in\S$. We claim that $\alpha$ has the asserted
    properties. Clearly $g\circ\alpha\circ g[X]\subseteq g[X]=A\cup E$; we show that $(g\circ\alpha\circ g)\inv[y]$ is large
    for all $y\in A\cup E$. Indeed, if $y\in E$, then $(g\circ\alpha)\inv[y]\subseteq B\subseteq
    A$. Thus, the preimage of $y$ under $g\circ\alpha\circ g$ is
    large. If $y\in A$, then $g\inv [y]$ is large and so is $g\inv
    [y]\setminus\gamma[D]$. Thus, $(g\circ\alpha)\inv [y]\setminus
    D$ is large as well. Hence, $(g\circ\alpha\circ g)\inv [y]$ is
    large which we wanted to show.
\end{proof}
\begin{prop}\label{PROP:someIlambdaIsBelow}
    Let $\G\subseteq\Oo$ be a monoid containing all bijections.
    Then either $\G\subseteq\A$ or there exists a cardinal
    $\lambda\leq\kappa$ such that $\I_\lambda\subseteq\G$.
\end{prop}
\begin{proof}
    This is an immediate consequence of Lemmas
    \ref{LEM:basicIlambdaFacts} and
    \ref{LEM:generousInfiniteRange}.
\end{proof}
The preceding proposition implies that when considering submonoids
$\G$ of $\Oo$ which contain the permutations, we can from now on
assume that $\I_\lambda\subseteq\G$ for some $\lambda$, since we
already treated the case $\G\subseteq\A$. We distinguish three
cases corresponding to the minimal $\lambda$ with the property
that $\I_\lambda\subseteq\G$: $\lambda=\kappa$, $0<\lambda<\kappa$
and $\lambda=0$.

\subsubsection{The case $\lambda=\kappa$}

Recall that $\G_\kappa$ consists of all functions $f\in\Oo$ with
the property that whenever $A$ is a large set then $f[X\sm A]$ is
co-large. Remember also that this is equivalent to $f$ being
either not almost surjective or almost injective.
\begin{prop}\label{PROP:PolGUnderPolGKappa}
    \begin{enumerate}
    \item{$\pol(\G_\kappa)$ is a maximal clone.}
    \item{If $\G$ is a submonoid of $\Oo$ such that
    $\S\cup \I_\kappa\subseteq\G$ and such that
    $\I_\lambda\subseteq\G$ for no $\lambda<\kappa$, then
    $\p{\G}\subseteq\p{\G_\kappa}$.}
    \end{enumerate}
\end{prop}
\begin{proof}
    (1) We will prove this together with the maximality of the other $\p{\G_\lambda}$ in Proposition \ref{PROP:PolGlambdaIsMax}.\\
    (2) Assume $f\nin\p{\G_\kappa}$; we show $f\nin\p{\G}$. Take
    $\alpha_1,\ldots,\alpha\nf\in\G_\kappa$ such that
    $g=f(\alpha_1,\ldots,\alpha\nf)\nin\G_\kappa$, that is, $g$ is almost
    surjective but not almost injective. Choose a co-large set
    $A\subseteq X$ such that $g[A]=g[X]$ is co-small. Because
    $\alpha_i\in\G_\kappa$,
    $\alpha_i[A]$ is co-large for $1\leq i\leq n_f$. Now fix for all
    $a\in A$ a large set $B_a$ such that $A\cap B_a=\{a\}$,
    $B_a\cap B_{a'}=\emptyset$ whenever $a\neq a'$ and such that
    $X=\bigcup_{a\in A} B_a$. This is possible since $A$ is
    co-large. Define for $1\leq i\leq n_f$ functions $\beta_i\in\Oo$
    by $\beta_i(x)=\alpha_i(a)$ whenever $x\in B_a$. It is clear
    that all $\beta_i$ are generous. Also, since
    $\beta_i[X]=\alpha_i[A]$ is co-large we have $\beta_i\in
    \I_\kappa\subseteq\G$ for all $1\leq i\leq n_f$. The function
    $f(\beta_1,\ldots,\beta\nf)$ is generous since it is constant on
    every $B_a$. Now $f(\beta_1,\ldots,\beta\nf)[X]\supseteq
    f(\beta_1,\ldots,\beta\nf)[A]=f(\alpha_1,\ldots,\alpha\nf)[A]=f(\alpha_1,\ldots,\alpha\nf)[X]$
    is co-small. Hence there exists a $\lambda<\kappa$ such that
    $f(\beta_1,\ldots,\beta\nf)\in \I_\lambda$, and since $\I_\lambda\nsubseteq
    \G$, we infer $f(\beta_1,\ldots,\beta\nf)\nin\G$ from Lemma \ref{LEM:basicIlambdaFacts} (3) which proves
    $f\nin\p{\G}$.
\end{proof}

\subsubsection{The case $0<\lambda<\kappa$}

We shall now investigate the case where $\G\nsupseteq\I_0$ but
$\G$ contains $\I_\lambda$ for some $0<\lambda<\kappa$. We collect
a couple of facts about the $\G_\lambda$ first. Recall that
$\G_\lambda$ consists of those functions $f$ for which it is true
that $|X\sm f[X\sm A]|\geq\lambda$ whenever $A\subseteq X$ is of
size $\lambda$. Recall also that for $\lambda=1$ or infinite this
is the case iff $f$ is $\lambda$-injective or not
$\lambda$-surjective.

\begin{lem}\label{LEM:GlambdaFacts} The following statements hold for all
$1\leq\lambda\leq\kappa$.
    \begin{enumerate}
    \item{If $g\in\On$ and $|X\sm g[X^n]|\geq\lambda$, then
    $g\in\p{\G_\lambda}$.}
    \item{$\G_\lambda$ is a submonoid of $\Oo$.}
    \item{$\G_n\supsetneqq\G_{n+1}$ for all $1\leq n< \aleph_0$.}
    \item{For $\lambda=1$ and for $\lambda\geq\aleph_0$,
    $\G_\lambda$ is a maximal submonoid of $\Oo$.}
    \end{enumerate}
\end{lem}
\begin{proof}
    (1) is obvious. For (2), let $f,g\in\G_\lambda$ and take an
    arbitrary $A\subseteq X$ with $|A|=\lambda$. Then $|X\sm g[X\sm
    A]|\geq\lambda$. Hence, $|X\sm
    f\circ g[X\sm A]|=|X\sm f[X\sm (X\sm g[X\sm A])]|\geq\lambda$ so that $f\circ g\in\G_\lambda$. It
    is clear that the identity map is an element of $\G_\lambda$ since
    it is injective.\\
    We prove (3). Observe first that the inclusion $X\sm f[X\sm(A\cup
    B)]\subseteq(X\sm f[X\sm A])\cup f[B]$ holds for all $A,B\subseteq X$ and all $f\in\Oo$. Now let
    $f\in \G_{n+1}$ for some $1\leq n<\aleph_0$. Take an arbitrary
    $n$-element subset $A$ of $X$. Choose any $a\nin A$. Then $n+1\leq
    |X\sm f[X\sm(A\cup \{a\})]|\leq|(X\sm f[X\sm A])|+ |f[\{a\}]|$ and
    so $n\leq |(X\sm f[X\sm A])|$. This proves $f\in\G_n$. It is
    obvious that $\G_n\neq\G_{n+1}$.\\
    The proof of (4) can be found in \cite{Ros74} (Proposition 5.2).
\end{proof}
\begin{lem}\label{LEM:existsLambda0}
    Let $1\leq\lambda\leq\kappa$. If $h\nin\G_\lambda$, then there exist a $\lambda_0 <\lambda$
    such that $\cl{\I_\lambda\cup\S\cup\{h\}}\supseteq \I_{\lambda_0}$. In particular,
    $\cl{\G_\lambda\cup\{h\}}\supseteq \I_{\lambda_0}$.
\end{lem}
\begin{proof}
    There exists $A\subseteq X$, $|A|=\lambda$ such that $|X\sm
    h[X\sm A]|<\lambda$. Set $\lambda_0=|X\sm
    h[X\sm A]|$. Choose a generous function $g$ with $g[X]=X\sm
    A$. Then $g\in\I_\lambda$ since $|X\sm g[X]|=|A|=\lambda$;
    thus, $h\circ g\in\cl{\I_\lambda\cup\{h\}}$. On the other
    hand, $h\circ g\in \I_{\lambda_0}$ and hence $\cl{\I_\lambda\cup\S\cup\{h\}}\supseteq
    \I_{\lambda_0}$ by Lemma \ref{LEM:basicIlambdaFacts} (3). The
    second statement is a direct consequence of the inclusion
    $\G_\lambda\supseteq\I_\lambda\cup\S$.
\end{proof}
\begin{lem}\label{LEM:theStrangeG}
    Let $B\subseteq X$, $|B|=\lambda_0<\lambda\leq\kappa$, and let
    $g\in\O\ut$ such that $g$ maps $(X\sm B)^2$ bijectively onto
    $X$ and such that $|g[B\mult X]\cup g[X\mult B]|<\kappa$. Then
    $g\in\p{\G_\lambda}$.
\end{lem}
\begin{proof}
    Let $\alpha,\beta\in\G_\lambda$ be given, and take an
    arbitrary $A\subseteq X$ of size $\lambda$. We have to show
    $|X\sm g(\alpha,\beta)[X\sm A]|\geq\lambda$. For
    $C=X\sm\alpha[X\sm A]$ we have $|C|\geq\lambda$. Thus, there
    exists some $c\in C\sm B$. Obviously, $g(\alpha,\beta)[\c A]\subseteq
    g[(\c \{c\})\mult X]$. But the conditions on $g$ yield that $g[(\c \{c\})\mult X]$ and
    $g[\{c\}\mult (\c B)]\sm (g[X\mult B]\cup g[B\mult X])$ are
    disjoint. Since $|g[\{c\}\mult (\c B)]|=\kappa$ and $|g[X\mult B]\cup g[B\mult
    X]|<\kappa$, this implies that $g(\alpha,\beta)$ misses
    $\kappa$ values on $X\sm A$ and hence,
    $g(\alpha,\beta)\in\G_\lambda$ and $g\in\p{\G_\lambda}$.
\end{proof}
\begin{prop}\label{PROP:PolGlambdaIsMax}
    \begin{enumerate}
    \item{$\p{\G_\lambda}$ is a maximal clone for all
    $1\leq\lambda\leq\kappa$.}
    \item{Let $\G\subseteq\Oo$ be a monoid containing all bijections as well as some $\I_\lambda$,
    where $0\leq \lambda<\kappa$, and let $\lambda$ be minimal
    with this property. If $\lambda>0$, then
    $\p{\G}\subseteq\p{\G_\lambda}$.}
    \end{enumerate}
\end{prop}
\begin{proof}
    (1) We show $\cl{\p{\G_\lambda}\cup\{h\}}=\O$ for an arbitrary $h\in\Oo\sm\G_\lambda$.
    By Lemma \ref{LEM:existsLambda0}, there
    exists $\lambda_0<\lambda$ such that
    $\I_{\lambda_0}\subseteq\cl{\G_\lambda\cup\{h\}}$. Now choose $B$
    and $g\in\p{\G_\lambda}$ as in Lemma \ref{LEM:theStrangeG}.
    Consider $\alpha: X\To (X\sm B)^2$ such that $\alpha$ takes every value twice.
    Clearly, $\alpha_1=\pi^2_1\circ\alpha$ and $\alpha_2=\pi^2_2\circ\alpha$ are elements of $\I_{\lambda_0}$.
    The function $p=g(\alpha_1,\alpha_2)=g\circ\alpha$ maps $X$ onto $X$ and takes
    every value twice as well. Therefore we can find a co-large
    set $A$ such that $p[A]=X$. Now fix a mapping $q: X\To A$
    so that $p\circ q$ is the identity map on $X$. Let an
    arbitrary $f\in\O$ be given. Then $q\circ f[X^{n_f}]\subseteq
    A$ is co-large which immediately implies $q\circ
    f\in\p{\G_\lambda}$. But then $f=p\circ(q\circ f)=f\in\cl{\p{\G_\lambda}\cup\{h\}}$ and
    so $\cl{\p{\G_\lambda}\cup\{h\}}=\O$ as $f$ was arbitrary.\\
    (2) First we claim that $\G\subseteq\G_\lambda$. Indeed, assume there exists $h\in
    \G\sm\G_\lambda$. Then, as $\I_\lambda\cup\S\subseteq\G$, by Lemma \ref{LEM:existsLambda0}
    there exists $\lambda_0<\lambda$ such that
    $\I_{\lambda_0}\subseteq\G$, in contradiction to the minimality of
    $\lambda$.\\
    Now let $f\nin\p{\G_\lambda}$ be arbitrary; we prove
    $f\nin\p{\G}$. There exist
    $\alpha_1,\ldots,\alpha\nf\in\G_\lambda$ such that
    $f(\alpha_1,\ldots,\alpha\nf)\nin\G_\lambda$. That is, there
    exists $A\subseteq X$ of size $\lambda$ with the property that
    $|\c f[\Gamma]|<\lambda$, where
    $\Gamma=\{(\alpha_1(x),\ldots,\alpha\nf (x)):x\in\c A\}$. Since
    $\alpha_i\in\G_\lambda$, $1\leq i\leq n_f$, for each $i$ there
    exists a set $B_i\subseteq X$, $|B_i|=\lambda$, such that
    $\alpha_i[\c A]\cap B_i=\emptyset$. Then $\Gamma\subseteq
    \Delta=(\c B_1)\mult \ldots \mult (\c B\nf)$. Choose $\beta: X\To \Delta$ onto and generous.
    Clearly $\beta_i=\pi^{n_f}_i\circ\beta\in\I_\lambda\subseteq\G$
    for all $1\leq i\leq n_f$. Now for all $C\subseteq X$ of size
    $\lambda<\kappa$ we have that $f(\beta_1,\ldots,\beta\nf)[\c
    C]=f[\Delta]\supseteq f[\Gamma]$ and so, as $|X\sm f[\Delta]|\leq|X\sm f[\Gamma]|<\lambda$,
    $f(\beta_1,\ldots,\beta\nf)\nin\G_\lambda\supseteq\G$. Hence,
    $f\nin\p{\G}$.
\end{proof}

\subsubsection{The case $\lambda=0$ and $\G\subseteq\F$}
In the following proposition we treat the case where
$\I_0\subseteq\G\subseteq\E\subseteq\F$. Recall that $\E$ consists
of those functions which are almost surjective (that is,
$\kappa$-surjective).
\begin{prop}\label{PROP:GunderE}
    \begin{enumerate}
    \item{$\p{\E}$ is a maximal clone.}
    \item{If $\G\subseteq\Oo$ is a monoid containing all bijections as well as
    $\I_0$, and if $\G\subseteq\E$, then $\p{\G}\subseteq\p{\E}$.}
    \end{enumerate}
\end{prop}
\begin{proof}
    (1) We prove that for any unary $h\nin\E$ we have
    $\cl{\p{\E}\cup\{h\}}=\O$. By definition $h[X]$ is
    co-large, so we can fix $A\subseteq X$ large and co-large such that
    $A\cap h[X]=\emptyset$. Choose any $g\in\Oo$ which maps $A$
    onto $X$ and which is constantly $0\in X$ on $\c A$. Then
    $g\in\E$ as it is onto. Moreover, $g\circ h$ is constantly
    $0$. Now let an arbitrary $f\in\On$ be given and define a
    function $\tilde{f}\in\O^{n+1}$ by
    $$
        \tilde{f}(x_1,\ldots,x_n,y)=\begin{cases}f(x_1,\ldots,x_n)&,y=0\\y&,\ow\end{cases}
    $$
    Then $\tilde{f}\in\p{\E}$. Indeed, this follows from the
    inclusion
    $\tilde{f}(\alpha_1,\ldots,\alpha_n,\beta)[X]\supseteq\beta[X]\sm\{0\}$ for arbitrary
    $\alpha_1,\ldots,\alpha_n,\beta\in\Oo$. Now
    $f(x)=\tilde{f}(x,0)=\tilde{f}(x,g\circ h(x_1))$ for all
    $x\in X^n$ and so $f\in \cl{\p{\E}\cup\{h\}}$.\\
    (2) Taking an arbitrary $f\nin\p{\E}$ we show that
    $f\nin\p{\G}$. There exist $\alpha_1,\ldots,\alpha\nf$ almost
    surjective such that $f(\alpha_1,\ldots,\alpha\nf)$ is not almost
    surjective. Consider a small set $A\subseteq X$ so that
    $A\cup\alpha_i[X]=X$ for all $1\leq i\leq n_f$. Let $\gamma$
    be a surjection from $\c A$ onto $X$ and define for $1\leq i\leq
    n_f$ functions
    $$
        \beta_i(x)=\begin{cases}\alpha_i\circ\gamma(x)&,x\in\c
        A\\x&,x\in A\end{cases}
    $$
    Clearly, all $\beta_i$ are surjective and
    $f(\beta_1,\ldots,\beta\nf)[X]=f(\alpha_1,\ldots,\alpha\nf)[X]\cup\{f(x,\ldots,x):x\in
    A\}$ is co-large. Fix any $\delta\in\I_0$. Obviously
    $\beta_i\circ\delta\in\I_0\subseteq\G$ and also
    $f(\beta_1\circ\delta,\ldots,\beta\nf\circ\delta)[X]$ is
    co-large. Thus
    $f(\beta_1\circ\delta,\ldots,\beta\nf\circ\delta)\nin\E\supseteq\G$
    so that we infer $f\nin\p{\G}$.
\end{proof}

In a next step we see what happens in the case
$\I_0\subseteq\G\subseteq\F$ and $\G\nsubseteq\E$. $\F$ is the set
of those functions which are almost surjective or constant.

\begin{prop}\label{PROP:GunderF}
    \begin{enumerate}
    \item{$\p{\F}$ is a maximal clone.}
    \item{If $\G\subseteq\F$ is a monoid which contains $\I_0$ as well as all
    bijections, and if $\G\nsubseteq \E$, then
    $\p{\G}\subseteq\p{\F}$.}
    \end{enumerate}
\end{prop}

\begin{proof}
    (1) can be found in \cite{Ros74} (Proposition 3.1).\\
    For (2), let $f\nin\p{\F}$ and fix
    $\alpha_1,\ldots,\alpha\nf\in\F$ satisfying
    $f(\alpha_1,\ldots,\alpha\nf)\nin\F$. Since $\G\nsubseteq\E$ but $\G\subseteq\F$, $\G$ must contain a constant function,
    and hence all constant functions as $\S\subseteq\G$. For those of the $\alpha_i$
    which are not constant we construct $\beta_i$ as in
    the proof of the preceding proposition, and for the constant
    ones we set $\beta_i=\alpha_i$. Observe that it is impossible that all $\alpha_i$ are
    constant. Choosing any $\delta\in\I_0$
    we obtain that for all $1\leq i\leq n_f$, $\beta_i\circ\delta$ is either constant
    or an element of $\I_0$, and hence in either case an element
    of $\G$. But as in the preceding proof,
    $f(\beta_1\circ\delta,\ldots,\beta\nf\circ\delta)\nin\F\supseteq\G$
    so that $f\nin\p{\G}$.
\end{proof}

\subsubsection{The case $\lambda=0$ and $\G\nsubseteq\F$} To conclude, we
consider submonoids $\G$ of $\Oo$ which contain the bijections as
well as $\I_0$, but which are not submonoids of $\F$. It turns out
that the polymorphism clones of such monoids are never maximal. We
start with a simple fact about such monoids.

\begin{lem}
    Let $\G\subseteq\Oo$ be a monoid containing $\S\cup\I_0$ such that
    $\G\nsubseteq\F$. Then $\chi=\{\rho\in\Oo:|\rho[X]|=2\text{ and }\rho \text{ is generous}\}\subseteq\G$.
\end{lem}
\begin{proof}
    Let $f\in\G\sm\F$. Since $f$ is not constant there exist
    $a\neq b$ in the range of $f$. Let $s:\c f[X]\To X$ be onto
    and generous and define $g\in\Oo$ by
    $$
        g(x)=\begin{cases}s(x)&,x\nin f[X]\\
                          a&,x=a\\
                          b&,\ow
        \end{cases}
    $$
    Then $g\in\I_0\subseteq\G$ and so $g\circ f\circ g\in \G$. On the other
    hand, $g\circ f\circ g\in\chi$ which proves the lemma since obviously any function of
    $\chi$ together with the permutations
    generate all of $\chi$.

\end{proof}

To prove that the remaining monoids do not yield maximal clones
via $\pol$, we are going to generalize the following completeness
criterion due to G. Gavrilov \cite{Gav65} (Lemma 31 on page 51)
for countable base sets.

\begin{lem}[G. Gavrilov]\label{LEM:deusExMachina}
    Let $X$ be countably infinite. If $\G\subseteq\Oo$ is a monoid containing $\S\cup\I_0\cup\chi$, and
    if $\H\subseteq\O$ is a set of functions such that
    $\cl{\Oo\cup\H}=\O$, then $\cl{\G\cup\H}=\O$.
\end{lem}

So we claim

\begin{prop}\label{PROP:generalizedGavrilov}
    Lemma \ref{LEM:deusExMachina} holds on all base sets of
    infinite regular cardinality.
\end{prop}

It follows immediately that $\p{\G}$ is not maximal for the
remaining monoids $\G$.

\begin{prop}\label{PROP:remainingPolsNotMaximal}
    If $\G\subseteq\Oo$ is a monoid such that $\S\cup\I_0\subseteq\G$ and
    such that $\G\nsubseteq\F$, then $\p{\G}$ is not maximal.
\end{prop}

\begin{proof}
    We have just seen that $\chi\subseteq\G$ so we can apply Proposition
    \ref{PROP:generalizedGavrilov}. Suppose towards contradiction that
    $\p{\G}$ is maximal. Since $\p{\G}\uo=\G\subsetneqq\Oo$ we
    have $\cl{\Oo\cup\p{\G}}=\O$. But then setting $\H=\p{\G}$ in
    the lemma yields that $\cl{\G\cup\p{\G}}=\O$, which is
    impossible as $\cl{\G\cup\p{\G}}=\p{\G}\neq\O$, contradiction.
\end{proof}

\subsection{The proof of Proposition
\ref{PROP:generalizedGavrilov}.} \label{SEC:almostUnary}

\begin{nota}
    We set $\L=\cl{\chi\cup\I_0\cup\S}$. Moreover, we write
    $\Const$ for the set of all constant functions.
\end{nota}
The following description of $\L$ is readily verified.
\begin{lem}
    $\L=\Const\cup\chi\cup\I_0\cup\S$. In words, $\L$ consists
    exactly of the bijections as well as of all generous functions
    which are either onto or take at most two values.
\end{lem}
\begin{lem}\label{LEM:universalFunction}
    Let $u\in\Oo$ be injective and not almost surjective. Then
    $\cl{\{u\}\cup\I_0} \supseteq \Oo$. In particular,
    $\cl{\{u\}\cup\L}\supseteq
    \Oo$.
\end{lem}
\begin{proof}
    Let an arbitrary $f\in\Oo$ be given. Take any $s: X\sm u[X] \To X$
    which is generous and onto. Now define $g\in \Oo$ by
    $$
        g(x)=\begin{cases}f(u\inv(x))&,x\in
        u[X]\\s(x)&,\ow\end{cases}
    $$
    Since $g\rest_{X\sm u[X]}=s$ we have $g\in\I_0$. Clearly, $f=g\circ
    u\in\cl{\{u\}\cup\I_0}$.
\end{proof}

\begin{defn} A function $f(x_1,\ldots,x_n)\in\On$ is \emph{almost unary} iff
    there exist a function $F:X\To \Pow(X)$ and some $1\leq k\leq n$ such that $F(x)$ is small for all
    $x\in X$ and such that for all
    $(x_1,\ldots,x_n)\in X^n$ we have $f(x_1,\ldots,x_n)\in F(x_k)$. We denote the set of all almost unary functions by
    $\U$.
\end{defn}

It is easy to see that on a base set of regular cardinality, $\U$
is a clone which contains $\Oo$. See \cite{Pin032} for a list of
all clones above $\U$; there are countably many, so in particular
$\U$ is not maximal. The reason for us to consider almost unary
functions is the following lemma.

\begin{lem}\label{LEM:LandNotAlmostUnary}
    Let $f\in\O\un\sm\U$ be any function which is not almost
    unary. Then $\cl{\{f\}\cup\L}\supseteq\Oo$.
\end{lem}

Observe that this lemma implies that $\p{\G}\subseteq\U$ for all
proper submonoids $\G$ of $\Oo$ which contain $\L$ and that we can
therefore conclude directly that these polymorphism clones are not
maximal. We will now prove Lemma \ref{LEM:LandNotAlmostUnary} by
showing that $\L$ together with a not almost unary $f$ generate a
function $u$ as in Lemma \ref{LEM:universalFunction}. We start by
observing that $\L$ and $f$ generate functions of arbitrary range.

\begin{lem}\label{LEM:largeAndColargeRange}
    Let $f\in\On\sm\U$. Then there exists a unary $g\in\cl{\{f\}\cup\L}$
    such that the
    range of $g$ is large and co-large.
\end{lem}
\begin{proof}
    We distinguish two cases.\\
    \textbf{Case 1.} For all $1\leq i\leq n$ and all $c\in X$ it
    is true that $f[X^{i-1}\mult\{c\}\mult X^{n-i}]$ is
    co-small. Then consider an arbitrary large and co-large  $A\subseteq
    X$. Set
    $\Gamma=f\inv[X\sm A]\subseteq X^n$ and let $\alpha: X\To\Gamma$ be onto. By the assumption for this case,
    $f[X^{i-1}\mult\{c\}\mult X^{n-i}]\sm A$ is still
    large for all $1\leq i\leq n$ and all $c\in X$. Thus
    the components $\alpha_i=\pi^n_i\circ\alpha$ are generous and
    onto; hence, $\alpha_i\in\I_0\subseteq\L$ for all $1\leq i\leq n$. But now
    $f(\alpha_1,\ldots,\alpha_n)[X]=f[X^n]\sm A$ is large and
    co-large so that it suffices to set $g=f\circ\alpha$.\\
    \textbf{Case 2.} There exists $1\leq i\leq n$ and $c\in X$
    such that $f[X^{i-1}\mult\{c\}\mult X^{n-i}]$ is co-large,
    say without loss of generality $i=1$. Since $f\nin\U$, there
    exists $d\in X$ satisfying that $f[\{d\}\mult X^{n-1}]$ is
    large. Choose $\Gamma\subseteq X^{n-1}$ large and co-large such that
    $f[\{d\}\mult\Gamma]$ is large and such that $f[\{c\}\mult
    X^{n-1}]\cup f[\{d\}\mult\Gamma]$ is still
    co-large. Take moreover $\alpha_2,\ldots,\alpha_n\in\I_0$ so that
    $(\alpha_2,\ldots,\alpha_n)[X]=X^{n-1}$. Now we define
    $\alpha_1\in\Oo$ by
    $$
        \alpha_1(x)=
        \begin{cases}d&,(\alpha_2,\ldots,\alpha_n)(x)\in\Gamma\\
        c&,\ow.
        \end{cases}
    $$
    Clearly, $\alpha_1\in\chi\subseteq\L$. Now it is enough to set $g=f(\alpha_1,\ldots,\alpha_n)$ and observe
    that $g[X]=f[\{c\}\mult (X^{n-1}\sm\Gamma)]\cup f[\{d\}\mult\Gamma]$
    is large and co-large.
\end{proof}
\begin{lem}\label{LEM:arbitraryRange}
    Let $f\in\On\sm\U$. Then for all $A\subseteq X$ there
    exists $h\in\cl{\{f\}\cup\L}$ with $h[X]=A$.
\end{lem}
\begin{proof}
    By Lemma \ref{LEM:largeAndColargeRange} there exists $g\in \cl{\{f\}\cup\L}$
    having a large and co-large range. Now taking any
    $\delta\in\I_0\subseteq\L$ with $\delta[g[X]]=A$ and setting
    $h=\delta\circ g$ proves the assertion.
\end{proof}
\begin{lem}\label{LEM:allGenerousFunctions}
    If $f\in\On\sm\U$, then $\cl{\{f\}\cup\L}$ contains all generous
    functions.
\end{lem}
\begin{proof}
    Let any generous $g\in\Oo$ be given and take with the help of the preceding lemma $h\in\cl{\{f\}\cup\L}$ with
    $h[X]=g[X]$. By setting $h'=h\circ\delta$, where
    $\delta\in\I_0\subseteq\L$ is arbitrary, we obtain a generous
    function with the same property. Now it is clear that there
    exists a bijection $\sigma\in\S\subseteq\L$ such that
    $g=h'\circ\sigma$.
\end{proof}

Now that we know that we have all generous functions we want to
make them injective. We start by reducing the class of functions
$f$ under consideration.

\begin{lem}\label{LEM:columnsWithManyDifferentValues}
    If $f\in\On\sm\U$ is so that for all $1\leq i\leq n$ and for all $a,b\in X$
    the set of all tuples
    $(x_1,\ldots,x_{i-1},x_{i+1},\ldots,x_n)\in X^{n-1}$ with
    $f(x_1,\ldots,x_{i-1},a,x_{i+1},\ldots,x_n)\neq
    f(x_1,\ldots,x_{i-1},b,x_{i+1},\ldots,x_n)$ is small, then
    $\cl{\{f\}\cup\L}\supseteq\Oo$.
\end{lem}
\begin{proof}
    Since
    $f\nin\U$ we can for every $1\leq i\leq n$ choose $c_i\in X$
    such that $f[X^{i-1}\mult\{c_i\}\mult X^{n-i}]$ is large.
    Choose moreover for every $1\leq i\leq n$ large sets $A_i\subseteq f[X^{i-1}\mult\{c_i\}\mult X^{n-i}]$ such that
    $\bigcup_{i=1}^n A_i$ is co-large and such that $A_i\cap
    A_j=\emptyset$ for $i\neq j$. Write each $A_i$ as a disjoint
    union of many large sets: $A_i=\bigcup_{x\in X} A^x_i$. Let $\triangleleft$ be any well-order of
    $X^n$ of type $\kappa$. Define
    $\Gamma\subseteq X^n$ by $x\in\Gamma$ iff there exists $1\leq i\leq
    n$ such that $f(x)\in A_i^{x_i}$ and whenever $y\triangleleft
    x$ and $y\in\Gamma$ then $f(x)\neq f(y)$. Observe that the latter
    condition ensures that $f\rest_\Gamma$ is injective.

    Now observe that for all $1\leq i\leq n$, all $c\in X$ and all large $B\subseteq A_i$
    we have that $f[X^{i-1}\mult\{c\}\mult X^{n-i}]\cap B$ is
    large. Indeed, say without loss of generality $i=1$ and set
    $D=\{(x_2,\ldots,x_n):f(c,x_2,\ldots,x_n)\neq f(c_1,x_2,\ldots,x_n)\}$.
    Then $D$ is small by our assumption. Now $|f[\{c\}\mult
    X^{n-1}]\cap B|\geq |f[\{c\}\mult
    (X^{n-1}\sm D)]\cap B| = |f[\{c_1\}\mult
    (X^{n-1}\sm D)]\cap B|=\kappa$. In particular, this observation
    is true for $B=A_i^c$. This implies that the set
    $\{x\in\Gamma: x_i=c\}$ is large for all $1\leq i\leq n$ and
    all $c\in X$. Moreover, $\Gamma$ itself is large.

    Therefore there exists a bijection $\alpha: X\To\Gamma$. By
    the preceding observation, the components
    $\alpha_i=\pi^n_i\circ\alpha$ are onto and generous, so
    $\alpha_i\in\I_0\subseteq\L$ for all $1\leq i\leq n$. Since $\alpha$
    is injective, $\alpha[X]=\Gamma$ and $f\rest_\Gamma$ is
    injective, we have that $g=f(\alpha_1,\ldots,\alpha_n)\in\cl{\{f\}\cup\L}$ is
    injective. Furthermore, $g[X]=f[\Gamma]\subseteq\bigcup_{i=1}^n A_i$
    is co-large. Whence
    $\Oo\subseteq\cl{\{g\}\cup\L}\subseteq\cl{\{f\}\cup\L}$ by Lemma
    \ref{LEM:universalFunction} and we are done.
\end{proof}

\begin{lem}\label{LEM:columnsWithManyDifferentValuesII}
    If $f\in\On\sm\U$ is so that for all
    $1\leq i\leq n$ there exist $c\in X$ and $S\subseteq
    X^n$ with $\pi^n_i[S]=\{c\}$ such that $f[S]$ large and such that for all $b\in X$ the set $\{x\in S: f(x)\neq
    f(x_1,\ldots,x_{i-1},b,x_{i+1},\ldots,x_n)\}$ is small, then $\cl{\{f\}\cup\L}\supseteq\Oo$.
\end{lem}
\begin{proof}
    Fix for every $1\leq i\leq n$ an element $c_i\in
    X$ and a set $S_i\subseteq X^n$ such that $\pi^n_i[S_i]=\{c_i\}$ and
    such that $f[S_i]$ large and such that for all $b\in X$ the
    set $\{x\in S_i: f(x)\neq
    f(x_1,\ldots,x_{i-1},b,x_{i+1},\ldots,x_n)\}$ is small. Set
    $A_i=f[S_i]$, $1\leq i\leq n$. By thinning out the $S_i$ we
    can assume that the $A_i$ are disjoint and that $\bigcup_{i=1}^n A_i$ is co-large. Now one follows the proof of
    the preceding lemma.
\end{proof}

\begin{lem}
    If $f\in\On\sm\U$, then there exists $g\in\cl{\{f\}\cup\L}$ having
    co-large range and with the property that $\{x\in
    X:|g\inv[x]|=1\}$ is large (that is, the kernel of $g$ has $\kappa$ one-element classes).
\end{lem}
\begin{proof}
    There is nothing to prove if $f$ satisfies the condition of Lemma
    \ref{LEM:columnsWithManyDifferentValuesII}, so assume it does not, and let $i=1$ witness this.
    Take $c\in X$ such that
    $f[\{c\}\mult X^{n-1}]$ is large and choose $S\subseteq X^n$ such
    that $\pi^n_1[S]=\{c\}$, such that $f[S]$ is still large and such that $f\rest_S$ is injective.
    By the lemma, there exists $b\in X$ such that $\{x\in S:
    f(x)\neq f(b,x_2,\ldots,x_n)\}$ is large. Thus, we can find a large $A\subseteq
    S$ with the property that $f[A]$ and $f[\{(b,x_2,\ldots,x_n): x\in
    A\}]$ are disjoint and such that the union of these two sets
    is co-large.
    Choose now generous $\alpha_2,\ldots,\alpha_n\in\Oo$ such that
    $(c,\alpha_2,\ldots,\alpha_n)[X]=A$. Since $\cl{\{f\}\cup\L}$ contains all generous functions by Lemma
    \ref{LEM:allGenerousFunctions}, we have $\alpha_j\in\L$ for $2\leq
    j\leq n$. Take a large and co-large $B\subseteq X$ such that
    $(c,\alpha_2,\ldots,\alpha_n)\rest_B$ is injective. Define
    $$
        \alpha_1(x)=
        \begin{cases}c&,x\in B\\
        b&,\ow
        \end{cases}
    $$
    and set $g=f(\alpha_1,\ldots,\alpha_n)$. Then $g\in\cl{\{f\}\cup\L}$ as $\alpha_1\in\chi\subseteq\L$.
    Clearly, $(\alpha_1,\ldots,\alpha_n)\rest_B$ is injective and so is
    $g\rest_B$. Since $g[B]$ and $g[X\sm B]$ are disjoint we have
    that $|g\inv[x]|=1$ for all $x\in g[B]$. Moreover, $g[X]\subseteq f[A]\cup f[\{(b,x_2,\ldots,x_n): x\in
    A\}]$ is co-large.
\end{proof}

\begin{lem}
    Let $f\in\On\sm \U$. If $h\in\Oo$ is a function whose
    kernel has at least one large equivalence class (that is, there
    exists $x\in X$ with $h\inv[x]$ large), then $h\in \cl{\{f\}\cup\L}$.
\end{lem}
\begin{proof}
    There exist a large $B\subseteq X$ and $b\in X$ such that
    $h[B]=\{b\}$. Let $g$ be provided by the preceding lemma.
    With the help of permutations of the base set we can assume that $|g\inv
    [x]|=1$ for all $x\in g[X\sm B]$. Since the range of $g$
    is co-large we can find $\delta: X\sm g[X] \To X$ onto and
    generous. Now define $m\in\Oo$ by
    $$
        m(x)=
        \begin{cases}
            \delta(x)&,x\nin g[X]\\
            b&, x\in g[B]\\
            h(g\inv(x))&, x\in g[X\sm B].
        \end{cases}
    $$
    Obviously $m\in \I_0\subseteq\L$ and $h=m\circ g\in\cl{\{f\}\cup\L}$.
\end{proof}

Having found many functions which $\cl{\{f\}\cup\L}$ must contain,
we are finally ready to prove Lemma \ref{LEM:LandNotAlmostUnary}.

\begin{proof}[Proof of Lemma \ref{LEM:LandNotAlmostUnary}]
    There are $c_1,\ldots,c_n\in X$ such
    that $f[X^{i-1}\mult\{c_i\}\mult X^{n-i}]$ is large for $1\leq
    i\leq n$. Take $B_1,\ldots,B_n$ large such that
    $\pi^n_i[B_i]=\{c_i\}$ for all $1\leq i\leq n$ and with the
    property that $f\rest_B$ is injective and $f[B]$ is co-large, where $B={\bigcup_{i=1}^n
    B_i}$. Let $\alpha: X\To B$ be any bijection. Since $\alpha_i\inv[c_i]$ is large for every
    component $\alpha_i=\pi^n_i\circ\alpha$, the preceding lemma yields $\alpha_i\in\cl{\{f\}\cup\L}$ for
    $1\leq i\leq n$. Whence, $g=f(\alpha_1,\ldots,\alpha_n)\in\cl{\{f\}\cup\L}$. But
    $g[X]=f[B]$ is co-large and $g$ is injective by construction;
    thus Lemma \ref{LEM:universalFunction} yields $\Oo\subseteq\cl{\{g\}\cup\L}\subseteq\cl{\{f\}\cup\L}$.
\end{proof}

This brings us back to our original goal.

\begin{proof}[Proof of Proposition \ref{PROP:generalizedGavrilov}]
    Since $\cl{\Oo\cup\H}=\O$, there must exist some
    $f\in\H\sm\U$. But then, since $\G\supseteq \L$, Lemma \ref{LEM:LandNotAlmostUnary} implies
    $\cl{\G\cup\H}\supseteq\Oo$ so that we
    infer $\cl{\G\cup\H}=\O$.
\end{proof}

\end{section}

\begin{section}{The proof of Theorem \ref{THM:allMaximalMonoids}}

We now determine on an infinite $X$ all maximal submonoids of
$\Oo$ which contain the permutations, proving Theorem
\ref{THM:allMaximalMonoids}. In a first section, we present the
part of the proof which works on all infinite sets; then follow
one section specifically for the case of a base set of regular
cardinality and another section for the singular case. Throughout
all parts we will mention explicitly whenever a statement is true
only on $X$ of regular or singular cardinality, respectively.

\subsection{The part which works for all infinite sets}
\begin{prop}\label{PROP:GlambdaIsMaximal}
    $\G_\lambda$ is a maximal submonoid of $\Oo$ for $\lambda=1$ and
    $\aleph_0\leq\lambda\leq\kappa$.
\end{prop}
\begin{proof}
    As already mentioned in Lemma \ref{LEM:GlambdaFacts}, the
    maximality of the $\G_\lambda$ for $\lambda=1$ or infinite has
    been proved in \cite{Ros74} (Proposition 5.2).
\end{proof}

The maximal monoids of Proposition \ref{PROP:GlambdaIsMaximal}
already appeared in the preceding section since they give rise to
maximal clones via $\pol$. We shall now expose maximal monoids
above the permutations which do not have this property. Recall
that $\M_\lambda$ consists of all functions which are either
$\lambda$-surjective or not $\lambda$-injective.

\begin{prop}\label{PROP:MlambdaIsMaximal}
    Let $\lambda=1$ or $\aleph_0\leq\lambda\leq\kappa$. Then $\M_\lambda$ is a maximal submonoid of $\Oo$.
\end{prop}
\begin{proof}
    We show first that $\M_\lambda$ is closed under composition.
    Let therefore $f,g\in\M_\lambda$, that is, those functions are
    either $\lambda$-surjective or not $\lambda$-injective; we claim that $f\circ g$ has either of these properties. It is
    clear that if $g$ is not $\lambda$-injective, then $f\circ g$
    has the same property. So let $g$ be $\lambda$-surjective. It is easy to see that if
    $f$ is $\lambda$-surjective, then so is $f\circ g$. So
    assume finally that $f$ is not $\lambda$-injective. We claim
    that $f\circ g$ is not $\lambda$-injective either. For
    $\lambda=1$ this is just the statement that if $f$ is not
    injective, and $g$ is surjective, then $f\circ g$ is not
    injective, which is obvious. Now consider the infinite case.
    There exist disjoint $A,B\subseteq X$ of size $\lambda$ such that
    $f[A]=f[B]$. Set $A'=A\cap g[X]$; $A'$ still has size
    $\lambda$ as $g$ misses less than $\lambda$ values. Clearly
    $B'=\{x\in B:\exists y\in A' (f(x)=f(y))\}$ has size $\lambda$
    as well and so does $B''=B'\cap g[X]$. But now for the sets
    $C=g\inv [A']$ and $D=g\inv [B'']$ it is true that $|C|,|D|
    \geq\lambda$, $C\cap D=\emptyset$, and $f\circ g[C]=f\circ
    g[D]$; hence $f\circ g$ is not $\lambda$-injective.\\
    Now we prove that $\M_\lambda$ is maximal in $\Oo$. Consider for
    this reason any $m\nin\M_\lambda$, that is, $m$ is
    $\lambda$-injective and misses at least $\lambda$ values. There
    exists $A\subseteq X$ so that $|X\sm A|<\lambda$ and such that the restriction of $m$ to $A$ is injective.
    Take any injection $i\in\Oo$ with $i[X]=A$. Then $i\in\M_\lambda$ as $i$ is $\lambda$-surjective. Now let
    $f\in\Oo$ be arbitrary. Define
    $$
        g(x)=
        \begin{cases}
            f((m\circ i)\inv(x))&,x\in m\circ i[X]\\
            a&, \ow\\
        \end{cases}
    $$
    where $a\in X$ is any fixed element of $X$. Being constant on the
    complement of the range of $m$, $g$ it is not $\lambda$-injective
    and whence an element of $\M_\lambda$. Therefore $f=g\circ
    m\circ i\in\cl{\M_\lambda\cup\{m\}}$ so that we infer
    $\cl{\M_\lambda\cup\{m\}}\supseteq\Oo$.
\end{proof}

\begin{lem}\label{LEM:allMonoidsAboveI0}
    There are no other maximal monoids above $\S\cup\I_0$ except the
    $\M_\lambda$ ($\lambda=1$ or $\aleph_0\leq \lambda\leq \kappa$).
\end{lem}

\begin{proof}
    Let $\G\supseteq \I_0\cup\S$ be a submonoid of $\Oo$ which is not contained in any of the
    $\M_\lambda$; we prove that $\G=\Oo$. To do this, we
    show that $\G$ contains an injective function $u\in
    \Oo$ with co-large range; then the lemma follows from Lemma
    \ref{LEM:universalFunction}. Fix for every $\lambda$ a
    function $m_\lambda\in\G\sm\M_\lambda$. Since $m_{\kappa}$ is
    $\kappa$-injective, there exists a cardinal $\lambda_1
    <\kappa$ and a set $A_1\subseteq X$ of size $\lambda_1$ such that
    the restriction of $m_\kappa$ to the complement of $A_1$ is
    injective. If $\lambda_1$ is infinite, then consider
    $m_{\lambda_1}$. Not being an element of $\M_{\lambda_1}$, $m_{\lambda_1}$ misses at least
    $\lambda_1$ values. Hence by adjusting it with a suitable
    permutation we can assume that $m_{\lambda_1}[X]\subseteq X\sm A_1$.
    There exists a cardinal $\lambda_2<\lambda_1$ and a subset
    $A_2$ of $X$ of size $\lambda_2$ such that the restriction of
    $m_{\lambda_1}$ to the complement of $A_2$ is injective.
    Hence, writing $\lambda_0=\kappa$ we obtain that
    $m_{\lambda_0}\circ m_{\lambda_1}\in\G$ is injective on $X\sm A_2$ and misses $\kappa$
    values. We can iterate this to arrive after a finite number of steps at a set
    $A_n$ of finite size $\lambda_n$ such that the restriction of
    $m_{\lambda_{0}}\circ \ldots\circ m_{\lambda_{n-1}}\in\G$ to $X\sm
    A_n$ is injective and misses $\kappa$ values. Since
    $m_1\nin\M_1$ is injective and misses at least one value we
    conclude that the iterate $m_1^{\lambda_n}\in\G$ is injective and
    misses at least $\lambda_n$ values. Modulo
    permutations we may assume that $m_1^{\lambda_n}[X]\subseteq X\sm
    A_n$. But now we have that $m_{\lambda_{0}}\circ \ldots\circ
    m_{\lambda_{n-1}}\circ m_1^{\lambda_n}\in \G$ is injective and misses $\kappa$ values,
    implying that $\G=\Oo$.
\end{proof}

\subsection{The case of a base set of regular cardinality}

We now finish the proof of Theorem \ref{THM:allMaximalMonoids} for
the case when $X$ has regular cardinality. The proof for this case
comprises Propositions \ref{PROP:GlambdaIsMaximal},
\ref{PROP:MlambdaIsMaximal}, \ref{PROP:AisMaximal} and
\ref{PROP:noOtherMaxMonoids}.

\begin{prop}\label{PROP:AisMaximal}
    If $X$ is of regular cardinality, then $\A$ is a maximal submonoid of $\Oo$.
\end{prop}

\begin{proof}
    This has been proved in \cite{Ros74} (Proposition 4.1).
\end{proof}

\begin{prop}\label{PROP:noOtherMaxMonoids}
    Let $X$ have regular cardinality. There exist no other maximal submonoids of $\Oo$ containing the permutations
    except those listed in Theorem \ref{THM:allMaximalMonoids} for the regular case.
\end{prop}
\begin{proof}
    Assume that $\G\supseteq \S$ is a submonoid of $\Oo$ not contained in any
    of the monoids of the theorem; we show that $\G=\Oo$. Indeed, since
    $\G\nsubseteq\A$, Proposition \ref{PROP:someIlambdaIsBelow} tells us that there exists a cardinal $\lambda
    \leq\kappa$ such that $\I_\lambda$ is contained in $\G$. Choose
    $\lambda$ minimal with this property. If $\lambda$ was greater than $0$, then
    $\G\subseteq\G_\lambda$ for otherwise Lemma \ref{LEM:existsLambda0}
    would yield a contradiction to the minimality of $\lambda$. But
    this is impossible as we assumed that $\G$ is not contained in
    any of the $\G_\lambda$, so we conclude that $\lambda=0$. Now
    Lemma \ref{LEM:allMonoidsAboveI0} implies that $\G=\Oo$.
\end{proof}

\subsection{The case of a base set of singular cardinality}

The only problem with base sets of singular cardinality is that
the set $\A$ is not closed under composition; in fact,
$\cl{\A}=\O$. A slight adjustment of the definition of $\A$ works
in this case. We will refer to results from preceding sections;
this might look unsafe since there we restricted ourselves to base
sets of regular cardinality. However, when proving the particular
results cited here we did not use the regularity of the base set.
The proof of Theorem \ref{THM:allMaximalMonoids} for singular
cardinals comprises Propositions \ref{PROP:GlambdaIsMaximal},
\ref{PROP:MlambdaIsMaximal}, \ref{PROP:A'isMaximal} and
\ref{PROP:noOtherMaximalMonoidsSingular}.

\begin{defn}
    A function $f\in\Oo$ is said to be \emph{harmless} iff
    there exists $\lambda < \kappa$ such that the set of all $x\in X$ for which
    $|f\inv [x]|>\lambda$ is small. With this definition,
    $\A'$ as defined in Theorem \ref{THM:allMaximalMonoids} is the set of all harmless functions.
\end{defn}

\begin{lem}
    $\A'$ is a monoid and $\A'\subseteq \A$. Moreover, $\A=\A'$ iff $\kappa$ is a
    successor cardinal.
\end{lem}
\begin{proof}
    It is obvious that $\A'\subseteq \A$ and that $\A=\A'$ iff $\kappa$ is a
    successor cardinal. To prove that $\A'$ is closed
    under composition, let $f,g\in\A'$; we show $h=f\circ g\in \A'$.
    There exist $\lambda_f,\lambda_g<\kappa$ witnessing that $f$
    and $g$ are harmless. Set $\lambda$ to be $\max(\lambda_f,\lambda_g)$; we
    claim that the set of $x\in X$ for which $|h\inv [x]|>\lambda$
    is small. For if $|h\inv [x]|>\lambda$, then either $|g\inv
    [x]|>\lambda$ or there exists $y\in g\inv
    [x]$ such that $|f\inv [y]|>\lambda$. Both possibilities occur
    only for a small number of $x\in X$ and so $h$ is harmless.
\end{proof}

\begin{lem}\label{LEM:notInA'MakesNotInA}
    Let $X$ have singular cardinality. If $g\nin\A'$, then $g$ together with $\S$ generate a function not in $\A$.
\end{lem}

\begin{proof}
    Set $\lambda < \kappa$ to be the cofinality of $\kappa$. Because $g$ is not harmless, there
    exist distinct sequences $(x^0_\xi)_{\xi < \lambda},\ldots,(x^\kappa_\xi)_{\xi < \lambda}$ of distinct elements of $X$
    such that $\bigcup_{\xi <\lambda} g\inv [x^\zeta_\xi]$
    is large for all $\zeta <\kappa$. Indeed, if $(\mu_\xi)_{\xi<\lambda}$ is any cofinal sequence of cardinalities
    in $\kappa$, then the fact that $g$ is not harmless allows us to pick for every $\xi<\lambda$ an element $x^0_\xi\in X$ such
    that $|g\inv[x^0_\xi]|>\mu_\xi$; it is also no problem to choose the elements distinct. This yields the first sequence and
    since with every sequence we are using up only $\lambda <\kappa$ elements, the definition of harmlessness ensures that we can
    repeat the process $\kappa$ times. By throwing away half of the sequences, we may assume that the set of all $y\in X$ which
    do not appear in any of the sequences is large.\\
    There exists a permutation $\alpha\in\S$ such that
    $g\circ\alpha(x_{\xi_1}^{\zeta_1})=g\circ\alpha (x_{\xi_2}^{\zeta_2})$ if and only if $\zeta_1=\zeta_2$,
    for all $\zeta_1,\zeta_2<\kappa$ and all $\xi_1,\xi_2
    <\lambda$. For we can map every sequence $(x_\xi^\zeta)_{\xi <\lambda}$ injectively into an equivalence
    class of the kernel of $g$ of size greater than $\lambda$;
    since there are many such classes every sequence can be assigned an
    own class, and we choose the classes so that a large number of classes are not hit at all.
    This partial injective mapping we can then extend to the permutation $\alpha$ as
    it is defined on a co-large set and has co-large range.\\
    Set $y^\zeta=g\circ\alpha(x_0^\zeta)$ for all $\zeta<\kappa$. Then
    the
    $y^\zeta$ are pairwise distinct and for all $\zeta<\kappa$ we have that
    $(g\circ\alpha\circ g)\inv[y^\zeta]\supseteq \bigcup_{\xi <\lambda} g\inv [x^\zeta_\xi]$
    is large. Hence, $g\circ \alpha \circ g\nin
    \A$.
\end{proof}

\begin{prop}\label{PROP:A'isMaximal}
    Let $X$ have singular cardinality. Then $\A'$ is a maximal submonoid of $\Oo$.
\end{prop}

\begin{proof}
    Let $g\in\Oo\sm\A'$. We know that $g$ together with $\A'$ generate a function not in $\A$.
    Then by Lemma \ref{LEM:generousInfiniteRange}, we obtain a function which is generous and has large range,
    call it $h$. Now take any $f\in\Oo$ such that $f\circ h[X]=X$
    which is injective on $h[X]$ and constant on $X\sm h[X]$. Then
    $f\in\A'$ and $f\circ h\in\I_0$. Thus,
    $\I_0\subseteq\cl{\{g\}\cup\A'}$ and since all injections are elements of $\A'$ we can apply
    Lemma \ref{LEM:universalFunction} to prove
    $\cl{\{g\}\cup\A'}\supseteq\Oo$.
\end{proof}

\begin{prop}\label{PROP:noOtherMaximalMonoidsSingular}
    Let $X$ have singular cardinality. There
    exist no other maximal submonoids of $\Oo$ containing the permutations
    except those listed in Theorem \ref{THM:allMaximalMonoids} for the singular case.
\end{prop}

\begin{proof}
    If $\G\supseteq\S$ is a submonoid of $\Oo$ which is not
    contained in $\A'$, then it is not contained in $\A$ by Lemma \ref{LEM:notInA'MakesNotInA}. From
    this point, one can follow
    the proof of Proposition \ref{PROP:noOtherMaxMonoids}.
\end{proof}

\end{section}
\end{chapter}\newpage
\pagestyle{myheadings}\markright{\uppercase{Bibliography}\hfill}

\newpage
\appendix
\pagestyle{myheadings}\markright{\uppercase{Curriculum
vitae}\hfill} \addcontentsline{toc}{chapter}{Curriculum vitae}
\chapter*{Curriculum Vitae}

\hspace*{4mm} {\it Personal Data}\vspace*{1mm}

\begin{tabular}{l}
Michael Pinsker\\
Born Nov. 4,~1977 in T\"{u}bingen, Germany.\\
Son of Wilhelm and Doris Pinsker, 2 brothers.\\
Austrian citizenship.\\
\end{tabular}\\

\vspace*{3mm} {\it School}\vspace*{1mm}

\begin{tabular}{p{40mm} p{92mm}}
09/1984 - 06/1987 & Primary School in Tü\"{u}bingen. \\
09/1987 - 06/1988 & Primary School in Vienna. \\
09/1988 - 06/1996 & Secondary School BRG Wenzgasse, Vienna.\\
\end{tabular}

\vspace*{3mm} {\it Civil service}\vspace*{1mm}

\begin{tabular}{p{40mm} p{92mm}}
10/1996 - 09/1997 & Civil service in Vienna\\
\end{tabular}

\vspace*{3mm} {\it Studies}\vspace*{1mm}

\begin{tabular}{p{40mm} p{92mm}}
10/1997 - 06/2002 & Study of Technical Mathematics at the
        Vienna University of Technology.
        Concentration in Set Theory and Universal Algebra.
        Diploma thesis ``Rosenberg's characterization of maximal clones'' written under the guidance of A.o. Prof. Martin Goldstern,
        Department of Algebra and Computer Science.\\
10/2002 - 09/2004 & Ph.D. student at the Vienna University of
        Technology under the supervision of A.o. Prof. Martin Goldstern,
        Subject: Clones on infinite sets.\\
\end{tabular}

\newpage {\it Research visits}\vspace*{1mm}

\begin{tabular}{p{40mm}p{92mm}}
02/2001 - 06/2001 & Technical University of Denmark within the
ERASMUS
    program.\\
03/2003 - 08/2003 & Free University and Humboldt University
Berlin, Germany.\\
09/2003 - 02/2004& Masaryk University Brno, Czech Republic.\\
\end{tabular}

\vspace*{3mm}{\it Teaching}\vspace*{1mm}

\begin{tabular}{p{40mm} p{92mm}}
Since 02/2000 & Teaching assistant at the Department of Analysis
and Technical Mathematics and the Department of Applied and
Numerical Mathematics, Vienna University of Technology.
\end{tabular}

\vspace*{3mm} {\it Awards}\vspace*{1mm}

\begin{tabular}{p{40mm} p{92mm}}
1999, 2000, 2002 & Scholarships for outstanding studies awarded by the Vienna University of Technology.\\
10/2002 - 09/2004 & DOC - Scholarship awarded by the Austrian
Academy of Sciences.
\end{tabular}

\vspace*{3mm} {\it Conferences and talks}\vspace*{1mm}

\begin{tabular}{p{40mm} p{92mm}}
11/2002 & Talk on ``Rosenberg's
classification of maximal clones on finite sets'' in the Vienna Algebra Seminar.\\
03/2003 & Participant of the AAA 65 conference,
Potsdam, Germany.\\
06/2003 & Talk on ``Clones containing all almost unary functions'' at Free University Berlin, Germany.\\
06/2003 & Talk on ``Clones on the natural numbers'' at the AAA 66
conference,
Klagenfurt, Austria.\\
09/2003 & Talk on ``Clones on infinite sets'' at the Summer School
on General Algebra and Ordered Sets,
Ko\v{s}icka Bela, Slovak Republic.\\
11/2003 & Two talks on ``Large clones on infinite sets'' at Masaryk University Brno, Czech Republic.\\
12/2003 & Talk on ``Clones above the unary functions'' in the
Vienna Algebra Seminar.\\
03/2004 & Talk on ``Maximal clones containing the permutations''
at the AAA 67 conference, Potsdam, Germany.\\
06/2004 & Talk on ``Monoidal intervals in the clone lattice'' at
the AAA 68 conference, Dresden, Germany.\\
07/2004 & Talk on ``Set theory in infinite clone theory'' at the
Logic Colloquium, Torino, Italy.
\end{tabular}
\end{document}